\documentstyle[eqsection]{article}
\newfont{\ind}{cmbx6}

\newcommand{\ovf}{\overline{\varphi}}
\newenvironment{pf}{\noindent{\sc Proof}.\enspace}{\rule{2mm}{2mm}\medskip}
\newenvironment{pfn}{\noindent{\sc Proof} \enspace}{\rule{2mm}{2mm}\medskip}
\newtheorem{theorem}{Theorem}[section] \newtheorem{lemma}{Lemma}[section]

\newtheorem{remark}{Remark}[section]
\newtheorem{remarks}{Remark}[section]
\newtheorem{definition}{Definition}[section]
\newcommand{\be}{\begin{equation}}
\newcommand{\ee}{\end{equation}}

\newcommand{\teta}{\theta}

\newcommand{\om}{\omega}

\newcommand{\ov}{\overline}

\newcommand{\wtilde}{\widetilde}

\renewcommand{\a }{\alpha }
\renewcommand{\b }{\beta }

\renewcommand{\d }{\delta }

\newcommand{\e }{\varepsilon }
\newcommand{\g }{\gamma}

\renewcommand{\l }{\lambda }
\renewcommand{\L }{\Lambda }
\newcommand{\m }{\mu }
\newcommand{\vphi}{\varphi }

\renewcommand{\o }{\omega }
\renewcommand{\O }{\Omega }

\newcommand{\dps}{\displaystyle}

\newcommand{\f }{\varphi }
\newcommand{\th }{\teta}

\newcommand{\fnl }{f_{nl} }
\newcommand{\bnl }{\beta_{nl} }

\newcommand{\dpr }{\partial}

\newcommand{\bc }{\ov{c} }

\newcommand{\oI }{\ov{I} }
\newcommand{\of }{\ov{\f} }
\newcommand{\oP }{\ov{P} }
\newcommand{\oQ }{\ov{Q} }
\newcommand{\oB }{\ov{B} }
\newcommand{\oD }{\ov{D} }
\newcommand{\og }{\ov{g} }
\newcommand{\oH }{\ov{H} }
\newcommand{\Et }{\wtilde{E} }
\newcommand{\dP }{\dot{P} }
\newcommand{\ddP }{\ddot{P} }
\newcommand{\dddot}[1]{\mathop{#1}\limits^{\dots}}
\newcommand{\wl}{\widetilde{\l}}

\begin{document}

\title{{\bf Drift in phase space: a new variational mechanism 
with optimal diffusion time}}

\author{Massimiliano Berti, Luca Biasco and Philippe Bolle}
\date{}
\maketitle

{\bf Abstract:}
We consider non-isochronous, nearly integrable,   
a-priori unstable Hamiltonian systems
with a (trigonometric polynomial) $O(\mu)$-perturbation 
which does not preserve the unperturbed tori. 
We prove the existence of Arnold diffusion with diffusion time 
$ T_d = O((1/ \mu) \log (1/ \mu ))$ by a variational method  
which does not require the existence
of ``transition chains of tori'' provided by KAM theory.
We also prove that our estimate of the diffusion time $T_d $  
is optimal as a consequence of a general stability result
derived from classical perturbation theory.
\footnote{Supported by M.U.R.S.T. Variational Methods and Nonlinear
Differential Equations.}
\\[2mm]
Keywords: Arnold diffusion, variational methods,
shadowing theorem, perturbation theory, nonlinear functional analysis
\\[1mm]
AMS subject classification: 37J40, 37J45.

\section{Introduction and main results}

Topological instability of action variables in multidimensional
nearly integrable Hamiltonian systems is known as Arnold Diffusion.
For  autonomous Hamiltonian systems with two degrees of freedom
KAM theory generically implies topological stability 
of the action variables, i.e. under the flow of the perturbed 
system the action variables  stay close to their initial values for all times. 
On the contrary, for systems with more than two degrees of freedom,
outside a large set of initial conditions provided by KAM theory,
the action variables may undergo a drift of order one in
a very long, but finite time called the ``diffusion time''. 
Arnold first showed up this instability phenomenon for a peculiar Hamiltonian
in the  famous paper \cite{Arn}.

As suggested by normal form theory near simple resonances,
the Hamiltonian models which are usually studied have the form 
$H (I, \vphi, p, q ) =  (I_1^2 /2 ) + \om \cdot I_2 + 
(p^2 / 2 ) + \e ( \cos q - 1 ) + \e \mu f ( I, \vphi,  p, q ) $
where $ \e $ and $ \mu $ are small parameters, $n := n_1 + n_2 $,
$ ( I_1, I_2, p ) \in  {\bf R}^{n}  \times 
{\bf R}$ are the action variables and $ ( \vphi, q ) = $ 
$ ( \vphi_1, \vphi_2, q) \in  {\bf T}^{n} \times 
{\bf T} $ are the angle variables. In Arnold's model $ I_1, I_2 \in {\bf R}$,
$ \om = 1 $, $ f ( I, \vphi, p, q ) = 
(\cos q - 1 )(\sin \vphi_1 + \cos \vphi_2 )$ and 
diffusion is proved for $ \mu $ 
exponentially small w.r.t.  $ \sqrt{\e} $.
Physically Hamiltonian $ H $ describes a system of $ n_1 $ ``rotators'' 
and $ n_2 $ harmonic oscillators weakly coupled with a pendulum
through a  perturbation term.

The mechanism proposed in \cite{Arn} to prove the existence 
of Arnold diffusion and thereafter become classical, is the following one.
For $ \mu = 0 $, the Hamiltonian system associated to $ H $ admits
a continuous family of $n$-dimensional 
partially hyperbolic invariant tori
$ {\cal T}_I$ = $\{ \vphi  \in  {\bf T}^n,
(I_1, I_2) = I, \ q = p = 0 \}$
possessing stable and unstable manifolds
$ W^s_0 ( {\cal T}_I) = W^u_0 ( {\cal T}_I) =$
$ \{ \vphi \in  {\bf T}^n,
(I_1, I_2) = I, \ (p^2/2) + \e (\cos q - 1) = 0 \}$.
The method used in \cite{Arn} to produce unstable orbits 
relies on the construction, for $\mu \neq 0, $ of {\it ``transition chains''} 
of perturbed partially hyperbolic  
tori $ {\cal T}_I^\mu $ close to $ {\cal T}_I $
connected one to another by heteroclinic orbits. Therefore in
general the first step 
is  to prove the persistence of such hyperbolic tori 
$ {\cal T}_I^\mu $ for $ \mu \neq 0 $ small enough, and to show that 
the perturbed stable and unstable manifolds $ W^s_\mu ({\cal T}_I^\mu )$
and $ W^u_\mu ({\cal T}_I^\mu) $ split and intersect transversally
(``splitting problem''). 
The second step is to find a transition chain of perturbed tori:
this is a diffucult task since, for general non-isochronous systems,  
the surviving perturbed tori $ {\cal T}^\mu_I $ are separated 
by the gaps appearing in KAM constructions.
Two perturbed invariant tori $ {\cal T}_I^\mu $ and $ {\cal T}_{I'}^\mu $
could be too distant one from the other,
 forbidding the existence of a heteroclinic intersection
between $ W^u_\mu ({\cal T}_I^\mu )$ and $ W^s_\mu ({\cal T}_{I'}^\mu )$:
this is the so called {\it ``gap problem''}. 
In \cite{Arn} this difficulty is bypassed   by the peculiar 
choice of the perturbation $ f ( I, \vphi, p, q ) = 
(\cos q - 1 ) f (\vphi) $, 
whose gradient vanishes on the unperturbed tori $ {\cal T}_I $, leaving 
them {\it all} invariant also for $ \mu \neq 0 $.  
The final step is to prove, by a ``shadowing argument'', the existence of 
a true diffusion orbit, close to a given transition chain 
of tori, for which the action variables $I$ 
undergo a drift of $O(1)$
in a certain time $ T_d $ called the {\it diffusion time}.

The first paper proving  Arnold diffusion in presence of
perturbations not preserving the unperturbed tori has been \cite{CG}. 
Extending  Arnold's analysis, it is proved in \cite{CG}  
that, if the perturbation is 
a trigonometric polynomial in the angles $ \vphi $, then,
in some regions of the phase space, the ``density'' of 
perturbed invariant tori is  high enough to allow the 
construction of a transition chain. 

Regarding the shadowing problem, geometrical method, see e.g.
\cite{CG}, \cite{M}, \cite{C1}, \cite{CrG},  
and variational  ones, see e.g. \cite{BCV}, 
have been applied, in the last years, 
in order to prove the existence of diffusion orbits shadowing 
a given transition chain of tori and to estimate the diffusion time.
We also quote the important papers \cite{Bs}-\cite{Bs1} which,
even if 
dealing with  Arnold's model perturbation only, 
have obtained, by variational methods, very good 
 diffusion time estimates 
and have introduced new ideas for studying the shadowing problem.
For isochronous systems new variational results concerning the
shadowing and the splitting problem have been obtained in
\cite{BBN}, \cite{BB2} and \cite{BB3}. 
\\[2mm]
\indent
In this paper we provide an {\it alternative mechanism} 
to produce diffusion orbits. This method is
not based on the existence of a transition chain of tori:
we avoid the KAM construction of the perturbed hyperbolic tori, 
proving directly the existence of 
a drifting orbit as a local minimum of an action functional.
At the same time our variational approach 
achieves the optimal diffusion time.  
We also prove that our  diffusion time estimate 
is the optimal one as a consequence of a general stability result,
proved via classical perturbation theory. 
As in \cite{CG} we deal with a perturbation  
which is a trigonometric polynomial in the angles and 
our diffusion orbits will not connect any two arbitrary
frequencies of the action space, even if we manage to connect
more frequencies than in \cite{CG}, proving the drift also
in some regions of the phase space  
where transition chains might not exist.
Clearly if the perturbation is chosen as in  Arnold's example 
we can drift in all the phase space with no restriction.
The results proved here have been announced in \cite{BBB}.

In this paper we will assume, 
as in Arnold's paper, the parameter $ \mu $ to be 
small enough in order to validate the so called 
Poincar\'e-Melnikov approximation,
when the first order expansion term in $ \mu $ for the splitting, 
the so called Poincar\'e-Melnikov function, is the dominant one.
For this reason, through this paper we will fix
the ``Lyapunov exponent'' of the pendulum $ \e := 1 $,
considering the so called ``a-priori unstable'' case.
Actually our variational shadowing technique is not 
restricted to the a-priori unstable case, but would
allow, in the same spirit of \cite{BBN}, \cite{BB2} and \cite{BB3}, 
once a ``splitting condition'' is someway proved,
to get diffusion orbits with the best diffusion time
(in terms of some measure of the splitting).  
\\[2mm]
\indent
We will consider nearly integrable
{\it non-isochronous} Hamiltonian systems defined by
\be\label{eq:Hamg}
{\cal H}_\mu =  \frac{I^2}{2} + \frac{p^2}{2} + (\cos q - 1) +
\mu f(I,  \vphi, p, q , t),
\ee
where $( \vphi, q , t) \in {\bf T}^d \times {\bf T}^1 \times {\bf T}^1 $
are the angle variables, $ ( I, p ) \in {\bf R}^d \times {\bf R}^1$
are the action variables and $ \mu \geq 0 $ is a small real  
parameter.
The Hamiltonian system associated with ${\cal H}_{\mu}$ writes
$$
\dot{\vphi}= I + \mu \partial_I f,  \quad
\dot{I}= -\mu \partial_{\vphi} f, \quad
\dot{q}= p + \mu \partial_p f,  \quad  
\dot{p} = \sin q - \mu \partial_q f.  \eqno({\cal S}_\mu)
$$
The perturbation $f$ is assumed to be a real trigonometric polynomial
of order $ N $ in $ \vphi $ and $ t $, namely\footnote{
$ \ov{f}_{n,l} (I, p, q) = f_{-n,-l} (I, p, q) $ for all 
$ (n,l) \in {\bf Z}^d \times {\bf Z} $ with 
$ |(n,l)| \leq N $ where
$\ov{z}$ denotes the complex conjugate of $z \in {\bf C}$.}
\be\label{trigpert}
f(I, \vphi, p, q, t) = \sum_{|(n,l)| \leq N} f_{n,l} (I, p, q) 
e^{ {\rm i}(n \cdot \vphi  + l t)}.
\ee
The unperturbed Hamiltonian system $({\cal S}_0) $ is completely integrable and
in particular the energy $ I_i^2/ 2 $ of each rotator is a
constant of the  motion.
The problem of {\it Arnold diffusion} in this context is whether,
for $ \mu \neq 0 $, there exist motions whose net effect is to transfer
$O(1)$-energy among the rotators. A natural complementary question regards 
the time of stability (or instability) 
for the perturbed system: what is the minimal time to 
produce an $O(1)$-exchange of energy, if any takes place, 
among the rotators?
\\[1mm]
\indent 
For simplicity, even if it is not really necessary, we assume $f$ to be 
a purely spatial perturbation, namely 
$f(\vphi, q, t) = \sum_{0\leq |(n,l)| \leq N} f_{n,l} (q) 
\exp ({\rm i}(n \cdot \vphi  + l t)) $. 
The functions $ f_{n,l} $ are assumed to be smooth.
\\[2mm]
\indent
Let  us define the ``resonant web'' $ {\cal D}_N $, 
formed by the frequencies
$ \om $ ``resonant with the perturbation''  
\be\label{eq:resonantweb}
{\cal D}_N := \Big\{ \om \in {\bf R}^d \ \Big| \ \exists (n, l)
\in {\bf Z}^{d+1} \ {\rm s. t.} \
0 < | ( n, l ) | \leq N \ {\rm and } \ \om \cdot n + l = 0 \Big\} = 
\cup_{0< |(n,l)| \leq N} E_{n,l}
\ee
where $ E_{n,l} := \{ \om \in {\bf R}^d \ | \ \om \cdot n + l = 0 \} $. 
Let us also consider the Poincar\'e-Melnikov primitive 
$$
\Gamma (\om, \teta_0  , \vphi_0 ) := - \int_{\bf R}
\Big[ f ( \omega t + \vphi_0, q_0 (t), t + \teta_0) -  
f ( \omega t + \vphi_0, 0 , t + \teta_0 ) \Big] \ dt, 
$$ 
where $ q_0 (t) = 4 \ {\rm arctan} (\exp t ) $ is
the  separatrix of the unperturbed  pendulum equation
$\ \ddot{q}=\sin q \ $ satisfying $ q_0 (0) = \pi $.

The next Theorem states  that, for any connected component
${\cal C} \subset  {\cal D}_N^c $, $ \om_I, \om_F \in {\cal C} $,  
there exists a solution of $ ({\cal S}_\mu) $
connecting a $O(\mu)$-neighborhood of $ \om_I $ in the action space to a 
$O(\mu)$-neighborhood of $ \om_F $, in the time-interval 
$ T_d = O((1 / \mu) | \log \mu |) $.

\begin{theorem}\label{thm:main}
Let ${\cal C}$ be a connected component of $ {\cal D}_N^c $, 
$\om_I, \om_F \in {\cal C} $ and let 
$ \gamma : [0,L] \to {\cal C} $ be a smooth embedding such that
$ \g (0) = \om_I $ and $\gamma (L) = \om_F $.
Assume that, for all $ \om := \gamma (s) $ $(s \in [0,L])$,  
$ \Gamma(\om, \cdot, \cdot) $ possesses a 
non-degenerate local minimum $ ( \teta_0^\om , \vphi_0^\om ) $. 
Then $\forall \eta > 0 $ there exists $ \mu_0 = \mu_0 (\g, \eta) > 0 $
and $C= C(\gamma) >0 $ such that $ \forall 0 < \mu \leq \mu_0 $ 
there exists a solution 
$ (I_\mu (t) , \vphi_\mu (t) ,   p_\mu (t) ,  q_\mu (t))$ of 
$ ({\cal S}_\mu) $ and two instants 
$ {\tau}_1 < {\tau}_2 $ such that $ I_\mu ( \tau_1 ) = \om_I + O(\mu) $,  
$ I_\mu ( \tau_2 ) = \om_F + O(\mu) $
and 
\be\label{totaltime}
|{\tau}_2 - {\tau}_1| \leq  
\frac{C}{ \mu } | \log \mu |.
\ee
Moreover ${\rm dist} (I_\mu (t), \gamma ([0,L]) ) < \eta $ 
for all $\tau_1 \leq t \leq \tau_2 $. 

In addition, the above result still holds for any 
perturbation $ \mu ( f + \mu {\wtilde f})$ with any 
smooth $\wtilde{f}(\vphi,q,t)$.
\end{theorem} 

We can also build diffusion orbits approaching
the boundaries of $ {\cal D}_N $ at distances as  
small as a certain power of $ \mu $: see for a precise statement 
Theorem  \ref{thm:tempobordo}.

Theorem \ref{thm:main} improves the corresponding result in \cite{CG}
which enables to connect two frequencies $ \om_I $ and $ \om_F $ 
belonging to the same connected component
${\cal C} \subset {\cal D}_{N_1}^c $ for $ N_1 = 14 d N $ and 
with dist$\{ \{ \om_I, \om_F \}, {\cal D}_{N_1} \} = O( 1 )$. 
Such restrictions of \cite{CG} in connecting 
the action space through diffusion orbits 
arise because transition chains could not exist in all 
$ {\cal C} \subset {\cal D}_N^c $ (see remark \ref{transchain}).
Unlikely our method enables to show up Arnold diffusion 
between any two frequencies 
$ \om_I, \om_F \in {\cal C} \subset {\cal D}_N^c $ and
along any path, since it does not require the existence of chains of 
true hyperbolic tori of (${\cal S}_\mu$). 
   
Theorem \ref{thm:main} also
improves the known estimates on the diffusion time.
The first estimate obtained
by geometrical method in \cite{CG}, is $ T_d = O( \exp{(1/\mu^2)} )$.
In \cite{M}-\cite{C1}-\cite{CrG}, 
still by geometrical methods, and in \cite{BCV}, 
by means of Mather's theory, the diffusion time 
has been proved to be just polynomially long in the splitting $ \mu $ 
(the splitting angles between the perturbed stable and unstable
manifolds ${\cal W}^{s,u}_{\m} ({\cal T}_{\om}^\m)$
at a homoclinic point are, by classical Poincar\'e-Melnikov theory, 
$ O(\mu )$).
We note that the variational method proposed by Bessi in \cite{Bs} had 
already given,  in the case of perturbations preserving 
all the  unperturbed tori, the  diffusion time estimate $ T_d = O(1/ \mu^2) $. 
For isochronous systems the estimate on the 
diffusion time  $ T_d = O((1/ \mu) |\ln \mu |) $ has  already been obtained in 
\cite{BBN}-\cite{BB2}. Very recently, in \cite{CrG}, the diffusion
time (in the non isochronous case) 
has  been estimated as $ T_d = O((1/\mu) |\log \mu |)$ by a method 
which uses ``hyperbolic periodic orbits''; 
however the result of \cite{CrG} is of local nature: 
the previous estimate holds only for 
diffusion orbits shadowing a transition chain close to some torus run with
diophantine flow. 

Our next statement (a stability result) concludes this quest for the minimal 
diffusion time $ T_d $: it shows the 
optimality of our estimate $ T_d = O( (1 / \mu) |\log \mu |)$.

\begin{theorem}\label{thm:Nec}
Let $ f ( I, \vphi, p, q, t) $ be as in (\ref{trigpert}), where the 
$f_{n,l}$ ($|(n,l)| \leq N$) are analytic functions. Then 
$\forall  \kappa, \ov{r}, \wtilde{r} > 0 $ there exist
$ \mu_1, \kappa_0 > 0 $ such that $ \forall \, 0 < \mu \leq \mu_1 $,
for any solution $ ( I(t), \vphi(t) , p(t) , q(t))$ of $({\cal S}_\mu) $
with  $| I(0) | \leq \ov{r} $ and   
$ | p(0) | \leq \wtilde{r} $, there results
\be\label{eq:stab}
|I(t)-I(0)| \leq  \kappa\ \qquad \ \forall\  t \ \ {\rm such} \ \
{\rm that} \ \  
|t| \leq \frac{\kappa_0}{\mu} \ln \frac{1}{\mu}.
\ee
\end{theorem}

Actually the proof of Theorem \ref{thm:Nec} contains much more
information: in particular the stability time (\ref{eq:stab}) is sharp
only for orbits lying close to the separatrices. On the other hand
the orbits lying far away from the separatrices are much more
stable, namely exponentially stable in time according to Nekhoroshev 
type time estimates, see (\ref{expstabf}) and (\ref{stabexpi}).
Indeed the diffusion orbit of Theorem \ref{thm:main} is found 
close to some pseudo-diffusion orbit whose $(q,p)$ variables move 
along the separatrices of the pendulum. 
\\[1mm]
\indent
As a byproduct of the techniques developed in this paper 
we have the following result (proved in section 6)
concerning  ``Arnold's example''  \cite{Arn}
where $ {\cal T}_{\om} := \{I= \om, \vphi \in {\bf T}^d, 
p = q = 0 \} $ are, for all $ \omega \in {\bf R}^d $, 
even for $ \mu \neq 0 $, invariant tori of $({\cal S}_\mu)$.

\begin{theorem}\label{thm:Arn}
Let $ f(\vphi, q, t) := (1 - \cos q) {\wtilde f} (\vphi, t) $. Assume that 
for some smooth embedding $ \gamma : [0,L] \to {\bf R}^d $, 
with $\g (0) = \om_I $ and $\gamma (L) = \om_F $, 
$\forall \om := \gamma (s) $ $(s \in [0,L])$,  
$ \Gamma(\om, \cdot, \cdot) $ possesses a 
non-degenerate local minimum $ (  \teta_0^\om , \vphi_0^\om ) $. 
Then $ \forall \eta > 0 $ there exists $ \mu_0 = \mu_0 ( \g , \eta ) > 0 $,
and $ C= C(\gamma) >0 $ such that $\forall 0 < \mu \leq \mu_0 $ 
there exists a heteroclinic orbit ($\eta$-close to $ \gamma $) 
connecting the invariant tori $ {\cal T}_{\om_I}$ and $ {\cal T}_{\om_F} $.
Moreover the diffusion time $ T_d $ needed
to go from a $\mu$-neighbourhood of
${\cal T}_{\om_I}$ to a $\mu$-neighbourhood of
${\cal T}_{\om_F}$ is bounded by
$  (C / \mu)  | \log \mu |$ for some constant $C$.
\end{theorem} 

The method of proof of Theorem \ref{thm:main} (and Theorem \ref{thm:Arn}) 
relies on a finite dimensional reduction 
of Lyapunov-Schmidt type, variational in nature, 
introduced in \cite{AB} and later extended in \cite{BBN},\cite{BB2} 
and \cite{BB3} to the problem of Arnold diffusion.
The diffusion orbit of Theorem \ref{thm:main} is found 
as a local minimum of the action functional 
close to some pseudo-diffusion orbit whose $(p,q)$ variables move 
along the separatrices of the pendulum. The pseudo-diffusion orbits, 
constructed by the Implicit Function Theorem, are
true solutions of  $ ({\cal S}_\mu) $ 
except possibly at some instants $ \teta_i $, for $i=1, \ldots, k $,
when they are glued 
continuously at the section $\{ q = \pi , \ {\rm mod} \ 2 \pi {\bf Z} \}$ 
but the speeds $({\dot \vphi}_\mu (\teta_i), {\dot q}_\mu (\teta_i) ) =  
(I_\mu (\teta_i), p_\mu (\teta_i) ) $ may  have a jump.
The time interval $ T_s = \teta_{i+1} - \teta_i $ is heuristically the
time required to perform a single transition during which the rotators 
can exchange $ O(\mu)$-energy, i.e. the action variables vary of $O(\mu)$.
During each transition we can exchange only $O(\mu )$-energy 
because the Melnikov contribution in the perturbed functional
is $O(\mu)$. Hence 
in order to exchange $O(1)$ energy the number of transitions required 
will be $k = O(1/ \mu)$.
\\[1mm]
\indent
We underline that the question of finding the optimal time and 
the mechanism for which we can avoid the construction of  transition chains
of tori are deeply connected. Indeed the main 
reason for which our drifting technique  avoids the construction 
of KAM tori is the following one: if
the time to perform a simple transition $ T_s $
is, say, just $ T_s = O( |\ln \mu | ) $  then, on such ``short'' 
time intervals, 
it is valid to approximate the pseudo diffusion orbits with  unperturbed 
solutions 
living on the stable and unstable manifolds of the unperturbed tori 
$ W^s ({\cal T}_\om ) = W^u({\cal T}_\om ) = 
\{ I = \om, \vphi \in  {\bf T}^d, p^2 / 2 + ( \cos q - 1) = 0 \}$,
when computing the value of the action functional. 
In this way we do not need 
to construct the true hyperbolic tori $ {\cal T}_\om^\mu $
(actually for our approximation we only need the time for a 
single transition to be $ T_s << 1/ \mu $).
\\[1mm]
\indent
The fact that it is possible to perform a single transition in a 
very short time interval like $ T_s = O(|\ln \mu |)$ is not obvius at all. 
In \cite{Bs} the time to perform a single transition,
in the example of Arnold, is $ O(1 / \mu )$. 
This transition time arises in order to ensure that the 
variations of the kinetic part of the action functional 
associated with  the rotators are small compared with 
the (positive definite) second derivative  
of the Poincar\'e-Melnikov primitive at its minimum point.
Unfortunately this time is too long  to use a simple 
approximation of the functional. 
The key observation that enables us to perform 
a single transition in a very short time interval  
concerns the behaviour of  the ``gradient flow'' 
of the unperturbed action functional of the rotators. This implies
a sort of a-priori estimate  satisfied by the minimal 
diffusion orbits, see remark \ref{rem:aprio}. 
We think that estimate (\ref{apiroriest}) is 
interesting in itself. In this way we can show that the variations of 
the action of the rotators are small enough, 
even on time intervals $T_s << 1/ \mu $, and do not ``destroy'' the minimum
of the Poincar\'e-Melnikov primitive. 

When trying to build a pseudo-diffusion orbit which 
performs single transitions in very short time intervals 
we encounter another difficulty linked with the ergodization time. 
The time to perform a single transition $T_s$ must be  long enough 
 to settle, at each instant $ \teta_i $, the projection
$(\theta_i, \f_i)$ of the 
pseudo-orbit on the torus ${\bf T}^{d+1}$ sufficiently close
to the minimum of the Poincar\'e-Melnikov function, i.e. the homoclinic point
(in our method 
it is sufficient to arrive just $O(1)$-close, independently of $ \mu$,
to the homoclinic point). This necessary request creates some
difficulty since our pseudo-diffusion orbit 
may arrive $O(\mu)$-close in the action space 
to resonant hyperplanes of  frequencies
whose linear flow does not provide a dense enough net of the torus. 
The way in which this problem is overcome is discussed in section
\ref{sec:pseudodiff}:
we observe a phenomenon of ``stabilization close to resonances''
which forces the time for some single transitions 
to increase. 
Anyway the total time required to cross these
(finite number of) resonances is still 
$ T_d = O( (1/ \mu) \log (1 / \mu)) $,
see (\ref{lengthest}) and the proof of Theorem \ref{thm:main}. 
This  discussion  enables us to prove optimal 
fast-Arnold diffusion in large regions of the phase 
space and allows to improve 
the local diffusion results of \cite{CrG}.

We need therefore some results on the ergodization time of the torus 
for linear flows possibly resonant but 
only at a ``sufficiently high order''.
We present these results in section \ref{sec:erg}. We point out that
the main result of this section, Theorem \ref{lem:ergo2}, implies
as corollaries Theorems B and D of \cite{BGW}, see remark \ref{rem:BGW}.
It is of independent interest and could possibly improve the other
results of \cite{BGW}.
\\[1mm]\indent
This work is  a further step of a reaserch line,
started in \cite{BBN}-\cite{BB2} and \cite{BB3}, 
for finding new mechanisms to prove Arnold diffusion.
We expect that the variational method developed in this paper 
could be suitably refined in order to prove the existence of 
drifting orbits in the whole action space and then to prove 
such results  for generic analytic perturbations too.
Another possible application of these methods could regard 
infinite dimensional Hamiltonian systems where the existence of 
``transition chains of infinite dimensional hyperbolic tori'' 
is quite far for being proved.
\\[1mm]
\indent
The paper is organized as follows:
in section 2  we perform the finite dimensional reduction and we 
define the variational setting.
In section 3 we provide a suitable development of the reduced 
action functional.
In section 4 we prove the new results on the ergodization time. 
In section 5 we define the unperturbed pseudo-orbit.
In section 6 we prove the existence of the 
diffusion orbit. In section 7 we prove the stability result, that
is to say the 
optimality of our diffusion time.
\\[2mm]
\indent
\textbf{Notations:}
Through this paper the notation $ a(z_1 , \ldots , z_k)
= O ( b( \mu ) )$ will mean that,
for a suitable positive constant $  C(\gamma, f) > 0 $, 
$ | a (z_1 , \ldots , z_p ) | \leq$
$ C(\gamma, f) | b (\mu) |$. 

\section{The variational setting and the 
finite dimensional reduction} \label{sec:3}

When the perturbation
$f(\vphi, q, t) = \sum_{|(n,l)| \leq N} f_{n,l} (q) 
exp ({\rm i}(n \cdot \vphi  + l t)) $ is purely 
spatial, \footnote{We will develop
all the computations for $f$. 
All the next arguments remain unchanged if
the perturbation is $ f + \mu {\wtilde f}$, see the proof
of Theorem \ref{thm:main}.}
system $({\cal S}_\mu)$ reduces to the second order system
\be\label{lagraeqmot}
\ddot{\vphi} = - \mu  \ \partial_{\vphi} f ( \vphi, q,t), \qquad 
- \ddot{q} + \sin{q} = \mu \ \partial_q f ( \vphi, q, t)
\ee
with associated Lagrangian
\be
{\cal L}_\mu (\vphi, \dot{\vphi}, q, \dot{q}, t) 
= \frac{{\dot \vphi}^2}{2} + 
\frac{{\dot q}^2}{2} + (1- \cos q) - \mu f( \vphi, q , t).
\ee
\noindent
Using the Contraction Mapping Theorem 
we will prove in lemma \ref{lem:connecting}
that, near the unperturbed solutions $(\om (t - \teta) + \vphi_0, 
q_0 (t- \teta) )$ living on the stable and unstable manifolds of 
the unperturbed tori $ {\cal T}_\om $,  
there exist, for $ \mu $ small enough, 
solutions of the perturbed system (\ref{lagraeqmot}) 
which connect the sections $\{ \f = \f^+, q = - \pi, t = \theta^+ \}$ and 
$ \{ \f = \f^-, q = \pi, t = \theta^- \}$
(under some assumptions). The diffusion orbit will be
a chain of such connecting orbits. 
\\[1mm]
\indent
We first introduce a few definitions and notations. 
For  $ \lambda := (\teta^+, \teta^-, \vphi^+, \vphi^-) \in 
{\bf R}^2 \times {\bf R}^{2d} $ with $\theta^+< $ $\theta^-$  
we define $T_{\lambda}:=\theta^--\theta^+$  and the ``mean frequency'' 
$ \om_{\lambda} \in {\bf R}^d $ as 
$
\om_{\lambda} := \dps \frac{\vphi^- - \vphi^+ }{ \teta^- - \teta^+ }.
$
The ``small denominator'' of  a  frequency $ \om \in {\bf R}^d $
is defined by  
\be\label{beta}
\b(\om):=\b_N(\om):=\min_{0<|(n,l)|\leq N}|n \cdot \om +l|.
\ee
$\beta(\om)$  measures how close the frequency $\om $ lies to the
resonant web  ${\cal D}_N$ defined in (\ref{eq:resonantweb}).
We use the abbreviation 
$\beta_{\lambda}$ for $\beta(\om_{\lambda})$.
We shall always assume through this paper that  
$\om $ stays in a fixed bounded set
containing the curve $ \gamma $.

For $T$ large enough, there exists a unique $T$-periodic solution $Q_T$
of the pendulum equation, of small positive energy with 
$ Q_T(0)= - \pi $, $Q_T(T)=\pi$. Moreover $Q_T$ satisfies
$ \forall t \in$ $[0,T/2) \cup (T/2,T]$,
$$
|\partial_T Q_T (t)| \leq K_1 e^{-K_2(T-t)} \quad ,
\quad  |\partial_T (Q_T(T-\cdot)) (t)| \leq K_1 e^{-K_2 (T-t)}
$$ and 
\begin{equation} \label{estQT}
|Q_T(t)- q_{\infty}(t)|+
|\dot{Q}_T(t)-\dot{q}_\infty (t)| \leq K_1 e^{-K_2 T} \  ,  \ 
|\dot{Q}_T (t)| \leq K_1 \max \{ e^{-K_2 t} ,  e^{-K_2 (T-t)} \},
\end{equation}
for some positive constants $K_1$ and $K_2$, where $q_\infty$
is defined by
$$
q_\infty (t)= q_0 (t) - 2\pi  \  {\rm if} \ t \in [0,T/2), \qquad 
 q_\infty (t) = q_0 ( t - T )  \  {\rm if} \ t\in (T/2,T].
$$ 

\begin{lemma} \label{lem:connecting} 
There exists  $\m_2 > 0$ and constants $C_0,C_1,\bc,c_1>0$ 
such that $\forall 0 < \m \leq \m_2 $, 
$\forall \lambda = ( \th^+, \th^-, \f^+, \f^-)$ such that 
$ C_0 \beta_{\l}^2 > \mu  $  and 
$C_1|\ln \m|\leq T_{\l}\leq C_0\b_{\l}/{\m}$  there exists a unique
solution $(\f_\mu (t), q_\mu (t) ) := 
(\f_{\mu,\l} (t), q_{\mu, \l}(t))$ of
(\ref{lagraeqmot}), defined for $t\in (\th^+-1,\th^-+1)$,
satisfying $\f_\mu (\theta^\pm) = \f^\pm $, 
$ q_\mu (\theta^{\pm})= \mp \pi $ and
\begin{equation}\label{stime}
\begin{array}{rl}
(i) & |\f_\mu (t) - \ov{\f}(t)| \leq  \bc\m(1+c_1\m T_\l^2)/{\b_\l^2}, \quad
\quad |\dot{\f}_\mu (t)-\o| \leq  \bc\m/{\b_\l}, \\
(ii) & |q_\mu (t)-Q_{T_\l}(t-\theta^+)| \leq  \bc\m, \quad \quad
|\dot{q}_\mu (t)-\dot{Q}_{T_\l}(t-\theta^+)| \leq \bc\m ,
\end{array}
\end{equation}
where $\ov{\f}(t):=\o_\l (t-\th^+)+\f^+$. 
Moreover $\f_{\m,\l}(t)$, ${\dot \f}_{\m,\l}(t)$,
$q_{\m,\l}(t)$ and ${\dot q}_{\m,\l}(t)$ are ${\cal C}^1$
functions of $(t,\lambda)$. 
\end{lemma}

The proof of lemma \ref{lem:connecting} is given in the Appendix.

\begin{remark}\label{remtranl}
Roughly, the meaning of the above estimates is the following.

1) We have imposed $ C_1 | \ln \mu | < T_\lambda := \teta^- - \teta^+ $ so
that by $(\ref{estQT})$,
on such intervals of time, the periodic solution $ Q_{T_\lambda} $
is $O(\mu)$ close to ``separatrices''  $ q_\infty $ 
of the unperturbed pendulum.

2) Estimate $(ii)$ implies that for 
$ t \approx ( \teta^+ + \teta^- ) \slash 2 $ the perturbed 
solution $ q_\mu $ may have $O(\mu)$ oscillations around the unstable
equilibrium of the pendulum $ q = 0 $, mod $ 2 \pi $, which is exactly
what one expects perturbing with a general $ f $. On the contrary 
for the class of perturbations considered in \cite{Arn} as 
$ f(\vphi, q, t ) = (1 - \cos q) f(\vphi, t) $
preserving {\it all} the invariant tori, estimate $(ii)$ 
can be improved, getting 
$ \max \{ | q_\mu (t) -$ $ Q_{T_{\l}}(t-\theta^+)|, $
$| {\dot q}_\mu (t) - {\dot Q}_{T_{\l}}(t-\theta^+)| \} =
O ( \mu \max \{ {\rm exp}(- C|t - \teta^+ |), {\rm exp}(- C|t - \teta^- |)\})$.

3) For $ \b_\l \approx \sqrt{\mu} $ estimate $(i)$ 
becomes meaningless: for a mean frequency $\om_\l$ such that  
$ n \cdot \om_\l  + l \approx \sqrt{\mu} $ 
for some $0 < |(n,l)| \leq N $ the perturbed 
transition orbits 
$\vphi_\mu $ are no more well-approximated by the straight lines 
$ \ov{\vphi} (t) := \vphi^+ + \om_\l (t - \teta^+ ) $. 
\end{remark}

\begin{remark}\label{transchain}
Let us define $ {\cal D}_N^\b := \{ \om \in {\bf R}^d \ | 
\ |\om \cdot n + l| > \b, \ \ \forall \ 0 < |(n,l )| \leq N \}$.
In \cite{CG} it is proved that hyperbolic invariant
tori $ {\cal T}_\om^\mu $ of system (${\cal S}_\mu$) exist for 
Diophantine frequencies $ \om \in  {\cal D}^{\b_1}_{N_1} $, for some 
$\b_1 = O(1)$ and some $ N_1 = O(d N) > N $, namely avoiding more 
``resonances with the trigonometric polynomial $f$'' than just $N$.  
The presence of such ``resonant hyperplanes $E_{n,l}$'' for 
$N < |(n,l)| < N_1 $ may be reflected  in estimate $(i)$
by the term $ \mu T_\l^2 $. However 
such term, for our purposes, can be ignored. 
From this point of view  lemma \ref{lem:connecting} could 
perhaps be interpreted 
as the first iterative step for
looking at invariant hyperbolic tori in the perturbed system bifurcating 
from the unperturbed one's .
\end{remark}

By lemma \ref{lem:connecting}, for $0<\mu \leq \mu_2$, we can define
on the set 
$$
\Lambda_{\m} := \Big\{ \l = (\th^+,\th^-,\f^+,\f^-) \ \Big| \
  \ C_0 \beta_{\l}^2 > \mu, \
 C_1 |\ln \m| \leq T_{\l} \leq \frac{C_0 \beta_\l}{\m } \Big\},
$$
the Lagrangian action functional $ G_\m : \Lambda_\m  \to {\bf R}$ as
\be \label{defiG}
G_\mu (\l ) =
G_\m(\theta^+,\theta^-,\f^+,\f^-):= \int_{\theta^+}^{\theta^-}
{\cal L}_{\mu}(\f_{\mu}(t),
\dot{\f}_{\m} (t) , q_{\mu}(t),\dot{q}_{\m}(t), t) \ dt.
\ee
We have
\begin{lemma} \label{lem:derG}
$G_\m$ is differentiable and (with the abbreviations $\f, q$ for 
$\f_\m , q_\m$)
$$
\nabla_{\f^+} G_\m(\lambda)= - \dot{\f} (\theta^+), \qquad
\partial_{\theta^+} G_\m (\lambda)= \frac{1}{2} |\dot{\f}
 (\theta^+)|^2 + \frac{1}{2} \dot{q}^2
 (\theta^+)+\cos q (\theta^+) -1 + \mu f(\f^+,\pi, \theta^+)
$$
$$
\nabla_{\f^-} G_\mu(\lambda)= \dot{\f} (\theta^-), \quad 
\partial_{\theta^-} G_\m(\lambda)=-\Big(\frac{1}{2} |\dot{\f}
 (\theta^-)|^2 + \frac{1}{2} \dot{q}^2
 (\theta^-)+\cos q (\theta^-) -1 + \mu f(\f^-, \pi, \theta^- )\Big).
$$
\end{lemma}
\begin{pf}
By lemma \ref{lem:connecting} the map 
$(\lambda,t) \mapsto ( \f_{\m,\l} (t) ,\dot{\f}_{\m,\l}(t), q_{\mu,\l} (t) ,
\dot{q}_{\m,\l}(t))$ is $ C^1 $ on the set 
$ \{ (\l,t)  \in $ $\Lambda_\m \times {\bf R} 
 \ | $ $ \  \theta^+ \leq t \leq \theta^- \}$.
Hence $G_\m$ is differentiable and
\begin{eqnarray*}
\partial_{\theta^+} G_\m(\lambda)& =& 
-{\cal L}_{\m} (\f^+, \dot{\f} (\theta^+),- \pi, \dot{q} (\theta^+), \teta^+) 
+ \int_{\theta^+}^{\theta^-} 
\dot{\f} (s) \cdot \partial_{\theta^+} {\dot \f} (s) +
\dot{q}(s) \partial_{\theta^+} {\dot q} (s) \ ds  \\
&+& \int_{\theta^+}^{\theta^-}
\sin q (s) \partial_{\theta^+} q (s)- \mu
\partial_{\f} f (\f (s),q (s),s) \cdot \partial_{\theta^+}
\f (s) - \mu \partial_{q} f (\f (s),q (s),s) 
\partial_{\theta^+} q (s) \ ds. 
\end{eqnarray*} 
Integrating by parts and using that $(q_{\m,\l}, \f_{\m,\l})$ 
satisfies (\ref{lagraeqmot}) in $(\theta^+, \theta^-)$, we obtain
$$
\partial_{\theta^+} G_\m(\lambda) = 
-{\cal L}_{\m} ( \f^+, \dot{\f} (\theta^+), -\pi , \dot{q} (\theta^+), \teta^+)
+ \Big[ \dot{q}(s) \partial_{\teta^+} q(s) 
+\dot{\f}(s) \cdot \partial_{\theta^+} \f (s)
\Big]_{\teta^+}^{\teta^-}.
$$
Now $q_{\m,\l}(\theta^+)=-\pi$ for all $\l$ hence 
$\dot{q}(\theta^+) + \partial_{\theta^+} q(\theta^+)=0$.
Similarly we get $\dot{\f}(\theta^+) + \partial_{\theta^+}
\f(\theta^+)=0$, $\partial_{\theta^+} q(\theta^-)=0 $,
$\partial_{\theta^+} \vphi (\theta^-)=0 $.
As a consequence 
$$
\partial_{\theta^+}G_\m (\lambda)=\frac{1}{2}
|\dot{\vphi}|^2(\theta^+) + \frac{1}{2}
\dot{q}^2(\theta^+) + (\cos q(\theta^+)-1)
+\m f(\f^+, \pi,\theta^+). 
$$ 
The other partial derivatives are computed in the same way.
\end{pf}

For $\beta>0 $ fixed,
denoting $ \l_i = (\theta_i,\theta_{i+1},\f_i,\f_{i+1})$, we define
on the set 
$$
\Lambda_{\m,k} :=
\Lambda_{\m,k}^{\b}:= \Big\{ \l = (\theta_1, \ldots, \theta_k, \f_1, 
\ldots , \f_k) \in {\bf R}^k \times {\bf R}^{kd} \ \Big| \ \forall 
\ 1\leq i \leq  k-1, \  \ \lambda_i \in \Lambda_\m \ , \ 
 \b_{\l_i} \geq \beta \Big\},
$$
the reduced action functional 
${\cal F}_\mu: \Lambda_{\m,k} \to {\bf R}$ as
\begin{eqnarray*}
{\cal F}_\mu ( \l ) & = & 
\om_I \vphi_1 - \frac{|\om_I|^2}{2} \teta_1 + 
\mu \Gamma^u (\om_I, \teta_1, \vphi_1) + \mu F (\om_I, \teta_1, \vphi_1) 
+ \sum_{i=1}^{k-1} G_\mu ( \l_i )\\
& - &  \om_F \vphi_k + \frac{|\om_F|^2}{2} \teta_k + 
\mu \Gamma^s (\om_F, \teta_k, \vphi_k) - \mu F (\om_F, \teta_k, \vphi_k)
\end{eqnarray*}
where 
\be\label{Gammau}
\Gamma^u (\om, \teta_0, \vphi_0 ) := -  \int_{- \infty}^0
\Big[ f ( \om t + \vphi_0, q_0 (t) , t + \teta_0 ) -  
f ( \om t + \vphi_0, 0 , t+ \teta_0 )) \Big] \ dt, 
\ee
\be\label{Gammas}
\Gamma^s (\om, \teta_0, \vphi_0 ) := - \int_0^{+ \infty}
\Big[ f ( \omega t + \vphi_0, q_0 (t), t + \teta_0 ) -  
f ( \om t + \vphi_0, 0 , t + \teta_0 )) \Big] \ dt, 
\ee
are called resp. the unstable and the 
stable Poincar\'e-Melnikov primitive, and
\be\label{risoterm}
F( \om, \teta_0, \vphi_0 ) := -f_{0,0} \teta_0 - \sum_{0 <|(n,l)| \leq N}
f_{n,l} \frac{e^{{\rm i}(n \cdot \vphi_0 
+ l \teta_0)}}{{\rm i}(n \cdot \om +l)},
\ee
 $ f_{n,l} := f_{n,l} (0) $ being the Fourier coefficients of 
$f(\vphi, 0, t)$.
\\[1mm]
\indent
Critical points of the
``reduced action functional'' $ {\cal F}_\mu $ 
give rise to  diffusion orbits 
whose action variables $I$  
go from a small neighbourhood of $ \om_I $ 
to a small neighbourhood of  $ \om_F $, 
as  stated in lemma \ref{lem:heter} below. The ``boundary terms''
$ \om_I \vphi_1 - \frac{|\om_I|^2}{2} \teta_1 + 
\mu \Gamma^u (\om_I, \teta_1, \vphi_1) + \mu F (\om_I, \teta_1, \vphi_1)$
and 
$ - \om_F \vphi_k + \frac{|\om_F|^2}{2} \teta_k + 
\mu \Gamma^s (\om_F, \teta_k, \vphi_k) - \mu F (\om_F, \teta_k, \vphi_k)$
have been added  also to enable us to find critical points of
$\cal{F}_{\mu}$ w.r.t. all the variables (including
$\theta_1, \vphi_1, \theta_k, \vphi_k$).

More precisely, for $\lambda = (\teta, \vphi) \in \Lambda_{\m,k}$
we define the pseudo diffusion solutions 
$(\f_{\m,\l}, q_{\m, \l })$ on the interval
$[\theta_1, \theta_k]$ by
$$
(\f_{\m, \l}(t),q_{\m, \l}(t)) := 
(\f_{\m,\l_i}(t),q_{\m,\l_i}(t)+2\pi(i-1))  \ \ {\rm for}
\  \  t\in [\theta_i, \theta_{i+1}],
$$
where $ (\f_{\m,\l_i}(t), q_{\m,\l_i}(t)) $ are given by 
lemma \ref{lem:connecting}.
The pseudo diffusion solutions
$(\f_{\m,\l}, q_{\m, \l })$ are then continuous functions which 
are true solutions of the equations of motion (\ref{lagraeqmot}) on each
interval $(\theta_i,\theta_{i+1})$, but the time derivatives
$({\dot \f}_{\m,\l}, {\dot q}_{\m, \l })$
may undergo a jump at time $\theta_i$. We have 

\begin{lemma}\label{lem:heter}
If $ \wl = ( \wtilde{\teta}, \wtilde{\vphi}) \in \Lambda_{\m,k} $
is a critical point of $ { \cal F}_\mu $, then
$(\vphi_{\mu ,\wl} (t), q_{\mu , \wl} (t))$
is a  solution of (\ref{lagraeqmot}) in the time interval
$(\wtilde{\teta}_1, \wtilde{\teta}_k)$. Moreover
$ {\dot \vphi}_\mu ( \wtilde{\teta}_1 ) =$ $ \om_I + O(\mu ),$  
$ {\dot \vphi}_\mu ( \wtilde{\teta}_k ) =$  $ \om_F + O(\mu ),$
{\it i.e.} 
$ (\vphi_{\mu, \wtilde{\l}} , q_{\mu,\wtilde{\l}} )$ 
is a diffusion orbit between $\om_I$ and $\om_F$
with diffusion time $T_d  =  |\wtilde{\teta}_k - \wtilde{\teta}_1| $. 
\end{lemma}
\begin{pf} 
By lemma \ref{lem:derG} if $\nabla_{\f_i} {\cal F}_{\m}
(\wl)=0$, then for $ 2 \leq i \leq $ $ k - 1 ,$ 
$\dot{\f}_{\m,\wl}( \wtilde{\theta}_i^-)=$ 
$\dot{\f}_{\m,\wl}(\wtilde{\theta}_i^+)$
and $\dot{\f}_{\m,\wl}(\wtilde{\theta}_1)=\om_I+O(\m)$,
$\dot{\f}_{\m,\wl}(\wtilde{\theta}_k)=\om_F+O(\m)$. Moreover,
if $\nabla_{\f_i} {\cal F}_{\m} (\wl)=0$ and
$\partial_{\theta_i} {\cal F}_{\m} (\wl)=0$ then (for $2\leq i
\leq$ $k-2$), 
$\dot{q}^2_{\m,\wl}(\wtilde{\theta}_i^+) = 
\dot{q}^2_{\m,\wl}(\wtilde{\theta}_i^-)$.
Now, by lemma \ref{lem:connecting} and $(\ref{estQT})$, 
${\dot q}_{\m, \wtilde{\l}} (\wtilde{\theta}_i^\pm)= \dot{q}_0(0)+O(\m)$. Hence
$\dot{q}_{\m,\wl}(\wtilde{\theta}_i^+) =
 \dot{q}_{\m,\wl}(\wtilde{\theta}_i^-)$
and the proof is complete.  
\end{pf}

\section{The approximation of the reduced functional}

In order to prove the existence of critical points of
the reduced action functional 
${\cal F}_\m$ thanks to the properties of the Poincar\'e-Melnikov
primitives $\Gamma(\om, \cdot, \cdot )$ we need an appropriate expression of 
${\cal F}_\m$, see lemma \ref{approxsum}. We shall express ${\cal F}_\m$ as
the sum of a function whose definition contains the 
$\Gamma(\om, \cdot, \cdot )$
(for which we can prove the existence of critical points)
and of a remainder whose derivatives are so small that it 
cannot destroy the critical points of the first function.

The first lemma gives an approximation of $G_\m$ (defined in $(\ref{defiG})$).

\begin{lemma}\label{lem:apprG}  
For $ 0 < \mu \leq \mu_3 $, for $\lambda \in \Lambda_{\mu}$  we have 
\be\label{eq:kinen} 
G_\m(\lambda) =  
\frac{1}{2} \frac{|\vphi^- - \vphi^+|^2}{(\teta^- - \teta^+)}
+\mu \Gamma^s (\om_\l, \theta^+ , \f^+) + \mu \Gamma^u (\o_\l , 
\theta^- , \f^-) - \mu \int_{\theta^+}^{\theta^-} 
f(\ov{\f}(t) , 0,t) \  dt  
+ R_0 (\mu, \l) 
\ee
where 
\be\label{eq:restkinen}
\nabla_\l R_0 (\mu, \l) = 
O \Big( \frac{\mu^2 (1+\mu T_\l^2)}{\b_\l^2} T_\l \Big).
\ee 
\end{lemma}
\begin{pf}
By lemma \ref{lem:connecting},
we can write $\f_{\m,\l}(t)=\ov{\f}(t)+v_{\mu,\l}(t)$,
$q_{\mu , \l} (t)=Q_{T_\l}(t-\theta^+)+w_{\mu,\l}(t)$, where 
$ v_{\mu,\l} (\teta^+ ) = v_{\mu, \l} (\teta^- ) =0 $, 
$ || \dot{v}_{\mu, \l}  ||_{L^{\infty}(\teta^+, \teta^-)} 
= O ( \mu \slash \b_\l ) $,
$|| v_{\mu, \l}  ||_{L^{\infty}(\teta^+, \teta^-)} 
= O ( (\mu \slash \b_\l^2)(1+\mu T_\l^2) )$
and $w_{\mu,\l} (\theta^+)=w_{\mu,\l} (\theta^-)=0$, 
$|| \dot{w}_{\mu,\l}  ||_{L^{\infty}(\teta^+, \teta^-)}+
|| w_{\mu,\l}  ||_{L^{\infty}(\teta^+, \teta^-)} 
= O ( \mu  ) $.

In the following, in order to avoid cumbersome notation, 
we shall use the abbreviations $v,w,Q$ for $v_{\m,\l}, w_{\m,\l},
Q_{T_\l}( \cdot -\theta^+)$, the dependency w.r.t. $\lambda$ and $\m$
being implicit.  We have
$$
\begin{array}{lll}
G_\m(\lambda)&=& \dps \int_{\theta^+}^{\theta^-} 
\frac{1}{2}|\dot{\ovf}(t)|^2 + \dot{\ovf}(t) \cdot \dot{v}(t)+
\frac{1}{2} |\dot{v(t)}|^2 +
\frac{1}{2}\dot{Q}^2(t) + \dot{Q}(t)\dot{w}(t)+
\frac{1}{2} \dot{w}^2(t) \\
&+& \dps  \int_{\theta^+}^{\theta^-}
[1-\cos (Q(t)+w(t))] - \mu f(\ovf(t)+v(t), Q(t)+w(t) , t) \ dt .
\end{array}
$$
Now since $ v(\theta^+) = v(\theta^-) = 0 $ and
$ w (\theta^+) = w(\theta^-) = 0 $,
$ \dps
\int_{\theta^+}^{\theta^-} \dot{\ovf}(t) \cdot \dot{v}(t) \ dt =
\int_{\theta^+}^{\theta^-} \om_\l \cdot \dot{v}(t) \ dt =0
$ and
$ \dps
\int_{\theta^+}^{\theta^-} \dot{Q}(t)\dot{w}(t) \  dt=
\int_{\theta^+}^{\theta^-} - \ddot{Q}(t) w(t) \ dt =
\int_{\theta^+}^{\theta^-} -(\sin Q(t)) w(t) \ dt. 
$
As a result, $G_\m (\lambda)=G^0_\m (\lambda)+ R_1(\lambda)$,
where
$$
G_\m^0(\l)=\int_{\theta^+}^{\theta^-} 
\frac{1}{2}|\dot{\ovf}|^2  + \frac{1}{2}\dot{Q}^2 + 
(1-\cos Q) - \mu f(\ovf, Q , t) ,
$$
$$
R_1(\lambda)=\int_{\theta^+}^{\theta^-} 
\frac{1}{2} |\dot{v}|^2+ \frac{1}{2} \dot{w}^2 +
(\cos Q-\cos (Q+w)-w\sin Q) - \mu f(\ovf+v, Q+w ,
t) + \mu f(\ovf, Q , t).
$$
We shall first prove that $|\nabla R_1 |=
O\Big( \frac{\mu^2 (1+\mu T_\l^2)}{\b_\l^2} T_\l \Big)$.
We have $\partial_{\theta^+} R_1=r_s1+r_2+r_3+r_4+r_5+r_6$ where
$$
\begin{array}{l}
r_1:=\dps \int_{\theta^+}^{\theta^-} 
\dot{v} \cdot \frac{d}{dt}(\partial_{\theta^+} v)
 - \mu \partial_{\f} f(\ovf+v, Q+w ,
t) \cdot (\partial_{\theta^+} v),   \\
r_2:=\dps \int_{\theta^+}^{\theta^-} 
\dot{w} \frac{d}{dt}(\partial_{\theta^+} w) +
\Big[ \sin(Q+w)-\sin Q - \mu \partial_{q} f(\ovf+v, Q+w ,
t) \Big] (\partial_{\theta^+} w),  \\
r_3:=\dps \int_{\theta^+}^{\theta^-} 
(-\sin Q+\sin (Q+w)-w\cos Q)
\partial_{\theta^+}Q,  \\
r_4:=\dps\mu \int_{\theta^+}^{\theta^-} \Big[ 
\partial_{\f} f(\ovf, Q,t) - \partial_{\f} f(\ovf+v, Q+w ,t) \Big] \cdot
\partial_{\theta^+} \ovf,   \\
r_5:=\dps \mu \int_{\theta^+}^{\theta^-} \Big[ \partial_{q} f(\ovf, Q,t)
- \partial_{q} f(\ovf+v, Q+w ,t) \Big] \partial_{\theta^+} Q,  \\
r_6:=\dps -\frac{1}{2} |\dot{v}(\theta^+)|^2 -\frac{1}{2} \dot{w}(\theta^+)^2. 
\end{array}
$$
Now $v$ and $w$ satisfy
$$
\left\{
\begin{array}{rll}
-\ddot{v}(t)&=& \mu  \partial_{\f}f(\ovf(t)+v(t),Q(t)+w(t),t)\\
-\ddot{w}(t)+\sin(Q(t)+w(t))&=& \mu  \partial_{q} f(\ovf(t)+v(t), Q(t)+w(t)
,t) + \sin Q(t).
\end{array}
\right.
$$
Moreover, deriving w.r.t. $\theta^+$ the equality
$v(\theta^+)=0$ we obtain that $(\partial_{\theta^+}v)(\theta^+)=
-\dot{v}(\theta^+)$. Similarly $ (\partial_{\theta^+}w)(\theta^+)=
-\dot{w}(\theta^+)$, $(\partial_{\theta^+}v)(\theta^-)=0$ and 
$(\partial_{\theta^+}w)(\theta^-)=0$. Therefore an
integration by parts gives $r_1=|\dot{v}(\theta^+)|^2$,
$r_2=\dot{w}(\theta^+)^2$ hence  
$|r_1|+|r_2|=O(\mu^2/\beta^2)$.

By the properties of $Q_T$, $\partial_{\theta^+}Q$ is bounded in
the interval $[\theta^+,\theta^-]$ by a constant independent of
$\l$. Moreover $ -\sin Q(t)+\sin (Q(t)+w(t))-w(t)\cos
Q(t)=O(w(t)^2)$. Therefore $r_3=O(\mu^2 T)$.

We have also, for some positive constant $c$,
$$
|r_4|+|r_5|\leq c \mu T \Big[ \sup_{t\in [\theta^+,\theta^-]}
|\partial_{\theta^+}Q(t)|+ |\partial_{\theta^+}\ovf(t)| \Big]
\Big[ \sup_{t\in [\theta^+,\theta^-]}(|v(t)|+|w(t)|) \Big].
$$
Since $\partial_{\theta^+}\ovf$ is bounded independently of
$\lambda$,  we have by lemma \ref{lem:connecting}
$|r_4|+|r_5|=O \Big( \frac{\mu^2 (1+\mu T_\l^2)}{\b_\l^2} T_\l
\Big)$.
Still by lemma \ref{lem:connecting}, $r_6 = O(\mu^2/\b^2)$. 
The estimate of the other derivatives of $R_1$ is obtained 
in the same way.
\\[2mm] indent
We now develop $ G_\m^0(\lambda)$ as  
$$
G_\m^0(\lambda)=\frac{1}{2} \frac{|\vphi^- - \vphi^+|^2}{(\teta^- - \teta^+)}
+\mu \Gamma^s (\om_\l, \theta^+ , \f^+) + \mu \Gamma^u (\o_\l , 
\theta^- , \f^-) - \mu \int_{\theta^+}^{\theta^-} 
f(\ov{\f}(t) , 0,t) \  dt +R_2(\l)+R_3(\l),
$$
where
\be\label{azpend}
R_2(\l)= \int_{\theta^+}^{\theta^-} \frac{1}{2}\dot{Q}^2(t) + 
(1-\cos Q(t)) \ dt = \int_{0}^{T_\l} \frac{1}{2} {\dot Q}_{T_\l}^2(t) + 
(1-\cos Q_{T_\l}(t)) \ dt,
\ee
$$
R_3 ( \l ) = \int_{\theta^+}^{\theta^-}  - 
\mu \Big[ (f(\ovf(t), Q(t) , t)
-f(\ovf(t),0,t) \Big] \ dt -\mu \Gamma^s (\om_\l, \theta^+ , \f^+)- 
\mu \Gamma^u (\o_\l , \theta^- , \f^-).
$$
There remains to prove estimate (\ref{eq:restkinen}) for
$ \nabla R_2 $ and $\nabla R_3$.
By (\ref{azpend})
$\partial_{\f^{\pm}} R_2=0$ and $\partial_{\theta^+} R_2
(\l)=-\partial_{\theta^-}R_2(\l)$ is the energy of the 
$T_\lambda$-periodic solution $Q_{T_\l}$
of the pendulum equation. Now this energy is $O(e^{-c_2 T_\l})$. Hence
(provided $C_1$ is large enough) $|\nabla R_2 (\l)|=O(\mu^2)$.

In order to estimate the derivatives of $R_3$, let us define 
$g(\f,q,t):=f(\f,q,t)-f(\f,0,t)$. We have
$$
R_3(\lambda)= \int_{\theta^+}^{\theta^-} -\mu g(\ovf(t),Q(t),t) \
dt -  \mu \Gamma^s(\o_\l, \theta^+,\f^+) - \mu \Gamma^u(\o_\l,
\theta^-,\f^-) =\mu (a_3(\l)+b_3(\l))
$$
where
$$
a_3(\l) := -\int_0^{T_\l /2}g(\o_\l t+ \f^+, Q_{T_\l}(t) ,t +\teta^+ ) \ dt
+ \int_0^{\infty} g(\o_\l t+ \f^+, q_0(t) ,t + \teta^+ ) \ dt ,
$$
$$
b_3(\l) := -\int_{-T_\l /2}^0 g(\o_\l t+ \f^-, Q_{T_\l}(t+T_\l) ,
t + \teta^-) \ dt
+ \int_{-\infty}^0  g(\o_\l t+ \f^-, q_0(t) ,t + \teta^-) \ dt.
$$
We have 
$$
a_3(\l)=
-\int_0^{T_\l /2} \Big[ g(\o_\l t+ \f^+, Q_{T_\l}(t) ,t + \teta^+) -
 g(\o_\l t+ \f^+, q_0(t) ,t + \teta^+ ) \Big] + 
\int_{T_\l /2}^{\infty} g(\o_\l t+ \f^+, q_0(t) ,t + \teta^+).
$$
Recalling that $\sup_{t\in (0,T/2)}|\partial_T Q_T (t)|=O(e^{-c_2T})$,
$\sup_{t\in (0,T/2)}| Q_T (t)-q_0(t)|=O(e^{-c_2T})$,
it is easy to see that the derivatives
of the first integral are $O(T_\l e^{-c_2 T_\l})=O(\mu)$
(still provided $C_1$ is large enough).
Moreover, using that $(|g(\o_\l t+\f^+,q_0(t),t)|+
|\partial_{\f} g(\o_\l t+\f^+,q_0(t),t)|+
|\partial_{t} g(\o_\l t+\f^+,q_0(t),t)|)=O(q_0(t)-2\pi)
=O(e^{-c_2 t})$ for $t\in (T_\l /2,+\infty)$, we find that 
the derivatives of the second integral are $O(\mu)$ as well. Hence 
$| \nabla a_3 (\l) | = O( \mu )$. The same estimate holds for
$b_3$. We then conclude that  
$\nabla R_3 ( \l )= O (\m^2) $, which completes the proof of lemma
\ref{lem:apprG}. 
\end{pf}

In section \ref{sec:difforbit}  we will look for a critical point of 
${\cal F}_\m$ in the set
\be \label{ens}
E := \Big\{ \l=( \theta_1,\ldots,\theta_k, \f_1,\ldots,\f_k ) \in
{\bf R}^{k} \times {\bf R}^{kd}  \ \Big| \  
\teta_i =  \ov{\teta}_i + b_i, \  \
\vphi_i =  \ov{\vphi}_i + a_i, \ \
|b_i|\leq 2\pi, \ |a_i| \leq 2\pi  \Big\},
\ee
where $k,\ov{\f}_i,\ov{\theta}_i$ will be defined in
section \ref{sec:pseudodiff}. It will result that 
$E \subset \Lambda_{\mu,k}$ (for some $ \b > 0 $ depending 
on the curve $ \gamma $). In particular, 
for all $ \lambda \in E$ 
\be\label{truesep}
C_1 | \ln \mu | \leq \teta_{i+1} - \teta_i 
< \frac{ C_0 \b_i }{\mu}, \quad \forall  i =1, \ldots, k -1,
\ee
where $\beta_i := \b_{\l_i} := \beta (\om_i)$ and 
$ \om_i := \om_{\l_i} := (\f_{i+1}-\f_i) \slash (\theta_{i+1}-\theta_i)$.
Moreover we will assume (see (\ref{onejump}))  
\begin{equation} \label{propphi}
|\ov{\om}_{i+1}-\ov{\om}_i| \leq \rho \m \qquad {\rm where}\ \ 
\ov{\om}_i := \frac{\ov{\f}_{i+1}-\ov{\f}_i}{\ov{\theta}_{i+1}
-\ov{\theta}_i} \ (1\leq i\leq k-1), \ \ \om_0:=\om_I, \  \
\om_k:=\om_F
\end{equation}
and $ \rho > 0 $ is a small constant to be chosen later
(see (\ref{defrho})).
For the time being, assuming (\ref{truesep}) and (\ref{propphi}),
we want to give a suitable expression of ${\cal F}_\mu $ in $E$. 
By lemma \ref{lem:apprG}, for $ \lambda \in E$, we have
\be \label{formulaFmu1}
\begin{array}{lll}
{\cal F}_{\m} ( \l ) & = & \dps \sum_{i=1}^{k-1}\frac{1}{2} 
\frac{|\f_{i+1}-\f_i|^2}{\theta_{i+1}-\theta_i}+ 
\om_I \vphi_1 - \om_F \vphi_k - \frac{|\om_I|^2}{2} \teta_1 
+ \frac{|\om_F|^2}{2} \teta_k  \\
 &+&\dps \sum_{i=1}^{k} \mu \Big( \Gamma^u (\om_{i-1} , \theta_i,\f_i)+
\Gamma^s (\om_{i} , \theta_i,\f_i) \Big) 
+ \m F(\om_I,\theta_1,\f_1) \\
&- & \dps \sum_{i=1}^{k-1}\m \dps \int_{\theta_i}^{\theta_{i+1}}
f(\om_i (t-\theta_i)+\f_i,0,t) \ dt 
- \m F(\om_F , \theta_k , \f_k) + \dps \sum_{i=1}^{k-1} 
R_0(\mu,\l_i),
\end{array}
\ee
where  $|\nabla_\l R_0 ( \mu,\l )| $ satisfies (\ref{eq:restkinen}).
We shall write ${\cal F}_\mu $  in an appropriate form
thanks to the following lemmas. The first one says
how close the ``mean frequencies'' $ \om_i $ are to the unperturbed 
$\ov{\om}_i$.
\begin{lemma}\label{lem:vici}
Let 
$\l = ( \teta_1, \ldots ,\teta_{k}, \vphi_1, \ldots , \vphi_{k})$ 
belong to $E$. Then 
\be\label{sepfre}
|\om_i - \ov{\om}_i | = 
O \Big( \frac{1 }{ \teta_{i+1} - \teta_i} \Big) = 
O \Big( \frac{1 }{ | \ln \mu |} \Big). 
\ee
Moreover 
\be \label{form:Gamma} 
\Gamma^u (\om_{i-1} , \theta_i,\f_i)+
\Gamma^s (\om_{i} , \theta_i,\f_i)=\Gamma
(\ov{\om}_i,\theta_i,\f_i) + R_4 (\lambda_i) , \ \ {\rm where} \ \
\nabla R_4=O(1/|\ln \m |).
\ee
\end{lemma}

\begin{pf}
Set $ \Delta \teta_i :=  \teta_{i+1} - \teta_i $, 
$ \Delta a_i :=  a_{i+1} - a_i $ and $ \Delta b_i :=  b_{i+1} - b_i $.
By an elementary computation we get
$ \om_i - \ov{\om}_i  = - \ov{\om}_i \Delta b_i \slash \Delta \teta_i  +
\Delta a_i \slash \Delta \teta_i. $
By the definition of $E$ and (\ref{truesep}), estimate (\ref{sepfre}) follows.

From the definition of $\Gamma^u, \Gamma^s$ 
and the exponential decay of $q_0$ it results that
$\partial_{\om}\Gamma^{u,s}$ is bounded by a uniform constant,
as well as its partial derivatives. Hence (\ref{form:Gamma})
is a straightforward consequence of (\ref{sepfre}) 
and of (\ref{propphi}).
\end{pf}

\begin{lemma} \label{lem:toriest}
For $ 0 < \mu \leq \mu_4 $ 
\be\label{nearto}
\mu F(\om_I, \teta_1, \vphi_1) - 
\sum_{i=1}^{k} \mu  \int_{\teta_i}^{\teta_{i+1}} 
f(\om_i(t-\theta_i)+\f_i, 0, t)  \ dt - \mu F(\om_F, \teta_k, \vphi_k) =  
 \sum_{i=1}^k R^i_5 (\mu,\lambda_{i-1},\lambda_i),
\ee
where, for all $i$ \footnote{In the cases $i=1,i=k$ we only have 
$ R_5^1=R_5^1(\mu,\teta_1,\vphi_1,\teta_2,\vphi_2) $ and 
$ R_5^k = R_5^k ( \mu,\teta_{k-1}, \vphi_{k-1},\teta_k,\vphi_k).$} 
\be\label{restnearto}
\nabla R_5^i 
( \mu, \teta_{i-1}, \vphi_{i-1}, \teta_i, \vphi_i,
\teta_{i+1}, \vphi_{i+1})
= O \Big( \frac{\mu}{\b_{i-1}^2(\teta_{i} - \teta_{i-1}) }+
  \frac{\mu}{\b_{i}^2(\teta_{i+1} - \teta_{i}) } + \frac{\mu
|\beta_i -\beta_{i-1}|}{\beta_{i-1}\beta_i} \Big).
\ee
\end{lemma}

\begin{pf}
We have
\begin{eqnarray*}
& \empty & 
- \int_{\teta_i}^{\teta_{i+1}} f(\vphi_i + \om_i (t-\teta_i), 0, t) \ dt 
=  F( \om_i, \teta_{i+1}, \vphi_{i+1}) - F(\om_i, \teta_i , \vphi_i) \\
& \empty & \qquad = \Big( F( \om_i, \teta_{i+1} , \vphi_{i+1} ) - 
F( \om_{i-1}, \teta_i , \vphi_i) \Big) 
+ \Big( F(\om_{i-1}, \teta_i ,  \vphi_i) - F(\om_{i}, \teta_i , \vphi_i) \Big),
\end{eqnarray*}
where $F(\om, \cdot, \cdot)$ is defined in (\ref{risoterm}). 
We obtain
$$
\mu F( \om_I, \teta_1, \vphi_1) -   
\sum_{i=1}^{k-1} \mu  \int_{\teta_i}^{\teta_{i+1}} 
f(\vphi_i + \om_i (t-\teta_i), 0, t) \ dt 
-  \mu F(\om_F, \teta_k, \vphi_k)  
= \sum_{i=1}^k {R}_5^i 
$$
where
\begin{eqnarray*}
R_5^i := R_5^i ( \mu, \teta_{i-1} ,\vphi_{i-1}
, \teta_i, \vphi_i, \teta_{i+1}, \vphi_{i+1}) & := &  
\mu\Big( F(\om_{i-1},  \teta_i ,\vphi_i) -  
F(\om_i, \teta_i , \vphi_i) \Big) \\
& = &  
- \mu\!\! \sum_{0 < |(n,l)| \leq N}\!\! f_{n,l} \frac{e^{{\rm i}
(n \cdot \vphi_i + l \teta_i)}}{{\rm i}}
\Big( \frac{1}{( n \cdot \om_{i-1} + l )} - \frac{1}{( n \cdot \om_i + l)} 
\Big) 
\end{eqnarray*}
Now we prove (\ref{restnearto}).
Let us consider for example
$ \partial_{\teta_i} {R}_5^i $. We have 
\begin{eqnarray}
\partial_{\teta_i} {R}_5^i & = & \mu \partial_{\teta_i}  
\Big( F(\om_{i-1}, \teta_i ,  \vphi_i) -  F(\om_i, \teta_i , \vphi_i) \Big)
\nonumber \\
& = & \mu\Big( \partial_\om F(\om_{i-1}, \teta_i ,  \vphi_i) 
.\frac{- \om_{i-1}}{(\teta_i - \teta_{i-1}) } -
\partial_\om F ( \om_i, \teta_i , \vphi_i) 
.\frac{\om_i }{(\teta_{i+1} - \teta_i) } \Big) \nonumber \\
& - & 
\mu\Big( \sum_{0< | (n,l) | \leq N}
f_{n,l} l  e^{{\rm i}(n \cdot \vphi_i + l \teta_i)} (
\frac{1}{(n \cdot \om_{i-1} + l )} -  \frac{1}{(n \cdot \om_{i} + l )}) \Big), \label{r4}
\end{eqnarray}
where
\be\label{derom}
\partial_\om F(\om, \teta_0, \vphi_0 ) = 
\sum_{0< | (n,l) | \leq N}
f_{n,l} \frac{n e^{{\rm i}(n \cdot \vphi_0 + l \teta_0)}}{{\rm i}(n \cdot \om + l )^2},
\ee
Estimate (\ref{restnearto}) follows immediately
from (\ref{r4}) and (\ref{derom}). The other partial derivatives
of $R_5^i$ can be estimated similarly.  
\end{pf}

Finally, to get a suitable expression of ${\cal F}_{\m}$,
we find convenient to introduce 
 coordinates $ ( b, c ) \in {\bf R}^{(1+d)k}$ defined by
(\ref{ens}) and  
\be\label{coordc}
c_i = a_i - \ov{\om}_i b_i, \qquad \forall i = 1, \ldots, k,
\ee
(we are just performing a linear change of coordinates adapted 
to the direction of the unperturbed flow at each $i$-transition    
$(b_i , a_i ) = b_i (1, \ov{\om}_i ) + ( 0, c_i )$). 

\begin{lemma} \label{lem:hjje}
We have 
\begin{eqnarray}
\sum_{i=1}^{k-1} \frac{1}{2} \frac{ | \vphi_{i+1} - 
\vphi_i|^2}{ (\teta_{i+1} - \teta_i) } +
\om_I \vphi_1 - \om_F \vphi_k - \frac{|\om_I|^2}{2} \teta_1 
+ \frac{|\om_F|^2}{2} \teta_k  & = &
\frac{1}{2}  \sum_{i=1}^{k-1} 
\frac{| c_{i+1} - c_i |^2}{ \Delta \ov{\teta}_i + (b_{i+1} - b_i )} \\ 
& + & 
\sum_{i=1}^k
R_6^i (\mu,\teta_i, \vphi_i,\teta_{i+1} ,\vphi_{i+1}),
\nonumber
\end{eqnarray}
where $ \Delta \ov{\teta}_i := \ov{\teta}_{i+1} - \ov{\teta}_i $ 
and\footnote{For $i=k$ we have $R_6^k=R_6^k(\mu,\teta_k , \vphi_k).$} 
\be\label{hetenkin} 
\nabla R_6^i  
(\mu,\teta_{i-1}, \vphi_{i-1},\teta_i,\vphi_i,\teta_{i+1} ,\vphi_{i+1} )
= O ( \Delta \ov{\om}_i ) = O ( \rho \mu ).
\ee
\end{lemma}

\begin{pf}
Let $\{ \gamma_i \}_{i=1, \ldots, k-1} $ be defined by 
$\vphi_{i+1} - \vphi_i = \ov{\om}_i (\teta_{i+1} - \teta_i ) +
\gamma_i$. We can write
 $ \om_I \vphi_1 - \om_F \vphi_k$ as 
\begin{eqnarray}
\om_I \vphi_1 - \om_F \vphi_k & = & \sum_{i=1}^{k-1}  
\Big( ( \ov{\om}_{i-1} - \ov{\om}_i ) \vphi_i -
\ov{\om}_i ( \vphi_{i+1} - \vphi_i ) \Big) + \vphi_k 
(\ov{\om}_{k-1} - \om_F ) \nonumber \\
& = & \sum_{i=1}^{k-1}  
\Big( ( \ov{\om}_{i-1} - \ov{\om}_i ) \vphi_i -
|\ov{\om}_i|^2 ( \teta_{i+1} - \teta_i ) - 
\ov{\om}_i \gamma_i \Big)
+ \vphi_k (\ov{\om}_{k-1} - \om_F ). \label{heterjump}
\end{eqnarray}
We can also write 
\be\label{enjump}
- \frac{|\om_I|^2}{2} \teta_1 + \frac{|\om_F|^2}{2} \teta_k 
=  \sum_{i=1}^{k-1}  
\Big( (\frac{|\ov{\om}_i|^2}{2} - \frac{|\ov{\om}_{i-1}|^2}{2} ) \teta_i +
\frac{|\ov{\om}_i|^2}{2} ( \teta_{i+1} - \teta_i ) \Big) +
\Big( \frac{|\om_F|^2}{2} - \frac{|\ov{\om}_{k-1}|^2}{2} \Big) \teta_k, 
\ee
\be\label{kinetic}
\sum_{i=1}^{k-1}
\frac{1}{2} \frac{|\vphi_{i+1} - \vphi_i|^2}{(\teta_{i+1} - \teta_i )} = 
\sum_{i=1}^{k-1} \frac{|\ov{\om}_i|^2}{2} ( \teta_{i+1} - \teta_i ) + 
\frac{1}{2} \frac{|\gamma_i|^2}{(\teta_{i+1} - \teta_i )} + 
\ov{\om}_i \gamma_i.
\ee
Summing (\ref{heterjump}), (\ref{enjump}) and (\ref{kinetic}) we get
$$
\sum_{i=1}^{k-1} \frac{1}{2} \frac{|\vphi_{i+1} - 
\vphi_i|^2}{ (\teta_{i+1} - \teta_i )} + 
\om_I \vphi_1 - \om_F \vphi_k - \frac{|\om_I|^2}{2} \teta_1  
+ \frac{|\om_F|^2}{2} \teta_k  =
\sum_{i=1}^{k-1} \frac{1}{2} \frac{|\gamma_i|^2}{(\teta_{i+1} - \teta_i)} +$$
\be\label{svill}
\sum_{i=1}^{k-1} 
\Big( \frac{|\ov{\om}_i|^2}{2} - \frac{|\ov{\om}_{i-1}|^2}{2} \Big) \teta_i
+ (\ov{\om}_{i-1} - \ov{\om}_i ) \vphi_i +
\vphi_k (\ov{\om}_{k-1} - \om_F ) +
\Big( \frac{|\om_F|^2}{2} - \frac{|\ov{\om}_{k-1}|^2}{2} \Big) \teta_k.
\ee
Substituting $\ov{\vphi}_i + a_i $ for $\vphi_i$ and 
$\ov{\teta}_i + b_i $ for $\theta_i$, we get
$\gamma_i = ( a_{i+1} - a_i ) - \ov{\om}_i ( b_{i+1} - b_i )$.
Moreover the non constant terms in the right handside of (\ref{svill}) 
(i.e. those depending on $a_i, b_i $) are the first one and
$$
\sum_{i=1}^k  (\ov{\om}_{i-1} - \ov{\om}_i) a_i +
\Big( \frac{|\ov{\om}_i|^2}{2} - \frac{|\ov{\om}_{i-1}|^2}{2} \Big) b_i =:
\sum_{i=1}^k R^i (\mu, \teta_i , \vphi_i)
$$
with  $\nabla R^i (\mu, \teta_i ,\vphi_i )  = O ( \Delta \ov{\om}_i ).$
Finally, expressing $\gamma_i $ in terms of $(b_i, c_i) $ we get
$\gamma_i = (a_{i+1} - a_i ) - \ov{\om}_i ( b_{i+1} - b_i ) = 
( c_{i+1} - c_i ) + b_{i+1} \Delta \ov{\om}_i$
and then from (\ref{svill}), developing the square, we get
(\ref{hetenkin}).
\end{pf}

From (\ref{formulaFmu1}) lemmas \ref{lem:vici}, \ref{lem:toriest} and
\ref{lem:hjje} we obtain the expression of ${\cal F}_\mu $ 
in the new coordinates $(b,c)$ required to apply 
the variational argument of section \ref{sec:difforbit}.

\begin{lemma}\label{approxsum}
There exists $ \mu_5, C_2 > 0 $ such that 
$ \forall\  0 < \mu \leq \mu_5 $,   if 
\be\label{betai}
\b_i \geq  C_2
\max\ \{\  \mu^{1/2}(\teta_{i+1}-\teta_i)^{1/2},\ 
\mu(\teta_{i+1}-\teta_i)^{3/2},\ 
(\teta_{i+1}-\teta_i)^{-1/2}\ \}
\ee
then 
\begin{eqnarray} 
{\cal F}_\mu ( b, c) & = &  \frac{1}{2} \sum_{i=1}^{k-1} 
\frac{|c_{i+1} - c_i|^2}{\Delta \ov{\teta}_i + (b_{i+1} - b_i)} + 
\mu \sum_{i=1}^k \Gamma (\ov{\om}_i, \ov{\theta}_i + b_i, 
\ov{\f}_i+\ov{\om}_i  b_i + c_i )+ R_7(b,c),\label{eq:somma}  \\ 
R_7(b,c) & :=  &  \sum_{i=1}^k R_7^i (\mu, b_{i-1}, c_{i-1},  b_i, 
c_i, b_{i+1}, c_{i+1}), \label{R7}
\end{eqnarray}
where\footnote{In the cases $i=1,i=k$ 
we  have $R_7^1=R_7^1(\mu,\teta_1, \vphi_1,\teta_2, \vphi_2)$ and 
$R_7^k=R_7^k(\mu,\teta_{k-1},\vphi_{k-1},\teta_k,\vphi_k).$}
\be\label{smalldef}
| \nabla  R_7^i | \leq C_2 \rho \mu.  
\ee
\end{lemma}

\begin{pf}
It is easy to see that (\ref{propphi}), (\ref{sepfre}) and 
(\ref{betai}) imply (provided $\mu$ is small enough)
that
\be \label{estagain}
\frac{\beta_{i-1}}{2} \leq \beta_i \leq 2 \beta_{i-1},  \quad  
\quad  |\beta_i-\beta_{i-1}| = O\Big(\frac{1}{\theta_i-\theta_{i-1}}
+ \frac{1}{\theta_{i+1}-\theta_i} + \mu \Big).
\ee
Noting that $\partial_{c_i} = \partial_{\vphi_i} $ and 
$\partial_{b_i} = \ov{\om}_i \partial_{\vphi_i} + \partial_{\teta_i} $,
estimate (\ref{smalldef}) follows 
from (\ref{eq:restkinen}), (\ref{form:Gamma}),
(\ref{restnearto}), (\ref{estagain}) and (\ref{hetenkin}).
\end{pf}

\section{Ergodization times}\label{sec:erg}

In order to define $ \ov{\f}_i, \ov{\theta}_i $ $ ( 1 \leq i \leq k )$
we need some results, stated in this section,
on the ergodization time of the torus 
${\bf T}^l := { \bf R}^l / {\bf Z}^l $
for linear flows possibly resonant but only at a ``sufficiently high level'',
\\[1mm]
\indent
Let $ \Omega \in {\bf R}^l $; it is well known that,   
if $ \Omega \cdot p \neq 0 $, 
$ \forall p \in {\bf Z}^l \setminus \{ 0\} ,$ then the trajectories
of the linear flow 
$\{ \Omega t + A \}_{t \in {\bf R}} $ are  dense on ${\bf T}^l $
for any initial point $ A \in  { \bf T}^l $. 
It is also intuitively clear that the trajectories of the linear flow 
$\{ \Omega t  + A \}_{t \in {\bf R}} $ 
will make an arbitrarly fine $\d$-net ($\d >0$)
if $ \Omega $ is resonant only at a sufficiently high level, namely if
$ \Omega \cdot p \neq 0 $, $ \forall p \in {\bf Z}^l $ with 
$ 0 < |p| \leq  M(\delta) $ for some  large enough $M(\delta)$.  
Let us make more precise and quantitative these considerations.

For any $ \Omega \in {\bf R}^l $ define the ergodization time 
$ T ( \Omega , \d )$ required to fill $ {\bf T}^l $ within $ \d > 0 $
as 
$$
T(\Omega, \delta)= \inf \Big\{ t\in {\bf R}_+ \ \Big| \ 
\forall x \in {\bf R}^l,
\  d(x, A+[0,t] \Omega + {\bf Z}^l) \leq \delta  \Big\},
$$ 
where $d$ is the Euclidean distance and $A$ some point of 
${\bf R}^l$.
$T(\Omega , \delta)$ is clearly independent of the choice of
$A$.  Above and in what follows, 
$\inf E$ is equal to $+\infty$ if $E$ is empty.  For $R>0$ let
$$
\alpha (\Omega , R)= \inf \Big\{ | p \cdot \Omega| \ \Big| 
\ p\in {\bf Z}^l, \ \ p \neq 0 \ , \ |p| \leq R \Big\}.
$$
\begin{theorem}\label{lem:ergo1}
$\forall l \in {\bf N} $ there exists a positive constant $ a_l $ 
such that, $\forall \Omega \in {\bf R}^l$,  
$\forall \delta >0$,  
$T(\Omega,\delta) \leq (\alpha (\Omega, a_l/\delta))^{-1}$. 
Moreover $ T ( \Omega, \delta ) \geq ( 1 / 4 ) 
\alpha( \Omega, 1/4\delta)^{-1} $.
\end{theorem}

In the above Theorem $\alpha^{-1}$ is equal to 0 if
$\alpha=+\infty$ and to $+\infty$ if $ \alpha=0 $. 

\begin{remark}\label{rem:BGW}
Assume that $ \Omega $ is a $C$-$\tau$ Diophantine vector, i.e.
there exist $ C > 0 $ and $\tau \geq l-1 $ such that
$ \forall k \in {\bf Z}^l$  
$|k \cdot \Omega| \geq C / |k|^{\tau}$.
Then $ \alpha (\Omega, R) \geq
C/R^{\tau}$ and so $ T (\Omega, \delta) \leq
a_l^{\tau}/ C\delta^{\tau}$. This estimate was
proved in Theorem D of \cite{BGW}. 
Also Theorem B of \cite{BGW} is an easy consequence of Theorem \ref{lem:ergo1}.
\end{remark}

Theorem \ref{lem:ergo1} is a direct consequence of more general 
statements, see Theorem \ref{lem:ergo2} and remark \ref{rem:nonerg}.
Let us introduce first some notations.
Let $\Lambda$ be a lattice of ${\bf R}^l$, {\it i.e.} a discrete
subgroup of ${\bf R}^l$ such that ${\bf R}^l/\Lambda$ has finite volume.
For all $\Omega \in {\bf R}^l$ we define
$$
T(\Lambda,\Omega, \delta)= \inf \Big\{ t\in {\bf R}_+ \ \Big| 
\ \forall x \in {\bf R}^l \  d(x, [0,t] \Omega + \Lambda) \leq \delta  \Big\}
$$ 
($T(\Lambda, \Omega , \delta)$ is the time required to have a 
$\delta$-net of the torus ${\bf R}^l / \Lambda$ endowed with the
metric inherited from ${\bf R}^l$). 
For $R>0$, let 
$$
\Lambda^*= \Big\{ p \in {\bf R}^l \ \Big| \  \forall \lambda \in \Lambda, \ 
p \cdot \lambda \in {\bf Z} \Big\} \quad {\rm and}
\quad \Lambda^*_R = \Big\{ p \in \Lambda^* \ \Big| \  0<|p|\leq R \Big\}
$$
($\Lambda^*$ is a lattice of ${\bf R}^l$ which is conjugated to 
$\Lambda$). We define
$$
\alpha(\Lambda,\Omega, R)= 
\inf \Big\{ |p \cdot \Omega | \ \Big| \  p \in \Lambda^*_R \Big\}.
$$
The following result holds:
\begin{theorem}\label{lem:ergo2}
$\forall l \in {\bf N}$ there exists a positive constant $a_l$ 
such that, for all lattice $\Lambda$ of ${\bf R}^l$, 
$\forall \Omega \in {\bf R}^l$, $\forall \delta >0 $,  
$T(\Lambda,\Omega,\delta) \leq (\alpha (\Lambda,\Omega, a_l/\delta))^{-1}$.
\end{theorem}

\begin{remark}\label{rem:nonerg}
It is fairly obvious that $T(\Lambda, \Omega, \delta ) \geq (1/4) 
\alpha ( \Lambda, \Omega, 1/4\delta)^{-1}$. Indeed, assume that 
$\Lambda^*_{1/4\delta} \neq \emptyset $ and let $ p \in
\Lambda^*_{1/4\delta} $ be such that $ p \cdot \Omega=
\alpha:=$ $\alpha ( \Lambda, \Omega, 1/ 4 \delta)$. Let $x \in {\bf R}^l$
satisfy $ p \cdot x = 1/2 $. Then $ \forall t \in [0, 1/4\alpha)$,
$\forall \lambda \in \Lambda$, 
$$
|x-(t\Omega+\lambda)| \geq \frac{|p\cdot (x- t \Omega - \lambda)| }{|p|} \geq 
4 \delta | p \cdot x- t p \cdot \Omega - p\cdot \lambda|,
$$
and $p\cdot x-p \cdot \lambda \in (1/2)+{\bf Z}$, whereas 
$|tp \cdot \Omega| = t \alpha < 1/4 $. 
Hence $| x - ( t \Omega + \lambda ) | > \delta $.
\end{remark}

In the next section we will apply  Theorem \ref{lem:ergo1} when 
$ \O = (\om, 1) \in {\bf R}^{d+1} $.
The proof of Theorem \ref{lem:ergo2} is given in the Appendix.  
We could give an explicit 
expression of $ a_l $. However it is not useful for our purpose and
the constants $ a_l $ which can be derived from our proof 
are certainly far from being optimal.

\section{The unperturbed pseudo-diffusion orbit} \label{sec:pseudodiff}

Consider the set $Q_M$ of ``non-ergodizing frequencies'' 
$$
Q_M := \Big\{ \om  \in {\bf R}^d \ \Big| \ 
\exists ( n, l) \in {\bf Z}^{d+1} \ {\rm with } \
0 < |( n, l) | \leq M, \ {\rm and} \ \om \cdot n + l = 0 \Big\} = 
\bigcup_{h\in S_M}E_h
$$	
where
$S_M:=\{h=(n,l)\in({\bf Z}^d\setminus \{ 0\})\times{\bf N}\ | \  0<|h|\leq M,
\ h\neq jh',\forall\, j\in{\bf Z},h'\in({\bf Z}^d\setminus \{ 0\})\times{\bf N} \}$
and
$ E_h=E_{n,l} := \{ \om \in {\bf R}^d \ | \ (\om,1) \cdot h 
=\om\cdot n+l = 0 \}. $
By Theorem \ref{lem:ergo1} (or Theorem \ref{lem:ergo2}, with 
$\Lambda = 2\pi {\bf Z}^{d+1}$), for $\delta>0$, if  
$\om$ belongs to 
\be
Q_M^c = \Big\{ \om \in {\bf R}^d \ \Big| \ \om \cdot n + l \neq 0, \ 
\forall 0 < |(n, l ) | \leq M \Big\},
\ee
with $M=8\pi a_{d+1}/\delta$, then the flow of  $(\om,1)$ provides
a $\delta/4$-net of the torus $ {\bf T}^{d+1} $.
 
Moreover if $ \om \notin Q_M $ then for all $(n,l)\in {\bf Z}^d 
\backslash \{ 0\} \times {\bf Z}$,
\be \label{distQM}
|n \cdot \om +l|=|n| {\rm dist} (\om,E_{n,l}) \geq {\rm dist}
(\om,E_{n,l}) \geq  {\rm dist} (\om, Q_M ) > 0.
\ee
 
By Theorem \ref{lem:ergo1} (or Theorem \ref{lem:ergo2}), we deduce from $(\ref{distQM})$  the estimate 
\be\label{timerg}
T( ( \om, 1 ), \d/4 ) \leq \frac{2\pi}{{\rm dist}(\om, Q_M )}
\ee
which  measures the divergence of the ergodization time $T( (\om,1) , \d )$ 
 as $\om $ approaches the  set $ Q_M $.

\begin{definition}\label{def:adm}
Given $ M > 0,$ a connected component ${\cal C}$ of ${\cal
D}_N^c$ and $\om_I , \om_F \in{\cal C}$, 
we say that an embedding $\g\in C^2([0,L],{\cal C})$ is a 
{\rm $Q_M$-admissible connecting curve between $\om_I$ and $\om_F$}
if the following properties are satisfied:
\begin{itemize}
\item[(a)]
$\g(0)=\om_I,$ $\g(L)=\om_F,$ $|\dot{\g}(s)|=1$ $\forall\, s\in (0,L),$

\item[(b)]
$\forall h=(n,l) \in S_M$, $\forall s\in [0,L]$ such that
$\gamma (s) \in E_h$, $n \cdot \dot{\gamma} (s) \neq 0$.
\end{itemize}
\end{definition}

Condition $(b)$ means that for all $h\in S_M$, $\gamma ([0,L])$
may intersect $E_h$ transversally only. It is easy to see 
that condition $(b)$ implies that ${\cal I}(\gamma)=\{ s\in [0,L] \ | \
\gamma (s) \in Q_M \}$ is finite and that there exists $\nu >0$
such that for all $s\in {\cal I}(\gamma)$, for all $h=(n,l)\in S_M$
such that $\gamma (s) \in E_h$, $|\dot{\gamma}(s) \cdot n|/|n| \geq 
\nu$. 

If a curve $ \a $ is not admissible we can always find ``close to it''  
an admissible one $\g$. Indeed the following lemma holds.

\begin{lemma}\label{admissible}
Let $ M>0 $, ${\cal C}$ be a connected component of $ {\cal D}_N^c $, 
$\om_I, \om_F \in {\cal C} $ and let 
$ \a \in C^2 ([0, L_0 ], {\cal C})$
be an embedding with $ \a(0)=\o_I $ and $\a(L_0)=\o_F$.
Then, $\forall\, \eta>0$, there exists a curve $ \g $, 
$Q_M$-admissible between $ \om_I $ and $ \om_F $, satisfying
${\rm dist}(\g(s),\a ([0,L_0]))<\eta, $ $\forall s\in[0,L].$
\end{lemma}

\begin{pf}
First it is easy to see that there exists an embedding $\a_1:
[0,L_1] \to {\cal C}$ such that $\a_1(0)=\o_I, \a_1(L_1)=\o_F$,
${\rm dist}(\a_1(s), \a ([0,L_0])) \leq \eta/4$ and 
$\forall \ h=(n,l) \in S_M$, $\o_{I} \notin E_h$ (resp.
$\o_{F} \notin E_h$)  or
$\dot{\a}_1(0) \cdot n \neq 0$ (resp. $\dot{\a}_1 (L_1) \cdot n \neq 0$).

Let $r>0, \nu_1>0$ be such that $\forall s \in [0,r] \cup [L_1-r,L_1]$,
 $\forall h=(n,l) \in S_M$, ${\rm dist}(\a_1 (s), E_h) \geq \nu_1$
or $|\dot{\a}_1 (s) \cdot n| \geq \nu_1$. Let $\phi: [0,L_1]
\to [0,1]$ be a smooth function such that $\phi (0)=\phi(L_1)=0$
and $\forall s \in [r,L_1-r] \  \phi(s)=1$.

We shall prove that for all $\e>0$ 
there exists $\o_\e \in {\bf R}^d$, 
$|\om_\e|< \e$, such that $\forall h = ( n, l ) \in S_M$,
for all $s\in [r,L_1-r]$ 
such that $\alpha_1(s) \in E_h + \om_\e$, 
$\dot{\a}_1 (s) \cdot n \neq 0$. For $h=(n,l)\in S_M$, let 
${\cal J}_h =\{ s\in [r, L_1-r] \ | \ n \cdot \dot{\a}_1(s)=0 \}$
and ${\cal V}_h=\{ \alpha_1 (s) - u \ | \ s \in {\cal J}_h, \
\ u \in E_h \}$. Let $ \psi_h: [r,L_1-r] \times E_h \to 
{\bf R}^d $ be defined by $ \psi_h (s,u) = \a_1 (s) - u $. 
$ D\psi_h (s,u) $ is singular iff $s \in {\cal J}_h$.
Therefore ${\cal V}_h$ is the set of the critical values of
$ \psi_h $ and by Sard's lemma, meas$ ( {\cal V}_h ) = 0 $. Hence for all
$ \e > 0 $ there exists $ \o_\e \in {\bf R}^d $ such that $ | \o_\e | < \e $,
$\om_\e \notin {\cal V}_h $ for all $h\in S_M$. Our claim follows.

Now we can define $ \alpha_2: [0,L_1] \to {\cal C}$
by $\alpha_2 (s)=\alpha_1(s)-\phi(s) \om_\e$. It is easy
to check that, provided $\e$ is small enough, $\alpha_2$
is an embedding which satisfies condition $(b)$. $\gamma$ is
obtained from $\alpha_2$ by a simple time reparametrization.
\end{pf}

\noindent 
If $ \Gamma ( \a (s), \cdot , \cdot )$ possesses, for each $ s $,
a non-degenerate
local minimum $(\theta_0^{ \a ( s ) }, \f_0^{\a ( s )})$, then, 
by the Implicit Function Theorem, along any curve $ \g $ sufficiently close
to $ \a $, $ \Gamma ( \g (s), \cdot , \cdot )$ possesses
local minima $(\theta_0^{\g(s)}, \f_0^{\g(s)})$ such that
\be \label{unifposdef}
D^2_{(\theta, \f)} \Gamma(\g(s) ,\theta_0^{\g(s)},\f_0^{\g(s)}) >
\lambda {\rm Id}, \qquad \forall  \ s \in [0,L],
\ee 
for some constant $\lambda >0 $ depending on $\a$. 
Therefore, by the above lemma,  it is enough to prove the 
existence of drifting orbits along admissible curves $\g $. Property
$(\ref{unifposdef})$ will be used in lemma \ref{lem:Meln}.

Given a $ Q_M $-admissible curve $\g$, let us call 
$s_1^*, \ldots , s_r^*$ the elements of ${\cal I}(\g)$, and
$\om_1^*=\g(s_1^*), \ldots , \om_r^*=\g(s_r^*)$ the corresponding
frequencies. Since, $\forall  m = 1, \ldots, r $, 
$( \theta_0^{\o_m^*} , \f_0^{\o_m^*})$ is a 
nondegenerate local minimum of $\Gamma(\o_m^*, \cdot, \cdot )$, there
is a neighborhood $W_m$ of $\o_m^*$ such that, 
$\forall \o \in W_m$, $ \Gamma(\o, \cdot )$ 
admits a nondegenerate local minimum $(\theta_0^\om , 
\f_0^\om )$, the map $\om \mapsto (\theta_0^\om , \f_0^\om )$
being Lipschitz-continuous on $W_m$.  
Therefore we shall assume without loss of
generality that for all $m =1, \ldots, r $,
\be \label{localmin}
\forall (\om, \om') \in (W_m\cap \gamma ([0,L]))^2 \  \
 |(\theta_0^{\o} , \f_0^{\o})-(\theta_0^{\o'} , \f_0^{\o'})|
\leq K|\om-\om'|.
\ee

It is easy to prove that,
if $\gamma $ is an admissible curve, there exists $ d_0 > 0 $ such that
\begin{itemize} 
\item[(*)]
$\{ s\in [0,L] \ | \ {\rm dist}(\g(s),Q_M) \leq d_0 \}$
is the union of a finite number of disjoint intervals 
$[S_1,S'_1], \ldots,$ $[S_r,S'_r]$; for all $m = 1, \dots, r $ 
each interval $[S_m,S'_m]$ intersects $ {\cal I}(\g) $ at a unique point
$s_m^*$ and $\g([S_m,S'_m])\subset W_m$. Moreover
$(s\mapsto {\rm dist}(\gamma (s), Q_M))$ is decreasing 
on $[S_m,s_m^*)$, increasing on $(s_m^*,S'_m]$, and 
${\rm dist}(\gamma (s), Q_M) \geq (\nu/2) |s-s_m^*|$
for all $s\in [S_m, S'_m]$.
\end{itemize}
\noindent
Now we are able to define the ``unperturbed transition chain'':
for some small constant 
$ \rho > 0 $ which will be specified later
we choose  $ k\in{\bf N}$ and $ k + 1 $ ``intermediate frequencies'' 
$$ {\om}_I =: \ov{\om}_0, \ov{\om}_1,\dots , 
\ov{\om}_{k-1},\ov{\om}_k := {\om}_F 
$$
with $\ov{\om}_i:=\g(s_i) $ for certain $0=:s_0<s_1<\ldots <s_{k-1}<s_k :=L $ 
verifying
\be\label{isaia}
\frac{\rho\mu}{2}\leq s_{i+1}-s_i \leq \rho \mu, 
\qquad \forall i= 0, \ldots, k - 1.
\ee
By (\ref{isaia}) there results that 
\be\label{kbumps}
\frac{L}{ \rho \mu } \leq k \leq \frac{2L}{\rho \mu}, 
\ee
moreover  it follows from $(a)$ that
\be\label{onejump}
| \ov{\om}_{i+1} - \ov{\om}_i | \leq \rho \mu, 
\qquad \forall i= 0, \ldots, k - 1.
\ee
This condition has been used before in lemma \ref{lem:hjje}. 
Given $ k $ time instants 
$ \ov{\teta}_1:= \teta_0^{\ov{\om}_1} 
< \ov{\teta}_2 < \ldots < \ov{\teta}_i < \ldots < \ov{\teta}_k$, 
we define the $\{ \ov{\vphi}_i \}_{i=1, \ldots, k} $ 
by the iteration formula
\be\label{eq:vphibar}
\ov{\vphi}_1 = \vphi_0^{\ov{\om}_1}, \qquad 
\ov{\vphi}_{i+1} = \ov{\vphi}_i + \ov{\om}_i 
(\ov{\teta}_{i+1} - \ov{\teta}_i ). 
\ee

\indent
The choice of the instants 
$\{ \ov{\teta}_i \}_{i=1, \ldots, k}$ is specified in the next lemma:
the main request is that 
$ (\ov{\teta}_i, \ov{\vphi}_i) $ must arrive $ \d $-close 
${\rm mod} \ 2 \pi {\bf Z}^{d+1} $, to the local 
minimum point $ ( \teta_0^{\ov{\om}_i} , \vphi_0^{\ov{\om}_i})$ 
of the Poincar\'e-Melnikov primitive $ \Gamma (\ov{\om}_i, \cdot, \cdot)$, 
see (\ref{phism})-(\ref{chism}). From (\ref{timerg}) we derive that if 
$ \ov{\om}_i $ is $ 1 \slash | \ln \mu |$ far 
from the set $Q_M$ of ``non-ergodizing frequencies'' we can 
reach this goal for  ``short'' time intervals 
$ \ov{\teta}_{i+1} - \ov{\teta}_i \approx |\ln \mu |$. In order to cross 
the set $Q_M$ of ``non-ergodizing frequencies'' we need to 
use longer time intervals $ \ov{\teta}_{i+1} - \ov{\teta}_i \approx 1 \slash 
{\rm dist} ( Q_M, \ov{\om}_i)$ if
$ \sqrt{\mu}/ | \ln \mu | < {\rm dist} ( Q_M , \ov{\om}_i ) < 
1 \slash | \ln \mu | $.
When the $ \ov{\om}_i $ are ``close'' (less than $ \sqrt{\mu}/ | \ln \mu | $-distant) 
to the set of non-ergodizing hyperplanes $ Q_M $ we choose again
$ \ov{\teta}_{i+1} - \ov{\teta}_i \approx | \ln \mu | $.
We also estimate in (\ref{lengthest}) the total time 
$ \ov{\theta}_k-\ov{\theta}_1 = \sum_{i=1}^k 
\ov{\teta}_{i+1} - \ov{\teta}_i $.

\begin{lemma} \label{choicephi}  
$ \forall \d > 0 $   there exists $\m_6 >0$ 
such that  $ \forall 0 < \m \leq \m_6$ there exist
$ \{ \ov{\teta}_i \}_{i=1,\ldots, k} $ with 
$\ov{\teta}_1 = \teta_0^{\ov{\om}_1} $  satisfying,
\begin{itemize}
\item $(i)$
if \ {\rm dist}$(\ov{\om}_i, Q_{M})> \frac{\sqrt{\mu}}{|\ln \mu |}  $ then 
\be\label{separ1}
\max \Big\{ C_1 |\ln \mu |, \frac{2\pi}{{\rm dist } 
( \ov{\om}_i , Q_{M} )} \Big\} 
< \ov{\teta}_{i+1} - \ov{\teta}_i < 2
\max \Big\{ C_1 |\ln \mu |, \frac{2\pi}{{\rm dist } 
( \ov{\om}_i , Q_M )} \Big\},
\ee
where $M=8\pi a_{d+1}/\delta$;
\item  $(ii)$
if {\rm dist}$(\ov{\om}_i, Q_{M}) \leq \frac{\sqrt{\mu}}{|\ln \mu |} $ then
$C_1 | \ln \mu | < \ov{\teta}_{i+1} - \ov{\teta}_i 
< 2 C_1 |\ln \mu |$,
\end{itemize}
and such that 
\be\label{phism}
{\rm dist}\Big( (  \ov{\teta}_i , \ov{\vphi}_i  ),
(  \teta_0^{\ov{\om}_i},\vphi_0^{\ov{\om}_i}) + 2 \pi {\bf Z}^{d+1} \Big) 
< \d, \quad 
\forall i = 1, \ldots, k,
\ee
where $\ov{\f}_1,\ldots,\ov{\f}_k$ are defined by $(\ref{eq:vphibar})$.
Equivalently, $ \forall i = 1, \ldots, k$, there exist 
 $ h_i \in {\bf Z}^{d+1}$ 
and $ \chi_i \in {\bf R}^{d+1}$ such that
\be\label{chism}
( \ov{\teta}_i , \ov{\vphi}_i,) = 
( \teta_0^{\ov{\om}_i} , \vphi_0^{\ov{\om}_i}) +
2 \pi h_i + \chi_i  \quad {\rm with} \quad  
| \chi_i | < \d.
\ee
Moreover there exists a constant $K(\gamma)$ such that 
\be \label{lengthest}
\ov{\theta}_k-\ov{\theta}_1 \leq K(\gamma)
\frac{|\ln \m|}{\rho \mu}. 
\ee
\end{lemma}

\begin{pf}
Let $\m_6 >0 $ be so small that $\sqrt{\m_6}/| \ln\m_6| < d_0$ and   
$\sqrt{|\ln \m_6|} \geq 32 \sqrt{C_1}/(\nu \sqrt{ \delta \rho}) $.

Let us define $ (\ov{\teta}_1, \ov{\vphi}_1) := 
( \teta_0^{\ov{\om}_1} , \vphi_0^{\ov{\om}_1})$.
Assume that $(\ov{\theta}_1,\ldots, \ov{\theta}_i)$ 
has been defined.
If dist$(\ov{\om}_i, Q_M) > \sqrt{\mu} \slash | \ln \mu |$ 
then by (\ref{timerg})
there certainly exists $(\ov{\teta}_{i+1}, \ov{\vphi}_{i+1})$ 
satisfying $(\ref{eq:vphibar})$,$(\ref{separ1})$, such that
$$
{\rm dist} \Big( ( \ov{\teta}_{i+1}, \ov{\vphi}_{i+1} ), 
(\teta_0^{\ov{\om}_{i+1}} , \vphi_0^{\ov{\om}_{i+1}})+ 
2\pi {\bf Z}^{d+1} \Big) < \d/4.
$$
We now consider the case in which $ \ov{\om}_i $ is close to some 
``non-ergodizing'' hyperplanes of $ Q_M $.
If dist$(\ov{\om}_{i-1},$ $ Q_M) > \sqrt{\mu} \slash |\ln \mu |$ and 
dist$(\ov{\om}_i, Q_M) \leq \sqrt{\mu} \slash |\ln \mu |$
we proceed as follows. 
We have $\ov{\om}_i=\gamma(s_i)$, with $s_i\in [S_q,S'_q]$
for some $q$, $1\leq q \leq r$. Moreover, by property $(*)$ there 
exists 
$ p^* \in {\bf N}$ such that $\{ j \in \{1, \ldots, k\} \ | \ 
s_j \in [S_q,S'_q] \ {\rm and} \ {\rm dist}(\ov{\om}_j,Q_M)
\leq \sqrt{\m}/|\ln \m| \}=\{i,\ldots,i+p^*-1\}$, and
$s_i \leq s_q^* \leq s_{i+p^*-1}$.
We shall use the abbreviations $s^*$ for $s_q^*$, and
$\o^*$ for $\o_q^*$.  We claim that  
\be\label{pstar}
1 \leq p^* \leq p := \Big[ \frac{\sqrt{\d}}{4\sqrt{C_1\rho\mu |\ln\mu|}} 
\Big].
\ee
In fact, by (\ref{isaia}) and  
$(*)$ 
$$
\frac{\nu  \rho}{4} \mu (p^*-1) \leq 
\frac{\nu}{2} [(s_{i+p^*-1}-s^*) + (s_* - s_i)]
\leq {\rm dist}(\ov{\o}_{i+p^*-1},Q_M)+
{\rm dist}(\ov{\o}_{i},Q_M) \leq 
2 \frac{\sqrt{\m}}{|\ln \m|}
$$ 
Hence $p^* \leq 8 (\nu \rho \sqrt{\m} | \ln \m| )^{-1}$, 
which implies $(\ref{pstar})$, by the choice of 
$\m_6$.

Now we can define the $ \ov{\teta}_{i+1},\ldots , \ov{\theta}_{i+p^*}$.
The flow of $(\om^* , 1)$, as any linear flow on a torus,
has the following property :
there exists $T^* (\om^*,\d)>0$ (abbreviated as $T^*$) such that  
any time interval of length $T^*$ contains   $t$ satisfying 
${\rm dist}((t \om^*,t) , 2\pi {\bf Z}^{d+1}) \leq \delta/4$.

Therefore (provided $C_1 |\ln \mu_6| > T^*$) we can define 
$\ov{\theta}_{i+1}, \ldots , \ov{\theta}_{i+p^*}$ such that
\be \label{sepres}
C_1 |\ln \m| \leq \ov{\theta}_{i+j+1}-\ov{\theta}_{i+j} \leq
2 C_1 |\ln \m|, \quad  \quad {\rm dist}\Big(
(\ov{\theta}_{i+j},\wtilde{\f}_{i+j}),
(\ov{\theta}_i, \ov{\f}_i)+2\pi {\bf Z}^{d+1} \Big) \leq \delta/4,
\ee
where $\wtilde{\f}_{i+j}=\ov{\f}_i+\o^*
(\ov{\theta}_{i+j}-\ov{\theta}_i)$.
For $1\leq j\leq p^*$, let
\be\label{vphibarint}
\ov{\vphi}_{i+j} =  \ov{\vphi}_{i} + 
\sum_{q=1}^j \ov{\om}_{i+q-1} (\ov{\theta}_{i+q}-
\ov{\theta}_{i+q-1}).  
\ee
We now check that for all $ j = 1, \ldots, p^* $, 
$(\ov{\teta}_{i+j}, \ov{\vphi}_{i+j})$, as defined in 
(\ref{sepres}) and (\ref{vphibarint}), satisfy  estimate 
(\ref{phism}), namely
\be\label{eq:vicin}
{\rm dist}_T \Big( (\ov{\teta}_{i+j}, \ov{\vphi}_{i+j}),  
( \teta_0^{\ov{\om}_{i+j}}, \vphi_0^{\ov{\om}_{i+j}} ) \Big) :=
{\rm dist} \Big( (\ov{\teta}_{i+j}, \ov{\vphi}_{i+j}),  
( \teta_0^{\ov{\om}_{i+j}}, \vphi_0^{\ov{\om}_{i+j}} )
+2\pi {\bf Z}^{d+1}\Big) \leq \d.
\ee
 We have by (\ref{vphibarint}) that
\begin{eqnarray*} 
{\rm dist}_T \Big( ( \ov{\teta}_{i+j}, \ov{\vphi}_{i+j}),  
(  \ov{\teta}_{i} , \ov{\vphi}_{i}) \Big)
& \leq & {\rm dist}_T \Big( 
(\ov{\teta}_{i+j}, \wtilde{\vphi}_{i+j}),   
(  \ov{\teta}_{i} ,\ov{\vphi}_{i}) \Big)  + \Big|
\sum_{q=1}^j   (\ov{\om}_{i+q-1} - \o^* ) (\ov{\theta}_{i+q}
-\ov{\theta}_{i+q-1})\Big| \\
& \leq & \d/4 + 2 C_1 | \ln \mu |
\sum_{q = 1}^{p^*} |s_{i+q-1}-s^* | \qquad ({\rm by} \ \ 
(\ref{sepres}) \ \ {\rm and} \ \ (a))   \\
& \leq & \d/4 + 2 C_1  | \ln \mu | p^* (s_{i+p^*-1}-s_i) \\
& \leq & \d/4 + 2 C_1  | \ln \m | p^2 \rho \m \leq 3\d /8, 
\end{eqnarray*}
by (\ref{isaia}) and  (\ref{pstar}).
Therefore, by (\ref{localmin}), 
\begin{eqnarray*}
{\rm dist}_T \Big( (\ov{\teta}_{i+j}, \ov{\vphi}_{i+j}),  
( \teta_0^{\ov{\om}_{i+j}}, \vphi_0^{\ov{\om}_{i+j}} )
\Big) & \leq & 
\frac{3\d}{8} + {\rm dist}_T 
\Big( (\ov{\teta}_{i}, \ov{\vphi}_{i}),  
( \teta_0^{\ov{\om}_{i}}, \vphi_0^{\ov{\om}_{i}} )
\Big)+  K |\ov{\o}_{i+j}-\ov{\o}_i| \\
&\leq & \frac{3\d}{8} + \frac{\d}{4}+  K\rho \mu p < \d
\end{eqnarray*}
by $(\ref{pstar})$, provided $\m_6$ has been chosen small enough.

There remains to prove $(\ref{lengthest})$. By $(*)$ we can
write
$$ A_m:=
\Big\{ s \in [S_m,S'_m] \ \Big| \ \frac{\sqrt{\m}}{|\ln \m|} \leq
{\rm dist}(\g (s),Q_M) \leq \frac{1}{2C_1 |\ln \m|} \Big\} = 
[U_m,V_m] \cup [V'_m,U'_m],
$$
with $S_m <U_m<V_m<s^*_m<V'_m<U'_m<S'_m$ (in the case when
$\om^*=\om_{I,F}$, $A_m$ is just an interval). Moreover, by $(a)$, 
$s_m^*-V_m , V'_m-s_m^* \geq \sqrt{\m}/ |\ln \m|$. 
Define $ A := \cup_{m=1}^r A_m$. We have 
$\ov{\theta}_k-\ov{\theta}_1 =\sigma_0+\sum_{m=1}^r \sigma_m$,
where
$$
\sigma_0 :=\sum_{1\leq i \leq k-1, s_i\notin A} (\ov{\theta}_{i+1} -
\ov{\theta}_i), \quad  \quad
\sigma_m := \sum_{1\leq i \leq k-1, s_i \in A_m} (\ov{\theta}_{i+1}
-\ov{\theta}_i).
$$
For $s_i \notin A$, $\ov{\theta}_{i+1}-\ov{\theta}_i \leq 2C_1
|\ln \m |$, hence $\sigma_0 \leq 2C_1 k |\ln \m | \leq 4C_1 L \ln \m /
(\rho \m)$.
For $ i\in A_m$,   $\ov{\theta}_{i+1}-\ov{\theta}_i
\leq $ $4\pi ({\rm dist}(\ov{\o}_i,Q_M))^{-1} \leq$ $ 8\pi/
(\nu |s_i-s_m^*|)$ by $(*)$, and  hence,
using that by (\ref{isaia}) $ s_{i+1} \geq s_i + \rho \mu /2$, 
$$
\sigma_m \leq  \frac{8\pi}{\nu}\sum_{1\leq i \leq k-1,s_i\in A_m} 
\frac{1}{|s_i-s_m^*|} \leq 
\frac{16\pi}{\nu \rho \mu }\sum_{1\leq i \leq k-1,s_i\in A_m} 
\frac{s_{i+1}- s_i}{|s_i-s_m^*|}. 
$$
Estimating the above sum with an integral we easily get
$$
\sigma_m \leq  \frac{8\pi}{\nu(s_m^*-V_m)} + 
\frac{16\pi}{\nu \rho \m} \int_{U_m}^{V_m} 
\frac{ds}{s_m^*-s} +  \frac{8\pi}{\nu(V'_m-s_m^*)}
+ \frac{16\pi}{\nu \rho \m} \int_{V'_m}^{U'_m} 
\frac{ds}{s-s_m^*}.
$$
$(\ref{lengthest})$ can be easily deduced by the bound on 
$s_m^*-V_m , V'_m-s_m^* $.
\end{pf}

In the next section we will prove the existence of a  diffusion 
orbit $( \vphi_\mu, q_\mu )$  close to the 
``unperturbed pseudo-diffusion orbit'' 
$ ( \ov{\vphi}(t), {\ov q}(t)): 
( \ov{\teta}_1, \ov{\teta}_k) \to {\bf R}^{d+1} $ defined, for 
$ t \in [\ov{\teta}_i, \ov{\teta}_{i+1}]$, as 
$ \ov{\vphi} (t) := \ov{\vphi}_i + \ov{\om}_i (t - \ov{\teta}_i) $ 
and $\ov{q}_{|[\theta_i,\theta_{i+1}]} := 
Q_{\ov{\teta}_{i+1} - \ov{\teta}_i}( \cdot - \ov{\teta}_i )$
(mod. $2\pi$). 

\section{The diffusion orbit} \label{sec:difforbit}

We need the following property of the Melnikov function 
$ \wtilde{\Gamma} (\om, \cdot, \cdot )$ defined w.r.t. to
the variables $(b,c)$ by
$$
\wtilde{\Gamma}(\om ,b,c):= \Gamma(\om, \theta_0^\om +b,
\f_0^\om+ b \om +c).
$$

\begin{lemma}\label{lem:Meln}
Assume that $\Gamma ( \om, \cdot, \cdot) $ 
possesses a non-degenerate local minimum in 
$ (\teta_0^\om, \vphi_0^\om)$. Then   
there exist $ r > 0 $, $ \ov{b} > 0 $, 
$ \nu_j  > 0 $ ($ j = 1, 2$) depending only on $\g$ 
such that $ \forall \om = \g ( s ) $, $ s \in [0,L] $ 
\begin{itemize}
\item $(i)$
$\partial_c \wtilde{\Gamma} (\om, b, c) \cdot c \geq \nu_2 > 0$
{\bf or} 
$| \partial_b \wtilde{\Gamma} (\om , b, c) | \geq \nu_1 > 0 $
for $| c | = r, | b | \leq \ov{b} $,
\item $(ii)$
$ \partial_b \wtilde{\Gamma} (\om, b, c) \times {\rm sign} (b) \geq \nu_1 > 0$ 
for $ |c| \leq r $ and $ b = \pm \ov{b} $.
\end{itemize}
\end{lemma}

\begin{pf}
We can assume that $(\ref{unifposdef})$ is satisfied.
Since $ \Gamma (\om, \cdot , \cdot ) $ possesses a non-degenerate 
minimum in $ ( \teta_0^\om, \vphi_0^\om) $, 
 $\wtilde{\Gamma}(\om, b, c)$  possesses in $(0,0)$ a non degenerate 
minimum.
Hence we write $ \wtilde{\Gamma} ( \om , b, c )$, up to a constant, as 
$ \wtilde{\Gamma} (\om , b, c) =  Q_2 (b,c) + Q_3 (b,c)$
where   
$ Q_2 (b,c) =: $ $ \b_\om b^2 / 2 +
 (\a_\om \cdot c) b + ( \g_\om c \cdot c) /2 $
is a  positive definite quadratic form ($\b_\om \in {\bf R}, \a_\om \in {\bf R}^d, 
\g_\om \in {\rm Mat}(d \times d)$) and $Q_3=O(|b|^3+|c|^3)$. 
More precisely, by $(\ref{unifposdef})$, there exists $\e >0$ such
that $\beta_\o >\e $, and $d_\om (c):= 
\b_\om (\g_\om c \cdot c) - (\a_\om \cdot c)^2 > \e |c|^2 $
for all $\om \in \g([0,L])$.
In addition, by the smoothness of $\Gamma$ and the fact that 
$\om=\g(s)$ lives in a compact subset of $ {\bf R}^d $, there
exists a constant $M$ such that, $\forall \om \in \g([0,L])$,
$|\alpha_\o|+|\beta_\o|+ |\gamma_\o| \leq M$, $| \nabla Q_3
(b,c)| \leq M(b^2+|c|^2)$.

We have 
$\partial_b Q_2 (b, c) =$ $\b_\om b + \a_\om \cdot c $ and
$\partial_c Q_2 (b, c) \cdot c = $ $b \a_\om \cdot c + (\g_\om c \cdot c)$.

\noindent
Let us define
$ \ov{\nu}_1 := \inf_{\om \in \g ([0,L]) } \e \slash (4 |\a_\om |) > 0 $
and 
$\ov{\nu}_2 := \inf_{\om \in \g ([0,L])} \e \slash (4 \b_\om) > 0 $.
Then consider  $ \nu_1 := \ov{\nu}_1 r $, $\nu_2=\ov{\nu}_2 r^2$
and $ \ov{b} := 
r \sup_{\om \in {\g}([0,L]) } (3 \ov{\nu}_1 + |\a_\om |) \slash \b_\om
$, $r\in (0,1]$. 
We now prove that, provided $ r > 0 $ has been chosen sufficiently small, 
conditions $(i)$ and $(ii)$ are satisfied with the  
above choice of the constants.
Indeed if 
$(|\a_\om \cdot c| + 2 \ov{\nu}_1 r) \slash \b_\om \leq |b| \leq \ov{b} $
and $ | c | \leq r$ then  
$  \partial_b \wtilde{\Gamma} (\om, b, c) \cdot {\rm sign}(b) \geq $ 
$ \b_\om |b| - |\a_\om \cdot c| - |\partial_b Q_3 (b,c)| \geq $ 
$ 2 \ov{\nu}_1 r - O(r^2 ) \geq \nu_1 $ for $r$ sufficiently small.
In particular this proves $(ii)$.
On the other hand if 
$ |b| < (|\a_\om \cdot c | + 2\ov{\nu}_1 r) \slash \b_\om $ and $ | c | = r $ 
then 
\begin{eqnarray*}
\partial_c \wtilde{\Gamma}(\om, b, c) \cdot c   
& = & b (\a_\om \cdot c) + (\g_\om c \cdot c) + \partial_c Q_3 (b, c) \cdot c  
\geq (\g_\om c \cdot c) - | b (\a_\om \cdot c) | + O(r^3) \\
& \geq &  \frac{\e r^2 + (\alpha_\o \cdot c)^2- |\alpha_\o \cdot
c|(|\a_\om \cdot c | + 2 \ov{\nu}_1 r)}{\b_\om}
+ O(r^3) \\
& \geq & \frac{\e  - 2\ov{\nu}_1 |\a_\om | }{\b_\om} r^2 + O(r^3) 
\geq \frac{\e }{ 2\b_\om} r^2 - O(r^3 )\geq 2\ov{\nu}_2 r^2+O(r^3).
\end{eqnarray*} 
Hence $(i)$ is satisfied for $r $ small enough.
\end{pf}

The partial derivatives of $\wtilde{\Gamma}$ are
Lipschitz-continuous  w.r.t. $(b,c)$ uniformly in
$\o \in \g([0,L])$. Therefore, by lemma \ref{lem:Meln},
there exists $\delta>0$ such that, $\forall \eta \in {\bf R}$
with  $|\eta| \leq \d$, $\forall \xi \in {\bf R}^d $
with $|\xi| \leq \d$, $\forall \o \in \g([0,L])$,
\be \label{defdel1}
\partial_c \wtilde{\Gamma} (\om, b+\eta, c+\xi) \cdot c \geq 3\nu_2/4 >
0  \quad  
{\bf or}  \quad
| \partial_b \wtilde{\Gamma} (\om , b+\eta, c+\xi) | \geq 3\nu_1/4 > 0
\quad 
{\rm for} \quad | c | = r, | b | \leq \ov{b} ,
\ee 
\be \label{defdel2}
 \partial_b \wtilde{\Gamma} (\om, b+\eta, c+\xi) \times {\rm sign} (b) \geq 3\nu_1/4 > 0 
\quad {\rm for} \quad   |c| \leq r \ \ {\rm and} \ \  b = \pm \ov{b} .
\ee
Moreover let us fix $\rho>0$ such that 
\be\label{defrho}
\rho \leq \min \{ \nu_1 / 2 , \nu_2 / r \} \slash ( 6 C_2 ),
\ee
where $C_2$ appears in $(\ref{smalldef})$.
These are the positive constants $(\d, \rho)$ 
that we use in order to define, for  $ 0 < \mu < \mu_6 $, 
$\ \ov{\o}_i$, $\ov{\theta}_i$, $\ov{\f}_i$ 
by lemma \ref{choicephi}.

Since $ \g ([0, L]) $ is a compact subset of $ {\cal D}_N^c $, 
$ \inf_{s \in [0,L] }\beta ( \g (s) ) > 0 $ and, 
by the choice of $ \ov{\theta}_i $, for $ \mu $
small enough $ ( \ref{betai} ) $ is satisfied. Therefore,
by lemma \ref{approxsum} and (\ref{chism}), there exists $ \mu_7 > 0 $ such that,
$\forall 0 < \m \leq \m_7 $, 
\be \label{Fmubc} 
{\cal F}_\mu ( b, c)  =   \frac{1}{2} \sum_{i=1}^{k-1} 
\frac{|c_{i+1} - c_i|^2}{\Delta \ov{\teta}_i + (b_{i+1} - b_i)} + 
\mu \sum_{i=1}^k \wtilde{\Gamma} (\ov{\om}_i, {\eta}_i+ b_i, 
{\xi}_i + c_i )+R_7,
\ee
where $| \eta_i | \leq \d $, $ | \xi_i | \leq \d $, $R_7$ 
is given by $(\ref{R7})$ and satisfies
$(\ref{smalldef})$.
\\[1mm]
\indent
We minimize the functional ${\cal F}_\mu $ on the closure of
$$
W := \Big\{ (b, c):= (b_1, c_1, \ldots, b_k, c_k)  
\in {\bf R}^{(d+1)k} \ \Big| \ | b_i | < \ov{b}, \ |c_i | < r,  \
\forall i = 1, \ldots, k \Big\}. 
$$
Since $\ov{W}$ is compact, ${\cal F}_\mu$ attains
its minimum in $ {\ov W} $, say at $( \wtilde{b}, \wtilde{c} )$.  
By lemma \ref{lem:heter}
the existence of the diffusion orbit will be proved 
once we show that $(\wtilde{b}, \wtilde{c} ) \in W $, see lemma 
\ref{inter}. Let us define for $ i = 1, \ldots, k - 1 $  
$$
w_i := w_i(b,c):=\frac{c_{i+1} - c_i }{ \teta_{i+1} - \teta_i  } =
\frac{c_{i+1} - c_i }{ \Delta \ov{\teta}_i + (b_{i+1} -b_i)},
$$ 
and $w_0 = w_k = 0 $.
From (\ref{eq:vphibar}) and (\ref{coordc}), 
$w_i$ can be written as 
\be\label{smallwi}
w_i = \frac{ \vphi_{i+1} - \vphi_i }{(\teta_{i+1} - \teta_i )} 
- \ov{\om}_i  - \frac{\Delta \ov{\om}_i b_{i+1} }{(\teta_{i+1} - 
\teta_i )} =
\Big( \om_i - \ov{\om}_i \Big) + O \Big( \frac{\mu }{| \ln \mu |} \Big).
\ee
By the expression of ${\cal F}_\mu $ in $(\ref{Fmubc})$ we have,
for all $i = 1, \ldots, k $,
\be\label{partialci}
\partial_{c_i} {\cal F}_\mu (b,c)  
= w_{i-1} - w_i + \mu \partial_c
\wtilde{\Gamma} (\ov{\om}_i, \eta_i + b_i, \xi_i + c_i ) + R_i
\ee
\be\label{partialb}
\partial_{b_i} {\cal F}_\mu (b, c)=
\frac{1}{2} \Big( |w_i |^2 - |w_{i-1}|^2 \Big) + \mu \partial_b 
\wtilde{\Gamma} ( \ov{\om}_i, \eta_i + b_i, \xi_i + c_i )
+ S_i
\ee
where 
$ R_i:= \partial_{c_i} R_7 $, $ S_i :=  \partial_{b_i} R_7 $
satisfy, by (\ref{smalldef}) and $(\ref{defrho})$
\be\label{restpiccolo}
| R_i |, |S_i| \leq \frac{\mu}{2} \min \Big\{ \frac{\nu_1}{2}, \frac{\nu_2}{r} \Big\}.
\ee
By (\ref{partialci})-(\ref{partialb}), a way to see critical points 
of ${\cal F}_\mu $ is to show that the terms $ w_{i-1} - w_i $ and 
$ |w_i |^2 - |w_{i-1}|^2 $ are small w.r.t the 
$O( \mu )$-contribution provided by the Melnikov function. 
By (\ref{sepfre}) $ |\om_i - \ov{\om}_i | = 
O (1 \slash ( \teta_{i+1} - \teta_i ) )$ and hence,
using  (\ref{smallwi}), an estimate for each $ w_i $ separately
is given by
$ w_i = O( 1 \slash | \ov{\teta}_{i+1} - \ov{\teta}_i | ) + 
O(\mu / | \ln \mu | ) $. 
Hence each $ | w_i | $ is $O(\mu)$-small if the
time  to make a transition 
$| \ov{\teta}_{i+1} - \ov{\teta}_i | = O( 1 \slash \mu )$,
as in \cite{Bs}. These time intervals are too large to 
obtain the approximation for the reduced action functional 
${\cal F}_\mu $ given in lemma \ref{approxsum} and $(\ref{Fmubc})$.
Therefore we need more refined estimates: 
the proof of Theorem \ref{thm:main} (and Theorem \ref{thm:Arn})
relies on the following crucial property
for $ \wtilde{w}_i:=w_i(\wtilde{b}, \wtilde{c}) $, 
satisfied by the minimum point $(\wtilde{b}, \wtilde{c})$.
\begin{lemma} We have
(for $ i = 1, \ldots, k,$) 
\be\label{wnearmu}
i) \ \ |\wtilde{w}_i - \wtilde{w}_{i-1}| = O ( \mu ), \quad  \quad
ii) \ \ | \wtilde{w}_i | = O \Big( \frac{\sqrt{\mu}}{\sqrt{{|\ln \mu
|}}} \Big).
\ee
\end{lemma}

\begin{pf} 
Estimate $(\ref{wnearmu})-i)$  is a straightforward consequence of 
(\ref{partialci}) and (\ref{restpiccolo}) 
if  $ | \wtilde{c}_i | < r $, since in this 
case $ \partial_{c_i} {\cal F}_\mu (\wtilde{b}, \wtilde{c}) = 0 $. 
We now prove that $(\ref{wnearmu}) i)$
holds also if $ | \wtilde{c}_i | = r $ for some $ i $.
Indeed if $ | \wtilde{c}_i | = r $ then
\be\label{alphal}
\partial_{c_i} {\cal F}_\mu (\wtilde{b}, \wtilde{c})
= \alpha_\mu \wtilde{c}_i \quad
{\rm for \ some} \quad \alpha_\mu \leq 0
\ee
(since $(\wtilde{b}, \wtilde{c})$ is a minimum point)
and then by (\ref{partialci}),(\ref{alphal}) and (\ref{restpiccolo})
we deduce 
\be\label{estidiff}
\wtilde{w}_{i-1} - \wtilde{w}_i = \alpha_\mu \wtilde{c}_i + O( \mu ). 
\ee
Let us decompose $ \wtilde{w}_{i-1}$ and $ \wtilde{w}_i $ in the ``radial'' and 
``tangent'' directions to the ball 
$ S_i = \{  |b_i | \leq \ov{b},  \ |c_i | \leq r \}$: 
\be\label{dec1}
\wtilde{w}_{i-1} = a_i \wtilde{c}_i 
+ u_i \qquad {\rm with } \ u_i \cdot \wtilde{c}_i =0
\ee
\be\label{dec2}
- \wtilde{w}_i = a'_i \wtilde{c}_i + u'_i, 
\qquad {\rm with } \ u_i' \cdot \wtilde{c}_i = 0.
\ee 
Since $| \wtilde{c}_{i-1}| \leq |\wtilde{c}_i| =r$,
$|\wtilde{c}_{i+1}| \leq |\wtilde{c}_i| = r $, there results that  
\be\label{positive}
a_i r^2 = \wtilde{w}_{i-1} \cdot \wtilde{c}_i \geq 0 \quad {\rm and } \quad
a'_i r^2 = - \wtilde{w}_i \cdot \wtilde{c}_i \geq 0,
\ee
so that  $ a_i, a'_i \geq 0 $.
Summing (\ref{dec1}) and (\ref{dec2}) and using (\ref{estidiff}) we obtain 
$$
( a_i + a'_i ) \wtilde{c}_i + (u_i + u'_i) = O(\mu)+ \alpha_\mu \wtilde{c}_i,
$$
with $ a_i, a'_i, -\alpha_\mu \geq 0$. This implies 
that  $ \a_\mu = O(\mu / r )$ and from equation 
(\ref{estidiff}) we get $(\ref{wnearmu}) i)$. 

We can now prove $(\ref{wnearmu})-ii)$. 
Let 
$ i_0 \in \{1, \ldots, k - 1\} $ be such that $ \forall 1\leq i
\leq k-1 $,  
$ |\wtilde{w}_{i_0} | \geq |\wtilde{w}_i | $.
For $ j \in \{ 1, \ldots, k -1 \}$, $ j \neq i_0 $ we can write
$ \wtilde{w}_j = \wtilde{w}_{i_0} + s_j $ with  
$ s_j = \sum_{i = i_0}^{j-1} ( \wtilde{w}_{i+1} - \wtilde{w}_i )$
and hence, by $(\ref{wnearmu}) i)$
\be\label{oralo} 
| s_j | \leq \sum_{i = i_0}^{j-1} 
| \wtilde{w}_{i+1} - \wtilde{w}_i|
\leq C \mu |j - i_0 |
\ee
for some constant $C > 0 $.  Hence
\be
\wtilde{c}_j - \wtilde{c}_{i_0} =  
\sum_{i=i_0}^{j-1} \wtilde{w}_i ( \wtilde{\teta}_{i+1} - \wtilde{\teta}_i )
=
\wtilde{w}_{i_0} (\wtilde{\teta}_{j} - \wtilde{\teta}_{i_0} ) + 
\sum_{i=i_0}^{j-1} s_i ( \wtilde{\teta}_{i+1} - \wtilde{\teta}_i )
\ee
and then by (\ref{oralo})
\be\label{quasifin}
\Big| \wtilde{c}_j - \wtilde{c}_{i_0} \Big| \geq
 |\wtilde{w}_{i_0}||\wtilde{\teta}_j - \wtilde{\teta}_{i_0}| -
C \mu |j - i_0 || \wtilde{\teta}_j - \wtilde{\teta}_{i_0}| =
\Big(  |\wtilde{w}_{i_0}| - C \mu |j - i_0 | \Big) 
| \wtilde{\teta}_j - \wtilde{\teta}_{i_0}|.
\ee
Since $| \wtilde{\teta}_{i+1} - \wtilde{\teta}_i| > C_1 | \ln \mu
|+O(1)$ (by (\ref{ens})),
$\forall i=1, \ldots, k-1$,  
$ | \wtilde{\teta}_j - \wtilde{\teta}_{i_0}| > C_1 |j -i_0| 
\cdot | \ln \mu | $. 
Take $ \ov{j} \in \{ 1, \ldots, k -1 \}$ such that
$ | \ov{j} - i_0 | = 
[(\sqrt{\mu} \sqrt{| \ln \mu |})^{-1}] +1 $
(such a $\ov{j}$ certainly exists since, by (\ref{kbumps}),
$k \approx 1 \slash \mu $ for $ \mu $ small).
Then  we obtain, using that  
$ |\wtilde{c}_i | \leq r $ for all $ i = 1, \ldots, k $,
$$
2r \geq \Big| \wtilde{c}_j - \wtilde{c}_{i_0} \Big| \geq
\Big(  |\wtilde{w}_{i_0}| - C \frac{\sqrt{\m}}{\sqrt{|\ln \m|}} -C\mu \Big) 
C_1 \frac{\sqrt{|\ln \mu|}}{\sqrt{\mu}},
$$
{\it i.e.} $|\wtilde{w}_{i_0}| \leq \dps 
\frac{(2r+CC_1)\sqrt{\mu}}{C_1
\sqrt{|\ln \m|}}+C\mu$.
 We have thus 
proved the important property $(\ref{wnearmu})-ii)$.
\end{pf}

\begin{remark}\label{rem:aprio}
By (\ref{smallwi}), 
$ (\wtilde{\om}_i - \ov{\om}_i ) = \wtilde{w}_i + 
O(\mu / |\ln \mu |) $, so that, by  (\ref{onejump}),  
(\ref{wnearmu}) implies
\be\label{apiroriest}
| \wtilde{\om}_i - \ov{\om}_i | 
=O \Big( \frac{\sqrt{\mu}}{ \sqrt{|\ln \mu|}} \Big), \qquad
| \wtilde{\om}_{i+1} - \wtilde{\om}_i | = O(\mu).
\ee
Note that, from (\ref{sepfre}),  we would just obtain 
$| \wtilde{\om}_i - \ov{\om}_i | =O( 1 / | \ln \mu |)$.
(\ref{apiroriest}) can be seen as an a-priori estimate satisfied
by the minimum point $(\wtilde{\teta}, \wtilde{\vphi}) $. 
\end{remark}

The following lemma proves the existence of a local minimum 
of the reduced action functional in the interior of $ W $ 
and hence of a true diffusion orbit.

\begin{lemma}\label{inter}
Let  $ ( \wtilde{b}, \wtilde{c})$ be a
minimum point of $ {\cal F}_\mu $ over $\ov{W} $. Then
$ ( \wtilde{b}, \wtilde{c}) \in W $, namely
\be\label{cir}
| \wtilde{c}_i | < r \qquad {\rm for \ all}  
\qquad i \in \{ 1, \ldots, k \}
\ee
and
\be\label{bint}
| \wtilde{b}_i | < \ov{b} \qquad {\rm for \ all} 
\qquad i \in \{1, \ldots, k \}. 
\ee
\end{lemma}

\begin{pf}
By (\ref{wnearmu}) we have
$ | |\wtilde{w}_{i+1}|^2 - |\wtilde{w}_i|^2 | \leq$ 
$ | \wtilde{w}_{i+1} - \wtilde{w}_i |  \cdot ( 
| \wtilde{w}_{i+1} |+ |\wtilde{w}_i |)  = O(\mu^{3/2}) $,
and hence, from (\ref{partialb}) we derive
\be\label{piccolb}
\partial_{b_i} {\cal F}_\mu (\wtilde{b}, \wtilde{c})= 
\mu \partial_b 
\wtilde{\Gamma} ( \ov{\om}_i, \eta_i + \wtilde{b}_i, \xi_i 
+ \wtilde{c}_i ) + O ( \mu^{3/2} ) + S_i.
\ee
Let us first 
assume  by contradiction that $\exists i$ such that 
$ |\wtilde{c}_i | = r $ and $|\wtilde{b}_i | < \ov{b} $. In this case 
we claim that
\be\label{minimump}
\partial_c \wtilde{\Gamma} ( \ov{\om}_i,  \eta_i + \wtilde{b}_i, 
\xi_i + \wtilde{c}_i) \cdot \wtilde{c}_i \leq \nu_2 \slash 2 \quad
{\rm {\bf and} } \quad 
| \partial_b \wtilde{\Gamma} 
( \ov{\om}_i, \eta_i + \wtilde{b}_i, \xi_i + \wtilde{c}_i ) | 
\leq \nu_1 \slash 2
\ee
contradicting $(\ref{defdel1})$, 
since $|\eta_i |, | \xi_i | \leq \d $. Let us prove (\ref{minimump}).
Since  $(\wtilde{b}, \wtilde{c})$ is a minimum point 
$$
\partial_{c_i} 
{\cal F}_\mu (\wtilde{b},\wtilde{c}) \cdot \wtilde{c}_i = 
 (\wtilde{w}_{i-1} - \wtilde{w}_i) \cdot \wtilde{c}_i + \mu \partial_c
\wtilde{\Gamma} (\ov{\om}_i, \eta_i + \wtilde{b}_i, 
\xi_i + \wtilde{c}_i ) \cdot \wtilde{c}_i + R_i \cdot \wtilde{c}_i
= \a_\mu \wtilde{c}_i \cdot \wtilde{c}_i =  \a_\mu r^2 \leq 0.
$$
By (\ref{positive}) and (\ref{restpiccolo}) it follows that 
$ \partial_c \wtilde{\Gamma} ( \ov{\om}_i, \eta_i + \wtilde{b}_i, 
\xi_i + \wtilde{c}_i) \cdot \wtilde{c}_i \leq \nu_2 \slash 2 $.
Moreover since $| \wtilde{b}_i | < \ov{b} $ we have
$\partial_{b_i} {\cal F}_\mu (\wtilde{b}, \wtilde{c})=0$, and 
 by (\ref{piccolb}),  (\ref{restpiccolo}) it follows that 
$| \partial_b \wtilde{\Gamma} ( \ov{\om}_i, \eta_i 
+ \wtilde{b}_i, \xi_i + \wtilde{c}_i ) | \leq \nu_1 \slash 2 $
(provided $\m$ is small enough).   
Estimate (\ref{minimump}) is then proved. As a result, if
$(\ref{bint})$ holds, so does $(\ref{cir})$.

Let us finally prove  (\ref{bint}). 
If by contradiction $\exists i$ with  $|\wtilde{b}_i | = \ov{b} $,  
by (\ref{piccolb}), (\ref{restpiccolo}) and since
$(\wtilde{b}, \wtilde{c})$ is a minimum point, 
arguing as before, we deduce that 
$ \partial_b \wtilde{\Gamma} 
( \ov{\om}_i, \eta_i +\wtilde{b}_i,  \xi_i + \wtilde{c}_i) 
{\rm sign}(\wtilde{b}_i) \leq \nu_1 /2 $. 
This contradicts (\ref{defdel2}) since 
$ |\eta_i |, | \xi_i | \leq \d $.  The lemma is proved.
\end{pf}

\begin{pfn}{\sc of Theorem } \ref{thm:main}.
Lemmas \ref{inter} and \ref{lem:heter} imply the 
existence of a diffusion orbit 
$ z_\mu (t) := ( \vphi_\mu (t), q_\mu (t), I_\mu (t), p_\mu (t))$ with 
${\dot \vphi_\mu} (\wtilde{\teta}_1 ) = \om_I + O(\mu) $ and 
${\dot \vphi_\mu} (\wtilde{\teta}_k ) = \om_I + O(\mu) $
($z_\mu (\cdot ) $ connects a $ O ( \mu ) $-neighborhood 
of ${\cal T}_{\om_I}$ to a 
$ O(\mu) $-neighborhood of $ { \cal T }_{ \om_F }$ in the time-interval 
$ (\tau_1, \tau_2 )$ where
$ \tau_1 := ( \wtilde{\teta}_1 + \wtilde{\teta}_2 ) \slash 2 $,
$ \tau_2 := (\wtilde{\teta}_{k-1} + \wtilde{\teta}_k ) \slash 2 $). 
The estimate on the diffusion time is a straightforward
consequence of (\ref{lengthest}) and the fact that
$\wtilde{\theta}_{1,k}=\ov{\theta}_{1,k}+O(1)$.
 That ${\rm dist}(I_\mu (t) , \gamma ([0,L]))<\eta$ for
all $t$, provided $\mu$ is small enough, results from 
$(\ref{apiroriest})$ and the estimates of lemma 
\ref{lem:connecting}. 

Finally we observe that, if the perturbation is 
$\mu( f + \mu \wtilde{f} )$, then lemma \ref{lem:connecting} still applies
with the same estimates. Moreover in the development of the 
reduced functional the term containing $ \mu^2 \wtilde{f} $  
gives, in time intervals 
$ \ov{\teta}_{i+1} - \ov{\teta}_i \leq const. |\ln \mu | / \sqrt{\mu} $, 
negligible contributions $ o (\mu) $. Therefore the same variational
proof applies.
\end{pfn}

\begin{pfn}{\sc of Theorem } \ref{thm:Arn}.
If the perturbation is of the form
$ f(\vphi, q, t ) = (1 - \cos q) f( \vphi, t) $, 
by remark \ref{remtranl}-$2)$, we  can prove that the
development (\ref{eq:somma}) holds along any path $ \g $ of
the action space  (without any condition as (\ref{betai})). 
Therefore the previous variational argument applies.
\end{pfn}

For $ \b > 0 $ small let 
${\cal D}_N^\b$ be the set of frequencies ``$\b$-non-resonant
with the perturbation'' 
$ {\cal D}_N^\b := \{ \om \in {\bf R}^d \ | 
\ |\om \cdot n + l| > \b, \ \ \forall \ 0 < |(n,l )| \leq N \}.$
If $ \b $ becomes small with $ \mu $ our estimate
on the diffusion time required to 
approach to the boundaries of ${\cal C} \cap {\cal D}_N^\b $ 
slightly deteriorates. 
In the same hypotheses as in  Theorem \ref{thm:main}
we have the following result.

\begin{theorem}\label{thm:tempobordo}
$\forall R > 0 $, $\forall \ 0 \leq a < 1 / 4$,
there exists $ \mu_8 > 0 $ such that $ \forall 0 < \mu \leq \mu_8 $, 
$ \forall \om_I, \om_F \in {\cal C} \cap 
{\cal D}_N^{\mu^a} \cap B_R ( 0 ) $  there exist 
a diffusion orbit
$( \vphi_\mu (t), q_\mu (t),$ $ I_\mu (t) , p_\mu (t))$ of $({\cal S}_\mu )$ 
and two instants 
$ {\tau}_1 < {\tau}_2 $ with $ I_\mu ( \tau_1 ) = \om_I + O(\mu) $,  
$I_\mu ( \tau_2 ) = \om_F + O(\mu) $ and  
\be\label{totaltime2}
|{\tau}_2 - {\tau}_1| = O ( 1 / \mu^{1+a} ).
\ee
\end{theorem} 

\begin{pf}
For simplicity we consider the case in which 
$ \beta(\om_I) = O(\mu^a) $ and $ \beta(\om_F)=O(1) $.
With respect to Theorem \ref{thm:main} we only need to prove 
the existence of a diffusion orbit connecting $\om_I$ to some fixed
$ \om^* $ lying in the same connected component of 
$ { \cal D }_N^c \cap B_R (0) $ containing $ \om_I $. 
In order to construct an orbit connecting $ \om_I $ to $ \om^* $ 
we can define $ \ov{\om}_i := \om_I + i (\om^* - \om_I )/ k $,
for $0 \leq i \leq k $ and $ k := [ |\om^* - \om_I |/ \rho \mu ] +1 $.
We obtain that $ \beta_j = \beta(\ov{\om}_j) \geq
C(\mu^a+j\rho\mu)$ for some $C>0$ and 
we choose $ \ov{\teta}_{j+1} - \ov{\teta}_j \geq const. \beta_j^{-2} $
verifying in this way the hypotheses of lemma \ref{approxsum}. 
If $ \om_I $ belongs to some $Q_M$ the transition  times
$|\ln \mu|/\sqrt{\mu}$ needed to cross $ Q_M $ (see lemma \ref{choicephi})
still satisfy (\ref{betai}). We finally obtain a diffusion time
$\ov{\teta}_{k} - \ov{\teta}_1 = 
\sum_{j=1}^{k-1}(\ov{\teta}_{j+1}- \ov{\teta}_j) = 
O(1 \slash \mu^{1+a}).$
\end{pf}

\section{The stability result and the optimal time}\label{sec:opt}

In this section we will prove, via classical perturbation theory, 
stability results for the action variables, implying, in particular,
Theorem \ref{thm:Nec}. 
We shall use the following notations: for  $l \in {\bf N} $, 
$ A \subset {\bf C}^l $ and $ r > 0 $,  
we define $ A_r := \{ z \in {\bf C}^l \ | \ {\rm dist} ( z, A ) 
\leq r \}$ and  $ {\bf T}^l_s :=
\{ z \in {\bf C}^l \  | \ | {\rm Im } \ z_j | < s, \ 
\forall\  1 \leq j \leq l \}$
(thought of as a complex neighborhood of ${\bf T}^l$).
Given two bounded open sets $B \subset {\bf C}^2,$  $ D\subset {\bf C}^l$ and 
 $ f ( I,\f,p,q ) $, real analytic function with holomorphic extension on 
$ D_{\sigma} \times {\bf T}^l_{s+\sigma} \times B_{\sigma}$ for
some $\sigma>0$, we define the following norm 
$\| f \|_{B,D,s}=\sum_{k\in  {\ind Z}^l} \sup_{(p,q)\in B \atop I\in D} 
| \hat f_k(I,p,q)| e^{|k|s}\,$ where
$ \hat f_k(I,p,q) $ denotes the $ k $-Fourier 
coefficient of the periodic function $ \f \to f(I,\f,p,q)$.
\\[1mm]
\indent
Let us  consider Hamiltonian $ {\cal H}_\mu $ defined in (\ref{eq:Hamg}) 
and assume that $ f ( I, \vphi, p, q, t) $, 
defined in (\ref{trigpert}), is a real analytic function, 
possessing, for some $ r, \ov{r}, \wtilde{r}, s > 0 $,
complex analytic extention on 
$\{ I\in {\bf R}^d \ | \  
| I | \leq \ov{r} \}_{r}\times {\bf T}^d_s\times
\{ p\in {\bf R}\ | \  | p|\leq \wtilde{r}\}_{r}\times {\bf T}_s 
\times {\bf T}_s$.

It is convenient to write Hamiltonian ${\cal H}_\mu $ in autonomous form.
For this purpose 
let us introduce the new action-angle variables 
$( I_0 , \f_0 ) $ with $ t = \f_0 $, that will still be denoted by
$ I:= ( I_0, I_1, \ldots, I_n)$ and $ \f:=(\f_0,\f_1, \ldots,\f_n)$.
Defining $ h (I) := I_0 + |I|^2/2 $ and $ E:=E(p,q):=p^2/2+(\cos q-1),$
${\cal H}_\mu $ is then equivalent to the autonomous Hamiltonian
\begin{equation}\label{hamiltoniana}
H := H( I, \f, p, q) := h(I) + E( p, q ) + \mu f( I, \f, p, q ).
\end{equation}
Clearly, Hamiltonian $ H $ 
is a real analytic function, with complex analytic extention on 
$$ 
\Big\{ I\in {\bf R}^{d+1}\ \Big| \  
| I| \leq \ov{r} \Big\}_{r} \times {\bf T}^{d+1}_s \times
\Big\{ p\in {\bf R}\ \Big| \  | p| \leq \wtilde{r} \Big\}_{r} 
\times {\bf T}_s.
$$
In the sequel we will denote by 
$z(t):=$ $(I(t),\f(t),p(t),q(t))$ the solution of the Hamilton 
equations associated to Hamiltonian (\ref{hamiltoniana})
with initial condition $z(0)=$ $(I(0),\f(0),p(0),q(0)).$
\\[1mm]
\indent
The proof of the stability of the action variables is divided in two steps:
\begin{itemize}
\item (i) ({\bf Stability far from the separatrices of the pendulum:})
prove stability in the region 
$$
{\cal E}_1 := {\cal E}_1^+\cup {\cal E}_1^-:=
\Big\{(I,\f,p,q) \ | \  E(p,q) \geq \mu^{c_d}\ \Big\} \cup 
\Big\{ (I,\f,p,q) \ | \ -2+\mu^{c_d}\leq E ( p, q ) \leq -\mu^{c_d}  \Big\}
$$
in which we can apply the Nekhoroshev Theorem obtaining actually stability
for exponentially long times,
\item (ii) ({\bf Stability close to the separatrices 
of the pendulum and to the elliptic equilibrium point:}) 
prove stability in the region    
$$
{\cal E}_2 := {\cal E}_2^+\cup {\cal E}_2^-:=
 \Big\{(I,\vphi, p, q)\ | \  - 2 \mu^{c_d} \leq E (p,q) 
\leq 2 \mu^{c_d}  \Big\} \cup 
\Big\{ (I,\vphi, p, q) \ | \  -2 \leq E (p, q) \leq -2 + 2 \mu^{c_d}  \Big\}
$$
in which we use some {\sl ad hoc} arguments,
\end{itemize}
where $ 0 < c_d < 1 $ is a positive constant
that will be chosen later on, see (\ref{eq:cd}).
\\[2mm]
\indent
We first prove (i). 
In the regions\footnote{$\Pi_{p,q}$ denotes the projection onto the 
$(p,q)$ variables.} ${\wtilde{\cal E}}_1^{\pm} := 
\Pi_{q,p} {\cal E}_1^{\pm} $ we first write 
the pendulum Hamiltonian $ E(p, q)$ in action-angle variables.
In the region\footnote{The case with $p<0$ is completely analogous.} 
${\wtilde{\cal E}}_1^+\cup\{ p>0 \}$
the new action variable $ P $ is defined by the formula
$$
P:=P^+(E):=\frac{\sqrt{2}}{\pi}\int^{\pi}_0\sqrt{E+ (1+\cos\psi)}\, d\psi.
$$
while in the region  ${\wtilde{\cal E}}_1^-$ the new action variable is
$$
P:=P^-(E)=\frac{2\sqrt{2}}{\pi}\int^{\psi_0(E)}_0\sqrt{E+(1+\cos\psi)}\, d\psi
$$
where  $\psi_0(E)$ is the first positive number such that
$E+(1+\cos\psi_0(E))=0$.
We will use the following lemma, proved in \cite{BC}, regarding 
the analyticity radii of these action-angle variables close to the
separatrices of the pendulum. 
\begin{lemma}
There exist intervals $D^{\pm}\subset {\bf R},$ symplectic transformations
$\phi^{\pm}=\phi^{\pm}(P,Q)$ real analytic on  $D^{\pm}\times{\bf T}$
with holomorphic extension on 
$D^{\pm}_{r_0}\times{\bf T}_{s_0}$ and  functions $E^{\pm}$
real analytic 
on $D^{\pm}$ with holomorphic extension on   $D^{\pm}_{r_0}$ such that 
$\phi^{\pm}(D^{\pm}\times{\bf T} )={\wtilde{\cal E}}_1^{\pm}$ and 
$$
E(\phi^{\pm}(P,Q))=E^{\pm}(P),
$$
with $ r_0 = const  \mu^{c_d}$ and $s_0 =const /|\ln \mu |$.
Moreover, for $ E $ bounded, the following estimates on the derivatives 
hold \footnote{If $f(x),g(x)$ are positive function,
with the symbol
$f\approx g$ we mean that $\exists\, c_1,c_2>0$ such that 
$c_1g(x)\leq f(x)\leq c_2g(x),\, \forall\, x.$ }
\begin{eqnarray}
\frac{dE^{\pm}}{dP}(P^{\pm}(E)) &\approx& 
\ln^{-1}(1+\frac{1}{\sqrt{|E|}})\label{E1}\\
\pm\frac{d^2E^{\pm}}{dP^2}(P^{\pm}(E)) 
&\approx& \frac{1}{|E|}\ln^{-3}(1+\frac{1}{\sqrt{|E|}})\label{E2}.
\end{eqnarray}
\end{lemma}
After this change of variables  Hamiltonian $H$
becomes
$$
H^{\pm}:=H^{\pm}(I,\f,P,Q):=h^{\pm}(I,P)+\mu f^{\pm}(I,\f,P,Q):=
h(I)+E^{\pm}(P)+\mu f^{\pm}(I,\f,P,Q)
$$
where $ f^{\pm}(I,\f,P,Q):=f(I,\f,\phi^{\pm}(P,Q)).$
\\[2mm]
\indent
{\bf Stability in the region ${\cal E}_1^+.$} \noindent
In the region ${\cal E}_1^+,$ the proof of the stability of the actions variables  follows
by a straightforward application of the Nekhoroshev Theorem as proved in
Theorem 1 of \cite{Po}. 
In order to apply such Theorem we need some definitions.
For $ l,m>0$, a function $h:=h(J)$ 
is said to be {\sl l,m-quasi-convex} on $ A \subset{\bf R}^{d+1},$ if 
at every point $J\in A$ at least one of the inequalities
$$
|\langle h'(J),\xi\rangle|>l |\xi |\ ,\ \ \ \langle h''(J)\xi,\xi\rangle\geq m |\xi |^2
$$
holds for each $\xi\in{\bf R}^{d+1}.$ 
Using the previous lemma it is possible to prove that, for every $\ov{r}>0,$
the Hamiltonian $h^{+}$ is {\sl l,m-quasi-convex} in the set 
$ S := D^{+}_{r_0} \times \{ I\in {\bf R}^{d+1}\ | \  
| I| \leq \ov{r} \}_{r_0} $ with $ l, m = O(1).$ In the previous set also 
holds  
$$
\| (h^{+})''\| =: M = 
O(\mu^{-c_d}\ln^{-3}(1/\mu))\ ,\ \ \ \| (h^{+})'\| =: \Omega_0 = O(1).
$$
Putting  $ \e := \mu \| f^{+} \|_{ S, s_0 } = O ( \mu ), $ 
$\e_0:=2^{-10}r_0^2m(m/11M)^{2(d+2)}=O(\mu^{2c_d(d+3)}\ln^{6(d+2)}(1/\mu)),$
$\alpha:=(1- 2c_d(d+3))/2(d+2)$ we obtain that, 
if the initial data 
$(I(0),\f(0),p(0),q(0))\in{\cal E}_1^{+},$ that is $P(0)\in D^{+},$
then 
\be\label{expstabf}
| I(t)-I(0)|\leq const.\mu^{\alpha}\ln^{-3}(1/\mu)\ ,\ \ \ {\rm for}\ \ \ 
|t|\leq const. \exp (const.\mu^{-\alpha}\ln^2(1/\mu )).
\ee
If $ c_d < 1 / 2 ( d + 3 ) $ then $ \a > 0 $ and  
we obtain stability for exponentially long times.
\\[2mm]
\indent
{\bf Stability in the region ${\cal E}_1^-.$} 
In  the region ${\cal E}_1^-$ we cannot use  
the Nekhoroshev Theorem as proved in
\cite{Po}, because $E^-$ is concave and so $h^-$ is not quasi-convex.
However we can still apply the Nekhoroshev Theorem 
in its original and more general form as 
proved in \cite{N} (see also \cite{N79});
in fact the function $ h^- $  proves to be {\sl steep} 
(see Definition 1.7.C. pag. 6  of \cite{N}).

\noindent
For simplicity we prove the {\sl steepness} of the function $h^-$ 
in the case $ d = 1$ only. In this case
$ h^- = h^-(I_0,I_1,P)= I_0 + I_1^2/2 + E^-(P).$
We need more informations on the function $E^-.$
In the following, in order to simplify the notation,
we will forget the apex $^-$ writing, for example, $E=E^-$ and  $P=P^-.$ 

By (1.11) of \cite{N}, since $ \nabla h^-\neq 0 $, 
a sufficient condition for $ h^- $ to be steep is  that the system
\begin{eqnarray}\label{sistema}
\eta_1+I \eta_2+E'(P)\eta_3 &:=& 0\nonumber\\
 \eta_2^2+E''(P) \eta_3^2&:=& 0\nonumber\\
E'''(P)\eta_3^3 &:=& 0
\end{eqnarray}
has no real solution apart from the trivial one $\eta_1=\eta_2=\eta_3=0$.

\noindent
Making the change of variable
$\psi=\arccos(1-\wtilde{E}+\xi\Et),$
where $\Et=E+2,$
we get\footnote{We will denote 
with `` $\dot{}$ '' the derivative with respect to $E,$ and with `` ${'}$ '' the derivative with respect to $P.$ }
\begin{equation}\label{leP}
\dP(E)=\int_0^1\,F_1(\xi;E)\, d\xi,\ \ \ \ddP(E)=3^{-1/2}\int_0^1\,F_2(\xi;E)\, d\xi,\ \ \ {\dddot P}(E)=\int_0^1\,F_3(\xi;E)\, d\xi,
\end{equation}
where
\begin{eqnarray}\label{leF}
F_1(\xi;E) &:=& \frac{\sqrt{2}}{\pi\sqrt{\xi}\sqrt{1-\xi}\sqrt{\Et\xi-E}}\nonumber\\
F_2(\xi;E) &:=& \frac{\sqrt{6}\sqrt{1-\xi}}{2\pi\sqrt{\xi}(\Et\xi-E)^{3/2}}\nonumber\\
F_3(\xi;E) &:=& \frac{3\sqrt{2}(1-\xi)^{3/2}}{4\pi\sqrt{\xi}(\Et\xi-E)^{5/2}}.
\end{eqnarray}
From the equation $E(P(E))=E,$
deriving with respect to $E,$ we obtain that 
$$
E'''(P(E))=-(\dP(E))^{-5}[\dP(E){\dddot P}(E)-3(\ddP(E))^2].
$$
We want to prove that 
\begin{equation}\label{E3}
E'''(P(E))<0. 
\end{equation}
for every $E$ with $ -2<E<0.$
This is equivalent to prove that
$ \dP(E){\dddot P}(E)>3(\ddP(E))^2.$
Using (\ref{leF})
we see that $F_1F_3=F_2^2$ and hence,
noting that $F_3(\xi;E)$ is not proportional to $F_1(\xi;E)$ 
for every $E$ fixed, we conclude that
$\int F_1\int F_3>(\int F_2)^2$
by a straightforward application of Cauchy-Schwarz inequality
and (\ref{E3}) follows from (\ref{leP}).

By (\ref{E3}) the unique solution of 
the system (\ref{sistema}) is the
trivial one $\eta_1=\eta_2=\eta_3=0$,
hence the function $h^-$ is steep.
It is simple to prove that 
the so called {\sl steepness coefficients} and {\sl steepness indices} 
(see again Definition 1.7.C. pag. 6 of \cite{N})
can be taken uniformly for $-2+\mu^{c_d}\leq E \leq -\mu^{c_d}:$
that is they do not depend on $\mu.$

\noindent
Now we are ready to apply the Nekhoroshev Theorem
in the formulation given in Theorem 4.4 of \cite{N}. 
In order to use the notations of \cite{N} we need
the following substitutions\footnote{We observe that we do not 
need to introduce the $(p,q)$ variables so in our case $C=+\infty$.}:
$$
(I,P)\rightarrow I,\ \  (\f,Q)\rightarrow \f,\ \  H^-\rightarrow H, 
\ \ h^-\rightarrow H_0,\ \   \mu f^-\rightarrow H_1,\ \ r_0\rightarrow \rho,
$$
$$ 
\{ I\in {\bf R}^{d+1}\ \ | \ \  | I|\leq \ov{r}\}\times D^-\rightarrow G, 
\qquad
\{ I\in {\bf R}^{d+1} \ | \  | I|\leq \ov{r}\}_{r_0}\times {\bf T}^{d+1}_{s_0}\times D^{+}_{r_0}\times {\bf T}_{s_0}
\rightarrow F.
$$
Defining $m:=\sup_F\|\frac{\partial^2 H_0}{\partial I^2}\|$
and remembering (\ref{E2}) and the definition of $r_0,$ we have 
\begin{equation}\label{mr}
m\leq const. \mu^{-c_d}\ln^{-3}(1/\mu),\ \ \ \ 
\rho = const. \mu^{c_d}.
\end{equation}
In order to apply the Theorem we have only to verify the following condition
\begin{equation}\label{condizioneNek}
M := \sup_F |H_1| < M_0
\end{equation}
where $M_0$ depends only on the  steepness coefficients and  steepness indices 
(which are independent of $\mu$) and on $m$ and $\rho$
(which  depend on $\mu$).
Moreover we use the fact that the dependence of $M_0$ on $m $ and $\rho$ is,
``polynomial'' (although it is quite cumbersome):
that is there exist constant $\wtilde{c}_d,\ov{c}_d>0$ 
such that $M_0(m,\rho)\geq const. m^{-\wtilde{c}_d}\rho^{\ov{c}_d}$ 
(see \S 6.8 of \cite{N79}). 
So condition (\ref{condizioneNek})
becomes, using (\ref{mr}),
$$
\mu\leq const.  \mu^{c_d(\wtilde{c}_d+\ov{c}_d)}\ln^{3\wtilde{c}_d}(1/\mu),
$$
which is verified choosing $c_d<(\wtilde{c}_d+\ov{c}_d)^{-1}.$

\noindent
Now we can apply the Nekhoroshev Theorem as formulated in Theorem 4.4 
of \cite{N}, obtaining that if $ (I(0),\f(0),p(0),q(0)) 
\in {\cal E}_1^- $ then 
\be\label{stabexpi}
|I(t)-I(0)|\leq {\rm d}/2:=M^b/2 =O( \mu^b )  \ \qquad \ \forall\ 
|t|\leq T:=\frac{1}{M}\exp \Big( \frac{1}{M} \Big)^a 
= O \Big( \frac{1}{\mu} \exp \Big(\frac{1}{\mu} \Big )^a \Big)
\ee
where $ a, b > 0 $ are some constants depending only on the 
steepness properties of $ H_0 $.
Finally, choosing  
\be\label{eq:cd}
c_d < \min\{ (2d+6)^{-1},\,(\wtilde{c}_d+\ov{c}_d)^{-1}\},
\ee
we have proved the exponential 
stability in the region ${\cal E}_1.$ \\[2mm]
{\bf Stability in the region ${\cal E}_2^+.$}
In the following we will denote $I^*:=(I_1,\dots ,I_d)$
the projection on the last $d$ coordinates.
We shall prove the following lemma 
\begin{lemma}\label{lemma1}
$\forall \kappa > 0, $ $ \exists \kappa_0, \mu_8 > 0 $ 
such that $ \forall \  0 < \mu \leq \mu_8,$
if $(I(t), \vphi(t), p(t), q(t)) \in {\cal E}_2^+ $ for $0<t\leq \ov{T}$, then 
$$
|I^*(t)-I^*(0)|\leq \frac{\kappa}{2} \ \qquad \ 
\forall \  t  \leq\min\{ \frac{\kappa_0}{\mu}\ln\frac{1}{\mu},\ov{T}\}.
$$
\end{lemma}
It is quite obvious that for initial conditions 
 $(I(0), \vphi(0), p(0), q(0)) \in {\cal E}_2^+$, Theorem \ref{thm:Nec}
follows from lemma \ref{lemma1} and the exponential stability in the region
 ${\cal E}_1$.  

In order to prove lemma \ref{lemma1} let us
define, for some fixed $ 0 < \d < \pi/4 $,
the following two regions in the phase space :
$ U := \{ (I,\f,p,q)| 
\ |q| \leq \d \  {\rm mod}\ 2\pi, \ |E(p,q)|\leq 2\mu^{c_d} \} $
and 
$ V := \{ (I,\f,p,q) | \ |q| > \d $   
mod $ 2 \pi, \ |E(p,q) | \leq 2\mu^{c_d} \}$.
We first note that\footnote{
In the following we will use $c_i$ to denote some positive constant independent on $\m.$
}
\be\label{primasta}
 \ z(t) \in V \ \ \forall t_1 < t < t_2, \
|q(t_1)|, |q(t_2)| = \d \ {\rm mod} \  2 \pi \ 
\Longrightarrow t_2 - t_1 < c_1,  \ 
|I (t_2) - I(t_1) |  \leq c_2(t_2-t_1) \mu.
\ee
Indeed in this case  $ \forall \, t_1 < t < t_2 $, 
$ c_3 \leq |\dot{q}(t)|\leq c_4 $.
This implies that $ t_2 - t_1 \leq c_1$ and then, 
integrating the equation of motion  
$ \dot{I}= - \mu \dpr_{\f} f$ in $(t_1, t_2)$,
we immediately get (\ref{primasta}).
We also claim that
\be\label{broop1}
\forall t_1 < t < t_2, \ z(t) \in U \
{\rm and} \ |q(t_1)|, |q(t_2)| = \d \ {\rm mod} \ 2 \pi
\ \Rightarrow 
t_2 - t_1 \geq c_5 |\ln\m | .
\ee
\noindent
We denote  with $t^i_U$ (resp. $t^i_V$) the $i$-th time for which  
the orbits enters in (resp. goes out from)  $U$, so that 
$t^i_U<t^i_V<t^{i+1}_U<t^{i+1}_V$
for $0\leq i\leq i_0$.
From (\ref{broop1})
it follows that 
$ i_0\leq c_6 \kappa_0 / \m$ and, from (\ref{primasta}),
that the time $T_V$ spent by the orbit in the region $V$  is bounded by
$c_7\kappa_0/\m.$

In order to prove (\ref{broop1})
we use the following normal form result for the pendulum 
Hamiltonian $ E ( p, q ) $
in a neighborhood of its hyperbolic equilibrium point
(see e.g. \cite{CG})
 
\begin{lemma}
There exist $ R, \wtilde{\d}>0,$ an analytic function $ g,$ with $g'(0)=-1$
and an analytic canonical transformation
$\Phi: B\longrightarrow\{ |p| \leq\wtilde{\d}\}\times\{|q| \leq\d\  {\rm mod}\ 2\pi \} $
where $B:=\{ |P|,|Q|\leq R\}$,
such that $ E ( \Phi(P,Q)) = g(PQ)$. 
\end{lemma}

In the coordinates $ ( Q, P ) $ the local stable and unstable manifolds
are resp. $ W^s_{loc} = \{ P = 0 \}$ and $ W^u_{loc} = \{ Q = 0 \}$ and 
Hamiltonian (\ref{hamiltoniana}) writes as 
$$
\widetilde{H}:=\widetilde{H}(I,\f,P,Q):=h(I)+g(PQ)+\mu\wtilde{f}(I,\f,P,Q)
$$
where $\wtilde{f}(I,\f,P,Q):=f(I,\f,\Phi(P,Q)).$

We are now able to prove  (\ref{broop1}). 
Certainly 
there exists an instant $ t_1^* \in [t_1, t_2 ) $ for which
$ (p( t_1^* ),q(t_1^*)) \in\Phi (B)$
but, $\forall  t_1 < t < t_1^*$,
$ (p( t ),q(t)) \notin\Phi (B)$.
It follows that, if we take the representant $q(t_1)\in [-\d,\d]$,
then $p(t_1^*)q(t_1^*)<0.$
We will denote with $ Z(t):=(I(t),\f(t),P(t),Q(t))=(I(t),\f(t),\Phi^{-1}(p(t),q(t)))$
the corresponding solution of the  
Hamiltonian system associated to $\widetilde{H}$. 
From the fact that $|q(t_1^*)|=\d$ or $(p(t_1^*), q(t_1^*)) \in
\partial \Phi (B)$ and that $|g(PQ)| \leq \mu^{c_d}$, $p(t_1^*)q(t_1^*)<0$,
it follows that 
$|P(t_1^*)|\leq c_8\m^{c_d}$ and $|Q(t_1^*)|\geq c_9.$

In the same way 
there exists an instant $t_2^*$ with $t_1< t_1^*<t_2^*< t_2$
for which $(P(t_2^*),Q(t_2^*))\in B$ but, $ \forall  t > t_2^* $ 
$ ( P(t), Q(t)) \notin B $;
in particular it results $ |P(t_2^*)|\geq c_{10}.$
We claim that $ t_2^* - t_1^* \geq c_{11} \ln(1/\mu) $.
Indeed $ P ( t ) $ satisfies the Hamilton's equation
$\dot{P}(t)=-g'(P(t)Q(t))P(t)
-\m\partial_{Q}\widetilde{f}(I(t),\f(t),P(t),Q(t))$
with initial condition $|P(t_1^*)|\leq c_8\m^{c_d}$.
Since   $ |P(t_2^*)|\geq c_{10}$, we can derive
from Gronwall's lemma 
 that  $ t_2^* - t_1^* \geq c_{11} \ln(1/\mu)$, 
which implies (\ref{broop1}).

By the following normal-form lemma 
there exists a close to the identity symplectic change of coordinates
removing the non-resonant angles $ \vphi $ in the perturbation
up to $ O( \mu^2 )$. It can be proved by standard
perturbation theory (see for similar lemmas section $\S 5$ of \cite{CG}).
\begin{lemma}\label{TFN}
Let $\b>0.$
There exist $ R, \rho > 0 $ so small  that, defining 
$ \l := \min_{|\xi|\leq R^2}|g'(\xi)|$,
$ S :=\max_{|\xi|\leq R^2}|g''(\xi)|$, then
$\l \geq 2SR^2$ and  
$\rho$ $\leq$ $\min\{\l/4N,$ $R^2/8s,$ $\b/2N$, $r\}.$ 
Let $ \L $ be a sublattice of ${\bf Z}^{d+1}.$
Let ${\cal D}\subset {\bf R}^{d+1}$ be bounded and 
$\b$-non-resonant mod $\L$, i.e. 
$\forall\, I\in{\cal D}$, $h\in {\bf Z}^{d+1}\setminus \L$,
$|h|\leq N$ it results $|(1,I^*)\cdot h|\geq\b.$
Suppose that 
\begin{equation}\label{condTFN}
\e := \mu\| \wtilde{f}\|_{B,D,s} \leq 2^{-11}\b_*\rho s,
\end{equation}
where\footnote{$B$ and $D$ are thought as complex domains, as in the sequel $\oB$ and $\oD$.}
$D:={\cal D}_{\rho}$, $\b_*:=\min\{\b,\l/2 \}.$
Then there exists an analytic canonical transformation
\begin{equation}
\begin{array}{rcll}
\Psi: & \oD\times {\bf T}^{d+1}_{s/4}\times\oB  & 
\longrightarrow &  D\times {\bf T}^{d+1}_{s}\times B \\[6pt]
 & (\oI,\of,\oP,\oQ) & \longmapsto & (I,\f,P,Q)
\end{array}
\end{equation}
with  $\oB:=\{ |\oP|,|\oQ|\leq R/8 \}$,
$\oD:={\cal D}_{\rho/4},$
such that 
$$
\oH:=\oH(\oI,\of,\oP,\oQ):=\widetilde{H}\circ\Psi=h(\oI)+\og(\oI,\of,\oP\oQ)+\ov{f}(\oI,\of,\oP,\oQ)
$$
with $\og(\oI,\of,\xi):=g(\xi)+f^*(\oI,\of,\xi)$,
$f^*(\oI,\of,\xi)=\sum_{h\in\L,\ |h|\leq N}f^*_h(\oI,\xi)e^{{\rm i}h\cdot\of}$
 and
$\| f^*\|_{\oB,\oD,s/4} \leq\e$.
Moreover the following estimates hold
\begin{equation}\label{IPQf}
|\oI-I| \leq \frac{2^4\e}{\beta_* s},\ \ \
|\oP-P|,|\oQ-Q| \leq \frac{2^5\e}{R\beta_*},\ \ \
\|\ov{f}\|_{\oB,\oD,s/4} \leq \frac{2^9\e^2}{\beta_*\rho s}.
\end{equation}
\end{lemma}

Let ${\cal L}$ be the (finite) set of the maximal 
sublattices $ \L = \langle h_1, \ldots, h_s \rangle
\subset {\bf Z}^{d+1}$ for some independent
$ h_i \in {\bf R}^{d+1}$ with $ |h_i| \leq N$ for 
$i=1, \ldots , s \leq d$.
For $\L\in{\cal L}$ we define 
the $\L$-resonant frequencies $R^{\L}:=\{ I^*\in{\bf R}^d\ | \  
(1,I^*)\cdot h=0,\ \forall\, h\in\L \}$ 
and the set of the $s$-order resonant frequencies
$Z^s:=\cup_{{\rm dim}\L=s}R^{\L}$.

Setting $ h_i = (l_i,n_i) $ with  $l_i\in{\bf R}$, $ n_i\in{\bf
R}^d$, we remark that  
if $ R^{\L} \neq \emptyset $ then $ n_1,\dots,n_s $ are independent.
We also define the $(d-s)$-dimensional linear subspace
(associated with the affine subspace $R^{\L}$) 
$L^{\L} := \cap_{i=1}^s n_i^{\perp} \subset {\bf R}^d$ and we 
denote by $\Pi^{\L}$ the orthogonal 
projection from ${\bf R}^d$ onto $L^{\L}.$

Since $ { \cal L} $ is a finite set, 
$ \alpha := \min_{\Lambda \in {\cal L} } \min_{n\in {\bf Z}^d, | n | \leq N,
\Pi^{\Lambda} n \neq 0} | \Pi^{\Lambda} n | $
is strictly positive.

We now perform a suitable version of the standard
``covering lemma'' in which the whole frequency space is covered by
non-resonant zones.
The fundamental blocks used to construct this covering will be
$r$-neighborhoods of any $R^{\L}$ i.e.
$R^{\L}_r:=\{ I^*\in {\bf R}^d\ | \ {\rm dist}(I^*,R^{\L})\leq r \}$ 
for suitable $r>0$ depending on  ${\rm dim}\L.$
Let $ r_d > 0 $ be such that $(d+1)r_d<c_{12} \kappa,$
for some $c_{12}$ sufficiently small to be determined.
For $ 1 \leq s \leq d-1$ we can define recursively 
numbers $ r_s $ sufficiently small 
such that  $0 < r_s < \a r_{s+1}/2N$, verifying\footnote{
Assumption (\ref{LL}) means that, in order to go from a 
neighborhood of a $(d-s)$-order resonance 
to a different one, we have to pass through an higher order dimensional one.}
\begin{equation}\label{LL}
{\rm dim}\L={\rm dim}\L'=s,\ \ R^{\L}\neq R^{\L'}\quad
\Longrightarrow \quad 
R^{\L}_{(s+1)r_s}\cap R^{\L'}_{(s+1)r_s}\subset \cup_{i=s+1}^d Z^i_{r_i}.
\end{equation} 
We also define,  for $1\leq s\leq d-1,$
$ S^0 := {\bf R}^d\setminus(\cup_{i=1}^d Z^i_{2r_i})$ and
$ S^s := Z^s_{(s+1)r_s}\setminus(\cup_{i=s+1}^d Z^i_{(s+2)r_i})$,
i.e. the $s$-order resonances minus the higher-order ones.
We claim that ${\bf R}^d=S^0\cup\dots\cup S^{d-1}\cup Z^{d}_{(d+1)r_d}$ is the covering that we need.
We also define
$ S^0 \subset S^0_*:={\bf R}^d\setminus(\cup_{i=1}^d Z^i_{r_i})$ and
$ S^s \subset S^s_* := Z^s_{(s+1)r_s}\setminus(\cup_{i=s+1}^d Z^i_{(s+1)r_i})$.

If the orbit lies near  a certain $R^{\L}$ (but far away from higher order resonances)
then the following lemma says that  the drift of the actions $I^*$ in the direction which is parallel to $R^{\L}$ is small. 
\begin{lemma}\label{silvia}
Suppose that $I^*(0)\in S^s$, $I^*(t)\in S^s_*$ and 
$|I^*(t)|\leq \bar{r}+r/2,$  $\forall\, 0\leq t\leq T^*$ for some
$T^*\leq \kappa_0|\ln \m|/\m$ and $0\leq s\leq d-1.$
Then, if $s\geq 1$, there exists a  sublattice $\L\subset{\bf Z}^{d+1}$,
${\rm dim}\L=s$ such that $I^*(t)\in R^{\L}_{(s+1)r_s}\setminus(\cup_{i=s+1}^d Z^i_{(s+1)r_i})$,  $\forall\, 0\leq t\leq T^*$.
Moreover if $\kappa_0$ is sufficiently small\footnote{In the case $s=0$ $\Pi^{\L}$ is simply the identity on ${\bf R}^d.$}
\begin{equation}\label{proiezione}
|\Pi^{\L}(I^*(t)-I^*(0))|\leq r_1/2\qquad \forall\ 0\leq t \leq T^*
\end{equation}
and hence, for $s\geq 1,$ $|I^*(t)-I^*(0)|\leq 2(s+1)r_s+r_1/2.$
In particular for $I^*(0)\in S^0$ we have that 
$|I^*(t)-I^*(0)|\leq r_1/2,$ $\forall\, 0\leq t \leq T^*$.
\end{lemma}
\begin{pf}
In the case $ s = 0 $ we take $\L=\{ 0\}.$
The existence of $\L$  is trivial because  $I^*(0)\in S^s$ and hence  $I^*(0)\in R^{\L}_{(s+1)r_s}$ for some $\L\in{\cal L}$ 
with ${\rm dim}\L=s.$
The fact that $I^*(t)\in R^{\L}_{(s+1)r_s}\setminus(\cup_{i=s+1}^d Z^i_{(s+1)r_i})$,  $\forall\, 0\leq t\leq T^*$,
follows from $I^*(t)\in S^s_*$,  $\forall\, 0\leq t\leq T^*$ and (\ref{LL}).
Now we want to apply lemma \ref{TFN} with $\b:=\a r_1/2$ and ${\cal D}:= R^{\L}_{(s+1)r_s}\setminus(\cup_{i=s+1}^d Z^i_{(s+1)r_i})$.
We have to verify that ${\cal D}$ is $\b$-non-resonant mod $\L.$
Fix $|h_0|\leq N$, $h_0=(l_0,n_0)\notin\L$ (resp. $\neq 0$ for $s=0$). 
We first estimate $|l_0+n_0\cdot I^*_0|$ for all 
$I^*_0\in{\cal D}_0:=R^{\L}\setminus(\cup_{i=s+1}^d Z^i_{(s+1)r_i}).$
If $\L':= \L \oplus \langle h_0 \rangle$
and $ n^*_0:=\Pi^{\L}n_0$ we have two cases: $ n_0^* \neq 0 $ 
or $ n_0^* = 0 $. In $ n_0^* \neq 0 $ 
we can perform the following decomposition: 
$I^*_0=I^*_1+v$ with $I^*_1\in R^{\L'}$, $v\in L^{\L}$ 
and moreover\footnote{We observe that ${\rm dist}(I^*_0, R^{\L'})=|v|.$} 
$v=\pm |v|n^*_0/|n^*_0|.$
Since $I^*_0\notin (\cup_{i=s+1}^d Z^i_{(s+1)r_i})$ then 
 $I^*_0\notin Z^{\L'}_{(s+1)r_{s+1}}$ and, hence 
$|v|\geq (s+1)r_{s+1}.$
Using the previous estimate, the fact that $I^*_1\in\L'$ and  $|n^*_0|\geq\a$,  we conclude that
\begin{equation}\label{amber}
|l_0+n_0\cdot I^*_0|=|(l_0+n_0\cdot I^*_1)+n_0\cdot v|=|n_0\cdot v|=|n^*_0\cdot v|=|v||n^*_0|\geq \a (s+1)r_{s+1}.
\end{equation}
Now we consider the case in which $ n_0^* = 0 $. In this case it is simple
to see that $ h_0 = (l', 0) + h $ where $ h \in \Lambda $  and 
$ l' \in {\bf Z} \setminus \{ 0 \} $. 
So $ | l_0 + n_0 \cdot I_0^* | = |l'| \geq 1 $.
Now we can prove that 
$|l_0+n_0\cdot I^*|\geq\b$ for all 
$I^*\in{\cal D}.$
In fact $I^*=I^*_0+u$ with $I^*_0\in{\cal D}_0$ and $|u|\leq (s+1)r_s.$
Using (\ref{amber}) and $r_s<\a r_{s+1}/2N$, we have
$$
|l_0+n_0\cdot I^*|\geq |l_0+n_0\cdot I^*_0|-|n_0\cdot u|\geq  \a (s+1)r_{s+1}-N(s+1)r_s\geq \a (s+1)r_{s+1}/2\geq\b,
$$
proving that  ${\cal D}$ is $\b$-non-resonant mod $\L.$
Finally we can verify (\ref{condTFN}) if $\m_8$ is sufficiently small. 
Now we are ready to apply lemma \ref{TFN} in order to prove (\ref{proiezione}).
Using (\ref{primasta}),  the fact that $f^*$ contains only the $\L$-resonant Fourier coefficients, 
 (\ref{IPQf}) and  Hamilton's equation for $\ov{H}$
we have
$$
|\Pi^{\L}(I^*(t)-I^*(0))|\leq c_2T_V\m + c_{13}\m^2(\kappa_0|\ln\m|/\m) +  c_{14}i_0\m 
\leq c_2c_7\kappa_0+c_{13}\m\kappa_0|\ln\m|+ c_{14}c_6\kappa_0 \leq r_1/2
$$
if $\kappa_0$ and $\m_8$ are sufficiently small.
\end{pf}

\begin{pfn}{\sc of lemma} \ref{lemma1}.
Suppose first that 
$|I^*(t)|\leq \ov{r}+r/2$ $\forall\, 0\leq t\leq\kappa_0|\ln \m|/\m $.
If $I^*(0)\in Z^{d}_{(d+1)r_d}$ and $I^*(t)\in Z^{d}_{(d+1)r_d}$ 
$\forall 0\leq t\leq\kappa_0|\ln \m|/\m $ then $|I^*(t)-I^*(0)|\leq 2(d+1)r_d$
and the lemma is proved if $c_{12}<1/4$.
Otherwise we can suppose that $I^*(0)\in S^s$ for some $0\leq s\leq d-1.$
If $I^*(t)\in S^s_*$ $\forall 0\leq t\leq\kappa_0|\ln \m|/\m $ then we can apply the lemma \ref{silvia}
proving the lemma for $c_{12}$ small enough.
Suppose that $\exists\, 0<T^*<\kappa_0|\ln \m|/\m$
such that  $I^*(t)\in S^s_*$ $\forall 0\leq t<T^*$ but  $I^*(T^*)\notin S^s_*$. 
We will prove that 
\begin{equation}\label{lanz}
I^*(T^*)\in S^0\cup\dots\cup S^{s-1}
\end{equation}
that means that the orbit can only enter in zones that are ``less'' resonant.
In fact by lemma \ref{silvia} we see that 
$I^*(T^*)\notin\cup_{i=s+1}^d Z^i_{(s+1)r_i}$,
moreover, since  $I^*(T^*)\notin S^s_*,$ we have that
$I^*(T^*)\notin Z^s_{(s+1)r_s}$ and hence
$I^*(T^*)\notin\cup_{i=s}^d Z^i_{(s+1)r_i}$.
If $I^*(T^*)\in S^0$ we have finished.
If $I^*(T^*)\notin S^0$
then $I^*(T^*)\in  \cup_{i=1}^{s-1} Z^i_{2r_i}$ $\subseteq  \cup_{i=1}^{s-1} Z^i_{(i+1)r_i}$.
If $I^*(T^*)\in S^1$ we have finished.
If $I^*(T^*)\notin S^1$
then $I^*(T^*)\notin Z^1_{2r_1}\setminus \cup_{i=2}^{d} Z^i_{3r_i}$ 
and hence  $I^*(T^*)\in \cup_{i=2}^{s-1} Z^i_{(i+1)r_i}.$
Iterating this procedure we prove (\ref{lanz}).

The conclusion is that if the order of resonance  changes along
the orbit, it can decrease only so that the orbit may eventually arrive
in the completely non resonant zone $S^0$ where there is stability.
Considering the ``worst'' case i.e. when   $I^*(0)\in Z^{d}_{(d+1)r_d}$
and the orbit arrives in $S^0$, summing all the contributions from lemma \ref{silvia}, 
we have that, if $c_{12}$ is sufficiently small,
\begin{equation}\label{laura}
|I^*(t)-I^*(0)|\leq  2(d+1)r_d+\sum_{s=1}^{d-1}(2(s+1)r_s+r_1/2)+r_1/2=  \sum_{s=1}^{d}2(s+1)r_s +dr_1/2\leq \kappa/2.
\end{equation}
In order to conclude the proof of the lemma
we have only to prove that if $|I^*(0)|\leq \ov{r}$ then $|I^*(t)|\leq \ov{r}+r/2$ $\forall\, 0\leq t\leq\kappa_0|\ln \m|/\m $.
This is an immediate consequence of (\ref{laura}) and of the fact that $\kappa\leq r.$
\end{pfn}
\\[2mm]
\indent
{\bf Stability in the region ${\cal E}_2^-.$}
If, for all $t \geq 0 \ $
$ (p(t), q(t)) \in {\cal E}_2^- $, then 
it follows easily that $|p(t)|,|q(t)-\pi|=O(\mu^{c_d/2}).$
Then, defining $ f_1(I,\f):=f(I,\f,0,\pi)$ and  
$ f_2 (I, \f, t):= \mu^{-c_d/2}[f(I,\f,p(t),q(t))-f_1(I,\f)]$, 
it results that
$|\partial_I f_2(I,\f;t)|,|\partial_{\f} f_2(I,\f;t)| \leq const$.
Clearly if $ (I(t), \f (t), q(t), p(t) )$ is a solution 
of (\ref{hamiltoniana}) then $ (I(t), \f (t))$
is solution of Hamiltonian
$$
H_1:=H_1(I,\f;t):=h(I)+\mu f_1(I,\f)+\mu^{1+(c_d/2)}f_2(I,\f;t).
$$
Now\footnote{For brevity we prove only the case in which $I(0)$ is in a non-resonant zone.
The resonant case can be treated as in  ${\cal E}_2^+.$ } 
one can construct, in the standard way, 
an analytic symplectic map $ \Phi : (\oI,\of) \to (I,\f) $ with
$| \oI - I | = O( \mu / \beta ),$ 
and two analytic functions $\ov{h},\ov{f}$
such that $[h+\mu f_1]\circ\Phi(\oI,\of)=\ov{h}(\oI)+\ov{f}(\oI,\of)$
with $\|\ov{f}\| = O( \mu^2 ) .$
Defining $f_3:=f_3(\oI,\of;t):=f_2(\Phi(\oI,\of);t)$ we also get that 
$ | \partial_{\oI} f_3(\oI,\of;t)|, 
|\partial_{\of} f_3(\oI,\of;t)|\leq const./\beta.$
The solutions of Hamiltonian $H_1$ 
are symplectically conjugated, via $ \Phi^{-1} $, 
to the solutions of the Hamiltonian
$$
H_2 := H_2(\oI,\of;t) := 
\ov{h}(\oI)+\ov{f}(\oI,\of)+\mu^{1+(c_d/2)}f_3(\oI,\of;t)
$$
for which we obtain, directly from Hamilton's equations, the estimates
$$
|\oI(t)-\oI(0)|\leq const \cdot \mu^{c_d/4}, \qquad 
\ \forall\ |t|\leq const \cdot \mu^{-1-c_d/4}.
$$
It follows that, if $ ( I(0),\f(0),p(0),q(0)) \in  {\cal E}_2^- $, then 
$$
| I(t) - I(0) | \leq 
|I(t)-\oI(t)| + |\oI(t)-\oI(0)| 
+|\oI(0)-I(0)| \leq  const.\mu^{c_d/4}, \qquad 
\forall\ |t|\leq const. \mu^{-1-c_d/4}
$$
(if at some instant $t$ the solution $z(t) $ escapes outside 
$ { \cal E }_2^-$ it is exponentially stable in time).
\\[2mm]
\indent
Finally, from the previous steps, 
we can conclude that there exists $ \mu_1 > 0 $ such that
$ 0 < \mu \leq \mu_1 $ Theorem \ref{thm:Nec} holds.

\section{Appendix}

\begin{pfn}{\sc of lemma } \ref{lem:connecting}.
We shall use the following lemma:
\begin{lemma} \label{lem:hyplin}
There exists $ T_0 > 0 $ such that, $\forall T\geq T_0 $, 
for all continuous $f: [-1,T+1]  \to {\bf R}$, there exists
a unique solution $h$ of
\begin{equation} \label{LPE}
-\ddot{h}+\cos Q_T(t)h =f, \qquad h(0)=h(T)=0.
\end{equation}
The Green operator 
$ {\cal G}: C^0([-1,T+1]) \to C^2([-1,T+1]) $ defined by 
$ {\cal G}(f):= h $, satisfies 
\begin{equation} \label{LPEest} 
\max_{t\in [-1,T+1]} |h(t)|+|\dot{h}(t)| \leq C \max_{t\in
[-1,T+1]}|f(t)|
\end{equation}
for some positive constant $C$ independent of $T$.
\end{lemma}
\begin{pf}
We first note that the homogeneous problem (\ref{LPE})
(i.e. $ f = 0 $) admits only the trivial solution $ h = 0 $.
This immediately implies the uniqueness of the solution of (\ref{LPE}).
The existence result follows by the standard
theory of linear second order differential equations.
We now prove that any solution $ h $ of
$(\ref{LPE})$ satisfies $(\ref{LPEest})$.
It is enough to show that 
$\max_{t\in[-1,T+1]} |h(t)| \leq C'  \max_{t\in[-1,T+1]} |f(t)|$. Indeed
we obtain by $(\ref{LPE})$ that 
$\max_{t\in[-1,T+1]} |h(t)|+|\ddot{h}(t)|\leq (2C'+1) \max_{t\in[-1,T+1]}
|f(t)|$ and, by elementary analysis, this implies $(\ref{LPEest})$
for an appropriate constant $C$.

Arguing by contradiction, we assume that there exist sequences 
$(T_n) \to \infty$, $(f_n)$, $(h_n)$  
such that 
$$
-\ddot{h}_n+\cos Q_{T_n}(t)h_n =f_n, \  \
h_n(0)=h_n(T_n)=0, \  \  |h_n|_n:=\max_{t\in [-1,T_n+1]} |h_n(t)|=1
, \ \  |f_n|_n \to 0.
$$
By the Ascoli-Arzela Theorem
there exists $ h \in C^2 ([-1,\infty), {\bf R} ) $ such that, 
up to a subsequence, $ h_n \to h $ in the topology of $C^2$ uniform 
convergence in $[-1, M ]$ for all $ M > 0 $. Since $ Q_{T_n} \to
q_0-2\pi$ uniformly in all bounded intervals of $[-1,\infty)$,
we obtain that 
\begin{equation} \label{LPElim}
-\ddot{h} + \cos q_0 (t) h = 0, \ \ 
h(0)=0, \   \  \sup_{t\in [-1,\infty)} |h(t)|\leq 1.
\end{equation} 
Now the solutions of the linear differential equation
in $(\ref{LPElim})$ have the form $ h = K_1 \xi + K_2 \psi$, 
where $(K_1,K_2) \in {\bf R}^2$, 
$\xi(t)=\dot{q}_0(t)=\frac{2}{\cosh t}$ and
$\psi(t) =\frac{1}{4}(\sinh t+\frac{t}{\cosh t})$ satisfies 
$\dot{\psi}\xi-\dot{\xi}\psi=1$. The bound on $h$ implies
that $K_2=0$ and $h(0)=0$ implies that $K_1=0$. Hence $h=0$.
In the same way we can prove that $h_n(\cdot -T_n) \to 0$
uniformly in every bounded  subinterval of $(-\infty , 1]$.

Now let us fix $\ov{t}$ such that for all $n$ large enough, for 
all $t\in [\ov{t}, T_n - \ov{t}]$, $\cos Q_{T_n} (t) \geq 1/2$
($\ov{t}$ does exist because of $(\ref{estQT})$).
By the previous step, for $ n $ large enough, there exists
a maximum point $ t_n \in (\ov{t}, T_n -\ov{t}) $ of
$ h^2_n(t)$, i.e.
$ h^2_n(t_n) = |h_n|^2_n = 1 $.  Then
$\dot{(h^2_n)}(t_n)=2h_n(t_n)\dot{h}_n(t_n)=0$
and $\ddot{(h^2_n)} (t_n)=2\ddot{h}_n(t_n) h_n(t_n)+
2\dot{h}^2_n(t) \leq 0$. By the differential equation satisfied
by $h_n$, we can derive from the latter inequality that
$\cos Q_{T_n} (t_n) h_n^2 (t_n) \leq f_n(t_n) h_n(t_n)$, i.e. 
$\cos Q_{T_n} (t_n)  \leq f_n(t_n)$, which, for $n$ large 
enough, contradicts the property of $\ov{t}$ and the fact that
$|f_n|_n \to 0$.
\end{pf}

Now we can deal with the existence result of lemma \ref{lem:connecting}.
Let 
$T:=(\th^--\th^+)$, $\omega=(\f^- -\f^+)/T$, 
$\ov{\f}(t):=\o(t-\th^+)+\f^+$.
In the following we call $c_i$ constants depending only on $f$.
We are searching for solutions $(\vphi, q) $ of (\ref{lagraeqmot})
with $ \vphi( \teta^\pm ) = \vphi^\pm  $,
$ q ( \teta^\pm ) = \mp \pi $, 
in the following form
$$
\left\{\begin{array}{ll}{\varphi(t)=\o(t-\theta^+)+\f^+ +
v(t-\theta^+)}\\ 
{q(t)=Q_T(t-\theta^+)+w(t-\theta^+).}\end{array}\right.
$$
Hence we need to find a solution, in the time interval  $ I:=[-1,T+1]$,
of the following two equations 
\begin{eqnarray} \label{eq:v,w}
\left\{\begin{array}{lll}{ \ddot{v}(t)=
-\m[F_{\f}(v,w)](t)}, & & {v(0)=v(T)=0},\\ 
{[L(w)](t)=[G(v,w)](t):=-[S(w)](t)+\m[F_q(v,w)](t)}, 
& & {w(0)=w(T)=0, }\label{eq}\end{array}\right.
\end{eqnarray}
where 
$$
\begin{array}{llll}
{[F_{\f}(v,w;\l,\m)](t)} 
&:=& \dpr_{\f}f(\o t+\f^++v(t), Q_T(t)+w(t),t+\theta^+),\\
{[F_q(v,w;\l,\m)](t)} &:=& \dpr_{q}f(\o t+\f^+ +v(t),Q_T(t)+w(t),t+\theta^+),\\
{[S(w)](t)} &:=& \sin (Q_T(t)+w(t))-\sin (Q_T(t))-\cos (Q_T(t))w(t),\\
{[L(w)](t)} &:=& -\ddot{w}(t) + \cos Q_T(t) w(t).
\end{array}
$$
We want to solve (\ref{eq:v,w}) as a fixed point problem.
By lemma \ref{lem:hyplin}, the second equation of (\ref{eq:v,w}) 
can be written $ w = K := {\cal G}( -S + \m F_q)$. 
Moreover the first equation (\ref{eq:v,w}) can be written
\begin{equation} \label{eqforv}
v(t)=J(t):= [J(v,w; \l, \m)] (t) :=
\ov{J}(t)-\frac{\ov{J}(0)(T-t)+\ov{J}(T)t}{T}, 
\end{equation}
where, setting $ F_{\f}(s) = F_{\f}(v(s), w(s)) $, 
$$
{[\ov{J}(v,w ;\l,\m)](t)} := -\mu \int_{T/2}^t \int_{T/2}^x F_{\f}(s) \
ds \ dx.
$$

Let us consider the  Banach space $Z=V\times W := {\mathcal
C}^1(I;{\bf R}^d) \times {\mathcal C}^1(I; {\bf R})$, endowed
with the norm $ \| z\| =\|(v,w)\| :=\max\{ \| v \|_V,\, \| w \|_W \}$,
defined by 
\begin{equation} \label{norms}
\| v\|_V:=\sup_{t\in I} \Big[ |v(t)| (1+c_1\m
T^2)^{-1}\b^2+|\dot{v}(t)|\b \Big], \quad 
\quad \| w\|_W:=\sup_{t\in I}\Big[|w(t)|+|\dot{w}(t)|\Big].
\end{equation}
A fixed point of the operator $ \Phi: Z \to Z $ defined 
$\forall z \in Z $ as
$ \Phi ( z ) := \Phi(z;\l,\m) := ( J(z), K(z) )$
is a solution of (\ref{eq:v,w}).
We shall prove in the sequel that $ \Phi $ is a contraction in the 
ball\footnote{If $X$ is a Banach space and $r>0$ 
we define $B_r(X):=\{x\in X; \| x\|\leq r\}.$}  
$D:= B_{\ov{c} \mu(Z)}$ for an 
appropriate choice  of $\ov{c},c_1,C_0$, provided 
$\mu$ is small enough.

We have $|[S(w)](t)|\leq w^2(t),$ so that 
$\forall t,$ $|[G (v,w)](t)|\leq \bc^2\m^2+ c_4\m.$
Now, choosing first $\bc$ sufficiently large 
and then $\m$  sufficiently small, we can conclude
using (\ref{LPEest}) that, if $ z \in D$,
$
\|K (z)\|_W\leq \bc\m /4.
$
Now we study the behaviour of $J.$
Let us first consider $\ov{J}$.
We define 
$$
f_{nl}(t):=f_{nl}(Q_T(t)+w(t)),\quad g_{nl}(t):=\fnl'(Q_T(t)+w(t)),
$$
$$
\alpha_{nl} :=n\cdot\f^++l\theta^+,\quad\bnl:=n\cdot\o+l.
$$
For $t \in [-1,T+1]$, $z\in D$, we want to estimate
$$
\dot{\ov{J}}(t)=-\mu \int_{T/2}^{t} F_{\f}=-\mu  \sum_{|(n,l)|\leq N} 
{\rm i} ne^{{\rm i}\alpha_{nl}}\int_{T/2}^{t}f_{nl}(s)
e^{{\rm i}n\cdot v(s)}e^{{\rm i}\bnl s}ds.
$$
Integrating by parts, we obtain

\begin{eqnarray}
-{\rm i}\bnl\int_{T/2}^{t}f_{nl}(s)e^{{\rm i}n\cdot v(s)}
e^{{\rm i}\bnl s}ds &=& 
f_{nl}(T/2)e^{{\rm i}n\cdot v(T/2)}e^{{\rm i}\bnl T/2}-
f_{nl}(t)e^{{\rm i}n\cdot v(t)}e^{{\rm i}\bnl t} \noindent \\
&+& \int_{T/2}^{t}g_{nl}(s)\dot{Q}_{T}(s)e^{{\rm i}n\cdot v(s)}
e^{{\rm i}\bnl s}ds  \label{secint} \\
&+& \int_{T/2}^{t}(g_{nl}(s)\dot{w} (s)
+f_{nl}(s){\rm i} n\cdot\dot{v} (s))e^{{\rm i}n\cdot v(s)}e^{{\rm i}\bnl s}ds 
\label{thdint}.
\end{eqnarray}

By $(\ref{estQT})$, the term (\ref{secint}) is bounded 
by $c_5 \max \{e^{-K_2t},e^{-K_2(T-t)}\}$. Hence, for $z\in D$, 
\begin{equation}\label{I1}
\int_{T/2}^{t}F_{\f} = u(t)- u(T/2)+ R(t),
\quad {\rm with} \quad |R(t)| \leq  \frac{c_6}{\beta} 
\Big[ \max \Big\{ e^{-K_2t},e^{-K_2(T-t)} \Big\} + \bc
(\m + \frac{\m}{\beta} )T \Big],
\end{equation}
where $u(t)=\sum (n/\beta_{nl}) e^{{\rm i}\alpha_{nl}} 
f_{nl}(t)e^{{\rm i}n\cdot v(t)}e^{{\rm i}\bnl t}$.

So we can write $\ov{J}(t)=j(t)+\mu (t-T/2) u(T/2)$, where 
$$
j(t)=\int_{T/2}^t - \m  u(s) \ ds +  \int_{T/2}^t - \m R(s) \ ds .
$$
By the bound of $R(t)$ given in (\ref{I1}), the second
integral can be bounded by $c_7 (\mu / \beta) [1+\ov{c} T^2 \mu/\beta]$.
Integrating once again by parts as above, 
we find that the first integral is bounded by
$c_8 (\mu / \beta^2)[1 + \bc (\m T/\beta)]$, hence, 
by the condition imposed
on $\mu T$, it can be bounded by $\mu \bc/8\beta^2 $, provided that
$C_0$ has been chosen small enough and $\bc$ is large enough.
Hence
$$
|j(t)| \leq \frac{\m \bc}{\beta^2} \Big[ \frac{c_7}{\bc}+c_7 \mu T^2
+\frac{1}{8} \Big].
$$
In addition
$$
\Big| \frac{d}{dt }{j} (t) \Big| = 
\mu |u(t)+R(t)| \leq c_{10} \frac{\mu \bc}{\beta}
\Big( \frac{1}{\bc}+ \frac{\mu T}{\beta} \Big). 
$$
As a result $||j||_V \leq \mu \bc/4$, provided $\bc$ and $c_1$
have been chosen large enough, $C_0$ small enough.

Now $\ov{J} (t)=j(t)+at+b$, where $a, b \in {\bf R}$, so that we may
replace  $\ov{J}$ with $j$ in $(\ref{eqforv})$.
Since $|J(t)| \leq |j(t)| +\max \{|j(0)|,|j(T)|\}(T+2)/T$
and $|\dot{J}(t)| \leq | dj(t) / dt |+(1/T) \int_{\-1}^{T+1} 
| dj (s)/ dt | \ ds $,  we obtain $||J||_V \leq
3||j||_V  \leq \mu 3\bc /4$.
We have finally proved that $\Phi$ maps $D$  into itself
(in fact into $B_{3\ov{c}\m/4}$). 

Now we must prove that $\Phi$ is a contraction.
$\Phi$ is differentiable and for $z=(v,w) \in D$, 
$(D\Phi (z) [h,g])(t)$ $ = (r(t),s(t))$,  $r$ and $s:
[-1,T+1] \to  {\bf R}$ being defined by
\begin{equation} \label{eqrs}
\ddot{r}(t)=a_1(t).h(t)+b_1(t)g(t), \  r(0)=r(T)=0, \ \
L(s)(t)=a_2(t).h(t)+b_2(t)g(t), \   s(0)=s(T)=0,
\end{equation}
where 
$$
\begin{array}{l}
a_1 (t)=-\mu \partial_{\f \f} f(\o t+\f^++v(t),
Q_T(t)+w(t),t+\theta^+), \   \\ 
b_1(t)=-\mu \partial_{\f  q} f(\o t+\f^++v(t), 
Q_T(t)+w(t),t+\theta^+), \qquad a_2(t)=-b_1(t), \\
b_2(t)=\cos(Q_T(t)+w(t))-
\cos Q_T(t) + \mu \partial_{qq} f (\o t+\f^++v(t), 
Q_T(t)+w(t),t+\theta^+).
\end{array}
$$
By the same arguments as above $ (A, B) \in V_1 \times V $ (where
$V_1:=C^1(I,{\bf R}^{d^2})$) defined by 
$$
\ddot{A}(t)=a_1(t), \   \  A(0)=A(T)=0, \quad \quad 
\ddot{B}(t)=b_1(t), \   \  B(0)=B(T)=0
$$
satisfy $||A||_{V_1}+||B||_V \leq c_{11} \bc \mu$ ($|| \
||_{V_1}$ being defined in the same way as $|| \ ||_V$).

Using an integration by parts, we can derive from 
$(\ref{eqrs})$ and the bound on $||A||_{V_1}+||B||_V$ that
\begin{equation} \label{estimr}
|\dot{r} (t)| \leq c_{12} \bc \frac{\m}{\b} \Big[ \Big(
\frac{1+c_1 \m T^2}{\b^2}||h||_V+ ||g||_W \Big) +
T \Big( \frac{||h||_V}{\b}+||g||_W \Big) \Big].
\end{equation}
Therefore, for $C_0$ small enough, $|\beta \dot{r}(t)|\leq 
1/8 \max\{||h||_V,||g||_W\}$. We derive also from
$(\ref{estimr})$ that
$$
|r(t)| \leq c_{13} \bc \Big[ \frac{\mu T}{\b^3}+\frac{c_1 \m^2
T^3}{\b^3} +\frac{\mu T^2}{\b^2} \Big] \max \{||h||_V,||g||_W \},
$$
which yields $ \beta^2 (1+c_1\mu T^2)^{-1} |r(t)| \leq 
c_{14} \bc ( \mu T /\beta+ (1/c_1)) \max \{||h||_V,||g||_W\} \leq \max
\{ ||h||_V,||g||_W \} / 8 $, provided $C_0$ is small enough
and $c_1/\bc$ is large enough.
Finally $||r||_V \leq \max \{||h||_V,||g||_W \}/4$.

Using the properties of $L$ and the fact that 
$$|a_2(t). h(t)+b_2(t)g(t)| \leq c_{15} \mu (1+ c_1 \mu T^2)/ \b^2 ||h||_V
+ c_{15} (|w(t)|+\mu) ||g||_W$$ we easily derive
$||s||_W \leq \max \{||h||_V,||g||_W \}/4$ (again provided that
$C_0$, more precisely $C_0 c_1$ is small enough). 
We have proved that for a good choice of $\bc ,c_1,C_0$,
$||D\Phi (z) [h,g]|| \leq ||(h,g)||/2$ for $z \in D$. Hence 
$\Phi$ is a contraction.
As a result, it has a unique fixed point $z_{\lambda}$ in $D$
(which in fact belongs to $B_{3\bc \m /4}$).
This proves existence.

Now there remains to prove that $\vphi_{\m, \l}(t), q_{\m,\l}(t)$
are $C^1$ functions of $(\l,t)$.  
Let $(\theta_0^+,\theta_0^-)$ be fixed with $T_0:= \theta^-_0 -
\theta^+_0$  and let
$\Lambda = \{ \lambda \  |  \  |\theta^+ - \theta^+_0| \leq 1/4, \ 
\  |\theta^- - \theta^-_0| \leq 1/4 \}$. For $\lambda \in
\Lambda$ $I_0:=[-1/2, T_0+1/2 ] \subset [-1, \theta^- - \theta^+ +1]$,
hence the restrictions $v_{\lambda}^0$ and $w_{\lambda}^0$
of $v_{\lambda} $ and $w_{\lambda}$ to $ I_0$ are well defined. 

Let $V_0\times W_0 := C^1(I_0 , {\bf R}^n) \times C^1 (I_0,{\bf R})$ 
be endowed with the norm $|| \ ||_0$ as defined in $(\ref{norms})$.
Define $\Psi : \Lambda \to V_0\times W_0$ by $\Psi (\lambda ) =
z^0_{\lambda}$. We shall justify briefly that $\Psi$ is
differentiable and that $ ||D\Psi|| \leq c_{16} \mu$. 
$z^0_{\lambda}$ is the unique solution in
$B_{\bc \m}$ of (\ref{eq}) (with $T=\theta^- - \theta^+$), 
which is equivalent to $(v_{\lambda} , w_{\lambda})=
\Phi (z_{\lambda}; \theta^+, \theta^-, \f^+ , \f^-,\mu)$,
where $\Phi : B_{\bc \m} \times \Lambda \times (0,\m_2) \to 
V_0 \times W_0 $ is smooth. 
Now, by the previous step,  $||D_z \Phi|| \leq 1/2$
everywhere, so that $I-D_z \Phi $ is invertible. Therefore, 
by the Implicit Function Theorem, $\Psi$ is $C^1$. This proves that 
$(\l,t) \mapsto \f_{\m,\l}(t)$ (resp. $(\l,t) \mapsto
q_{\m,\l}(t)$) and $(\l,t) \mapsto \dot{\f}_{\m,\l}(t)$
(resp. $(\l,t) \mapsto \dot{q}_{\m,\l}(t)$) have continuous
partial derivatives w.r.t. $\lambda$ in the set 
$\{ (\l,t) | -1/2+\theta^+ < t < 1/2+\theta^- \}$, and
by the standard theory of differential equations, 
these  partial derivatives have continuous extensions on
$\{ (\l,t) | -1+\theta^+ < t < 1+\theta^- \}$.
Finally, by $(\ref{lagraeqmot})$, $\ddot{\f}_{\m,\l}$ and
$\ddot{q}_{\m,\l}$ depend continuously on $(\l,t)$.
\end{pfn}

\begin{pfn}{\sc of Theorem} \ref{lem:ergo2}.
In order to prove Theorem \ref{lem:ergo2} we need a preliminary
lemma. Observe that  $\Lambda^*_R$ is a finite set which is
symmetric with respect to the origin. Hence,
 if it is not empty there exists $p \in \Lambda^*_R$ such that
$p\cdot \Omega= \alpha (\Lambda, \Omega , R)$. 

\begin{lemma}\label{lem:ergopf}
Assume that $\Lambda^*_R \neq \emptyset$ and let 
$p\in \Lambda^*_R$  be such that    $p\cdot\Omega =\alpha:=
\alpha (\Lambda, \Omega, R)$. Assume moreover that $\alpha >0$ and define
 $E:=[p]^{\bot}$. Then $\Lambda_0:=\Lambda \cap E$ is a lattice of $E$.
In addition : 

\vspace{3mm}

$(i)$ $\displaystyle \frac{\alpha}{\beta |p|} \leq
\dps \frac{2}{R} \ $,  where
$
\beta = \inf \{ |q \cdot \Omega|; 
q\in (\Lambda_0)^*_{\sqrt{3}R/2}\}$, 
$ (\Lambda_0)^* = \{ q \in E \ | \  \forall  x \in \Lambda_0 \
q \cdot x \in {\bf Z}\}$. \\[2mm]
In particular $\alpha \leq 2\beta$.

\vspace{3mm}

$(ii)$ $\alpha(\Lambda, \Omega, \sqrt{7}R/2) \leq \beta$.
\end{lemma}
\begin{pf}
Since $\Lambda$ is a lattice, it is not contained in $E$.
Hence $p \cdot \Lambda$ is a non trivial subgroup of ${\bf Z}$, 
$p \cdot \Lambda = m {\bf Z}$ for some integer $m\geq 1$,
which implies that $p/m \in \Lambda^*$. But $p/m \cdot \Omega = \alpha /
m $ and $|p/m| \leq R$, hence by the definition and the
positivity of $\alpha$, $m=1$. As a result  
 there exists $\ov{x} \in \Lambda$ such 
that $p \cdot \ov{x}=1$. Obviously  $ \Lambda_0 + {\bf Z} \ov{x} 
\subseteq \Lambda$. On the other hand all $x\in \Lambda$ can be written
as $x=(x\cdot p)\ov{x}+y$, where $ y \in \Lambda$, $y \cdot p=0$, {\it i.e.}
$y\in \Lambda_0$. So the reverse inclusion holds and 
we may write $\Lambda = \Lambda_0 + {\bf Z} \ov{x}$.  As a consequence
$\Lambda_0$ is a lattice of $E$ and
$$
\Lambda^{*}=\{ r \in {\bf R}^l \ | \ r \cdot \Lambda_0 \subset 
{\bf Z} \ {\rm and}  \ r \cdot \ov{x} \in {\bf Z} \} =
\{ q+ap \ : \ q\in \Lambda_0^* \  , \
\ a \in {\bf Z}-q \cdot \ov{x} \},
$$
$$
\Lambda^*_R = \{ q+ap \ : \ q\in \Lambda_0^* \  , \
\ a \in {\bf Z} - q \cdot \ov{x} \ , \ 
0< |q|^2+ a^2 |p|^2 \leq R^2\}.
$$
If $\beta = +\infty$ there is nothing more to prove.
If $\beta < +\infty$, let $q\in (\Lambda_0)^*_{\sqrt{3}R/2}$ be such 
that $q \cdot \Omega=\beta$. Let
$$
 S=\{ a\in {\bf R} \ : \ q+ap \in \Lambda^*_R  \} = \{ a\in {\bf R} \ : \
a \in {\bf Z}-q \cdot \ov{x} \ , \ |a| \leq (R^2-|q|^2)^{1/2}/|p| \}.
$$
Since $|q|^2 \leq 3R^2/4$,
$S \supseteq  S':=({\bf Z}-q \cdot \ov{x}) \cap [-R/2|p|,R/2|p|]$.
Hence by the definition of $ \alpha$, 
for all $a \in S'$, $|(q+ap) \cdot \Omega|= |\beta+ a \alpha|\geq \alpha$,
 {\it i.e.} 
$\beta/\alpha \notin (-1-a,1-a)$.
  
As $|p| \leq R$, the interval $[-R/2|p|,R/2|p|]$ has length
$\geq 1$ and must intersect $ ({\bf Z}-q \cdot \ov{x}) $. Therefore 
$S' \not= \emptyset$, more precisely $S'=\{ u,u+1, \ldots , u+K\}$,
 for some integer $K \geq 0$, where
$u=\inf S'$. As a result,
$$
\beta/\alpha \notin \bigcup_{k=0}^K ( -1-u-k,1-u-k ) 
= ( -1-u-K, 1-u ).
$$
Now $S' \cap [-1/2,1/2] \neq \emptyset$, hence $ u + K \geq -1/2$ and 
$-1-u- K < 0$.  As a consequence $ \beta / \alpha \geq 1-u $. Since   
$[-R/2|p|,-R/2|p|+1] \subseteq [-R/2|p|,R/2|p|]$ intersects
${\bf Z} - q \cdot \ov{x}$,  $u \leq -R/2|p|+1$.
Therefore 
$\beta/\alpha \geq R/2|p|$, which is $(i)$. 
In particular, since $|p| \leq R$, $\alpha \leq 2\beta$.

\vspace{2mm}

Finally  there exists  $a\in [-1,0) \cap ({\bf Z} - q \cdot \ov{x})$; $q+ap \in
\Lambda^*$, and $|q+ap|^2=|q|^2+a^2 |p|^2 \leq 3R^2/4 + R^2 =
7R^2/4$. Hence $q+ap \in \Lambda^*_{\sqrt{7}R/2}$.
We have $|(q+ap) \cdot \Omega|=|\beta + a \alpha| \leq \beta$, because
$-1\leq a \leq 0$ and $\alpha \leq 2\beta$. This proves $(ii)$.
\end{pf}

Now we turn to the proof of Theorem \ref{lem:ergo2}
We first prove that the statement is true  for $l=1$,
with $a_1=1/2$. Here $\Lambda=\lambda_0
{\bf Z}$ for some $\lambda_0>0$, and $\Lambda^*=(\lambda_0)^{-1}
{\bf Z}$. We can assume without loss of generality that  $\Omega >0$.
If $ \lambda_0 < 2\delta$, then for all $x\in {\bf R}$,
$d(x, \Lambda) < \delta$. Hence $T(\Lambda, \Omega, \delta)=0$
 
If $ \lambda_0 \geq 2\delta$, then it is  easy to see
that $T(\Lambda, \Omega, \delta)=(\lambda_0-2\delta)/\Omega \leq
\lambda_0/ \Omega$. On the other hand,
$1/\lambda_0 \in \Lambda^*_{1/2\delta}$ and 
$\alpha(\Lambda, \Omega ,1/(2\delta))=\Omega / \lambda_0$. 
The result follows.

\vspace{2mm}

Now we assume that the statement holds true up to
dimension $l-1$ $(l\geq 2)$. We shall prove it 
in dimension $l$. 

Fix $R>0$ and define $\delta_R=(4 a_{l-1}^2/3  +4)^{1/2}/R$. 
We claim that: 
\\[2mm]
{\it $(a)$ If $\Lambda_R^*= \emptyset$ then $T(\Lambda, \Omega ,
\delta_R)=0$.
\\[2mm]
$(b)$   If $\Lambda_R^* \neq \emptyset$, let 
$p \in \Lambda^*_R$ be such that $p \cdot \Omega=\alpha:=
\alpha(\Lambda,\Omega,R)$, and define $\beta$ as in lemma \ref{lem:ergopf}.
Then
$$
T(\Lambda,\Omega,\delta_R) \leq \max \{ \alpha^{-1}, \beta^{-1} \}.
$$}

Postponing the proof of $(a)$ and $(b)$, we show how to define
$a_l$. In the case $(b)$, by lemma \ref{lem:ergopf} (ii),
$T(\Lambda, \Omega, \delta_R) \leq \alpha (\Lambda, \Omega,
\sqrt{7}R/2)^{-1}$.
This estimate obviously holds  in the case $(a)$ too. Hence for
all $R>0$, 
$$
T(\Lambda,\Omega, (4 a_{l-1}^2/3 +4)^{1/2}/R) \leq \alpha (\Lambda, \Omega,
\sqrt{7}R/2)^{-1}. 
$$ 
As a consequence, the statement of Theorem \ref{lem:ergo2}  holds with
$a_l=(\sqrt{7} (4a_{l-1}^2/3+4)^{1/2} /2)$.

\vspace{2mm}

There remains to prove $(a)$ and $(b)$. First assume that 
$ \Lambda^*_R = \emptyset$. Let $p \in \Lambda^* \setminus \{ 0 \}$ 
be such that for all $ p'\in \Lambda^* \setminus \{ 0 \} $, 
$ |p| \leq |p'|$. Then $|p|>R$. 
Let $ E $, $ \Lambda_0$ be defined from $p$ as in lemma
\ref{lem:ergopf}.
 
Arguing by contradiction, we assume that
$(\Lambda_0)^*_{\sqrt{3}R/2} \neq \emptyset$. By the same arguments
as previously
there  exist $q\in (\Lambda_0)^*_{\sqrt{3}R/2}$ and
$a\in [-1/2, 1/2]$ such that $q+ap \in \Lambda^*$.
But $|q+ap|^2=|q|^2+a^2 |p|^2 \leq (3/4)R^2+|p|^2/4 < |p|^2$ and this
 contradicts the definition of $p$. Hence 
$(\Lambda_0)^*_{\sqrt{3}R/2} = \emptyset$  and by the iterative 
hypothesis, all point of $E$ lies at a distance from $\Lambda_0$
less than $2a_{l-1}/\sqrt{3}R$.

From the proof of lemma \ref{lem:ergopf}, there exists $\ov{x} \in
\Lambda$ such that $p \cdot \ov{x}=1$ and $\Lambda=\Lambda_0 + {\bf Z} \ov{x}$.
 Therefore for all $x\in {\bf R}^l$, there
is $x'\in x+\Lambda$ such that $ |x' \cdot p| \leq 1/2$. This implies
that $d(x',E) \leq 1/(2|p|) \leq 1/(2R)$ and hence that
$d(x',\Lambda_0) \leq (4 a_{l-1}^2/3+1/4)^{1/2}/R \leq \delta_R$. 
Hence the distance from any point of ${\bf R}^l$ to $\Lambda$ is
not greater than $\delta_R$. This completes the proof of $(a)$.

\vspace{2mm}

Next assume that $\Lambda_R^* \neq \emptyset$ and
let $p$ be as in lemma \ref{lem:ergopf}. Define $\alpha $ and $\beta$ in the
same way as in lemma \ref{lem:ergopf}. Let $x\in {\bf R}^l$.
Again $\Lambda=\Lambda_0+ {\bf Z} \ov{x}$ for some
$\ov{x} \in \Lambda$ such that $p \cdot \ov{x}=1$, hence there exists
$x'\in x+\Lambda$ such that $p\cdot x' \in [0,1)$. We
have
$$
x'=y+ \frac{w}{|p|^2} p  \  ,  \ \Omega= U+\frac{\alpha}{|p|^2}p,
$$  
with $y,U \in E=[p]^{\bot}$, $w=p \cdot x'\in [0,1)$. We shall assume that 
$\alpha>0$ (if $\alpha=0$, there is nothing to prove). Let
$\ov{t}=w/ \alpha$, and consider the time interval defined by
$$J=[0,1/\beta] \ \ {\rm if} \  \ov{t} < 1/\beta, \quad  \quad
J = [\ov{t}-1/\beta, \ov{t}] \  \ {\rm if} \  \ov{t} \geq 1/\beta.$$
 $J\subset [0,
\max \{1/\beta , 1/\alpha \}]$, and it is enough to prove that 
there exists $t \in J$ such that $ d(x', t \Omega + \Lambda_0) \leq 
\delta_R $. The length of $ J $ is not less than $ 1 / \beta $. 
Hence by the iterative hypothesis, there exists $ t \in J $
such that $d(y, tU + \Lambda_0 ) \leq 2a_{l-1}/(\sqrt{3}R)$
(notice that for all $q\in \Lambda_0^*$, $q \cdot U=q \cdot \Omega$, so that
the linear flow $(tU)$ creates a $2a_{l-1}/(\sqrt{3}R)$-net of
$E/ \Lambda_0$ in time $\beta^{-1}$).
We have
$$
d( x', t \Omega + \Lambda_0 )^2 = 
\Big( \frac{(t-\ov{t})\alpha}{|p|} \Big)^2 + d( y, tU + \Lambda_0)^2
\leq \Big( \frac{\alpha}{\beta |p|} \Big)^2+ \frac{4a_{l-1}^2}{3R^2}.
$$
Hence, by lemma \ref{lem:ergopf} $(i)$ , $d(x', t\Omega + \Lambda_0) \leq 
(4a_{l-1}^2/3+4)^{1/2} /R$. This completes the proof 
of $(b)$.
\end{pfn}

\noindent
{\it Massimiliano Berti and Luca Biasco, S.I.S.S.A., Via Beirut 2-4,
34014, Trieste, Italy, berti@sissa.it, biasco@sissa.it }.
\\[2mm]
{\it Philippe Bolle,
D\'epartement de math\'ematiques, Universit\'e
d'Avignon, 33, rue Louis Pasteur, 84000 Avignon, France,
philippe.bolle@univ-avignon.fr}

\end{document}